\documentclass[11pt,a4paper,reqno]{amsart}

\usepackage{times} 
\usepackage{amsmath} 
\usepackage{amssymb}  
\usepackage{amsthm}
\usepackage{latexsym}
\usepackage{amsfonts,bbm}
\usepackage{xcolor}
\usepackage{mathtools}
\usepackage{enumerate}
\usepackage{cite}
\usepackage{rotating}
\usepackage{pdflscape}
\usepackage{makecell,hhline}
\usepackage{tikz}
\usepackage{graphicx,caption}
\usepackage{todonotes}
\usepackage{appendix}
\usepackage{placeins}

\newcommand{\changed}[1]{\textcolor{black}{#1}}

\newtheorem{thm}{Theorem}[section]
\newtheorem{cor}[thm]{Corollary}
\newtheorem{lem}[thm]{Lemma}
\newtheorem{prop}[thm]{Proposition}
\theoremstyle{definition}
\newtheorem{defn}[thm]{Definition}

\newtheorem{rem}[thm]{Remark}

\newtheorem{ass}[thm]{Assumptions}

\usepackage[foot]{amsaddr}
 
\numberwithin{equation}{section}
\allowdisplaybreaks

\newcommand{\sbvek}[2]{\left[\begin{smallmatrix}#1\\#2\end{smallmatrix}\right]}
\newcommand{\spvek}[2]{\left(\begin{smallmatrix}#1\\#2\end{smallmatrix}\right)}
\newcommand{\sbmat}[4]{\left[\begin{smallmatrix}#1 & #2 \\ #3 & #4\end{smallmatrix}\right]}

\newcommand{\bvek}[2]
{
	\begin{bmatrix}
		#1\\
		#2
	\end{bmatrix}
}

\newcommand{\pvek}[2]
{
	\begin{pmatrix}
		#1\\
		#2
	\end{pmatrix}
}

\newcommand{\re}{\operatorname{Re}}
\newcommand{\IM}{\operatorname{Im}}
\newcommand{\dom}{\operatorname{dom}}
\newcommand{\im}{\operatorname{im}}
\newcommand{\Div}{\operatorname{div}}
\newcommand{\grad}{\operatorname{grad}}
\DeclareMathOperator{\Refl}{{\reflectbox{R}}}

\newcommand{\diag}{\operatorname{diag}}

\newcommand{\tF}{\mathtt F}
\newcommand{\tG}{\mathtt G}
\newcommand{\tK}{\mathtt K}
\newcommand{\tL}{\mathtt L}
\newcommand{\tC}{\mathtt C}
\newcommand{\tD}{\mathtt D}
\newcommand{\tR}{\mathtt R}
\newcommand{\tS}{\mathtt S}

\newcommand{\frakA}{\mathfrak A}
\newcommand{\frakB}{\mathfrak B}
\newcommand{\frakC}{\mathfrak C}
\newcommand{\frakD}{\mathfrak D}
\newcommand{\frakT}{\mathfrak T}
\newcommand{\frakI}{\mathfrak I}
\newcommand{\frakO}{\mathfrak O}

\newcommand{\Id}{I}
\newcommand{\setdef}[2]{\left\{ #1 \left\vert\vphantom{#1} #2 \right.\right\}}

\newcommand{\N}{\mathbb{N}}
\newcommand{\R}{\mathbb{R}}
\newcommand{\C}{\mathbb{C}}
\newcommand{\Uad}{\mathcal{U}_\mathrm{ad}}
\newcommand{\<}{\langle}
\renewcommand{\>}{\rangle}

\title{Linear-quadratic 
optimal control for infinite-dimensional input-state-output systems}

\author[T.\ Reis]{Timo Reis$^1$}

\address[]{$^1$ Institute of Mathematics, Technische Universit\"at Ilmenau, Ilmenau, Germany}

\author[M.\ Schaller]{Manuel Schaller$^2$}
\address[]{$^2$ Faculty of Mathematics, Chemnitz University of Techology, Chemnitz, Germany}

\begin{document}
\maketitle
\begin{abstract}
We examine the minimization of a quadratic cost functional composed of the output and the final state of abstract infinite-dimensional evolution equations in view of existence of solutions and optimality conditions. While the initial value is prescribed, we are minimizing over all inputs within a specified convex subset of square integrable controls with values in a Hilbert space. 
The considered class of infinite-dimensional systems is based on  the system node formulation. Thus, our developed approach includes optimal control of a wide variety of linear partial differential equations with boundary control and observation that are not well-posed in the sense that the output continuously depends on the input and the initial value. 
We provide an application of particular optimal control problems arising in energy-optimal control of port-Hamiltonian systems. Last, we illustrate the our abstract theory by two examples including a non-well-posed heat equation with Dirichlet boundary control and a wave equation on an L-shaped domain with boundary control of the stress in normal direction.
\end{abstract}

\smallskip
\noindent \textbf{Keywords.} linear-quadratic optimal control, infinite-dimensional systems, partial differential equations, boundary control, input constraints

\smallskip
\noindent \textbf{Mathematics subject classications (2020).} 49J20, 49N10, 49K27, 93C25

\bigskip

\section{The objective}\label{sec:objective}
\noindent In this work, we consider dynamic optimal control of abstract evolution equations.
We aim to minimize the cost functional for a specified linear system characterized by the input function $u:[0,T]\to U$, the state $x(t)\in X$ at time $t\in [0,T]$, and the output $y:[0,T]\to Y$, where $U$, $X$, and $Y$ are Hilbert spaces. The cost functional to be minimized is given by
\begin{align}\label{eq:erstegleichung}
\frac{1}{2}\int_0^T\|y(t)-y_\mathrm{ref}(t) \|^2_{Y}{\rm d}t+ \frac{1}{2}\|Fx(T)-z_f\|_Z^2,  
\end{align}
where $y_\mathrm{ref}\in L^2([0,T];Y)$, $z_f\in Z$ ($Z$ is another Hilbert space), and the bounded operator $F:X\to Z$ are given.
The minimization is performed while enforcing the constraint that the input function $u$ belongs to a specified closed and convex set $\mathcal{U}_\mathrm{ad}$. Additionally, the system is initialized with a predefined initial value $x(0)=x_0$.

Minimization is performed subject to a~linear input-state-output system, typically \changed{described}, in a simple representation, \changed{by} $\dot{x}(t)=Ax(t)+Bu(t)$, $y(t)=Cx(t)+Du(t)$.
However, this representation lacks generality,
especially when handling boundary control and
observation scenarios. Instead, we adopt the
system node framework as developed by 
\changed{{\sc Staffans} in \cite{Staffans2005}, based on the considerations by  
	{\sc {\v{S}}muljan} \cite{vsmuljan1986invariant} and {\sc Salamon} \cite{Sala87} on unbounded input and output operators.} 

That is, we consider systems of the form
\begin{equation}
\spvek{\dot{x}(t)}{y(t)}
= \sbvek{A\&B\\[-1mm]}{C\&D} \spvek{{x}(t)}{u(t)},\;\;x(0)=x_0,\label{eq:ODEnode}\end{equation}
where $A\&B:X\times U\supset \dom(A\&B)\to X$,
\changed{$C\&D:\dom(A\&B)\to Y$} are linear operators, where detailed specifications regarding their properties will be outlined in the upcoming section.

\noindent A central aspect is that domain of these operators is not necessarily a Cartesian product of subspaces of $X$ and $U$. This particular setup allows for the incorporation of boundary control, among other considerations. We emphasize that we are not assuming well-posedness, which refers to the existence 
of some constant $c>0$, such that the solutions
of \eqref{eq:ODEnode} (the precise definition of the solution concept will be provided in the following section) fulfill
\begin{equation}
\|y\|_{L^2([0,T];Y)}+\|x(T)\|_X\leq c\,\big(\|u\|_{L^2([0,T];U)}+\|x_0\|_X\big).
\label{eq:wp}
\end{equation}This condition establishes existence and continuous dependence of the state and output upon the initial value and input. While well-posedness is a valuable property in analysis, it poses certain challenges. Firstly, it can be difficult to verify in specific cases. Secondly, and perhaps more crucially, it excludes a range of essential cases: For instance, systems like the heat equation with Dirichlet boundary control and Neumann observation, which are prevalent in real-world applications, fall outside the scope of well-posed systems. Thirdly, though well-posedness significantly simplifies the analysis of the problem, the authors believe that it is not a natural assumption in optimal control. This is because the formulation of the optimization problem itself ensures that the optimal control input corresponds to an output that is square integrable, rather than a distribution. Nevertheless, since our results are also novel for well-posed system, throughout this work we will briefly offer some observations on how our presented theory simplifies for well-posed systems.

It is important to note that our presented theory also encompasses a quadratic penalization of the input. That is, for some $c>0$, the minimization of
\begin{equation}
\frac{1}{2}\int_0^T\|y(t)-y_\mathrm{ref}(t) \|^2_{Y}+c\|u(t)\|^2_{Y}{\rm d}t+ \frac{1}{2}\|Fx(T)-z_f\|_Z^2.\label{eq:costfun2}
\end{equation}
subject to \eqref{eq:ODEnode} can be led back to to an optimal control problem with cost functional of type as in \changed{\eqref{eq:erstegleichung}}
by artificially extending the output of \eqref{eq:ODEnode}. Namely, this can be performed by \changed{setting $\widetilde{y}_\mathrm{ref}:=\spvek{y_\mathrm{ref}}{0}$, and} considering the system
\begin{equation}
\spvek{\dot{x}(t)}{\widetilde{y}(t)}
= \sbvek{A\&B}{\sbvek{C\&D}{\;0\,\sqrt{c}\Id}} \spvek{{x}(t)}{u(t)}.\label{eq:ODEnodeext}\end{equation}

Before delving into the material, we present a brief overview of existing results regarding linear-quadratic optimal control of infinite-dimensional systems. A system node approach to optimal control, as pursued in \cite{OpmeerStaffans2014,OpmeerStaffans2019}, involves minimizing $\|u\|^2_{L^2}+\|y\|^2_{L^2}$ on the positive and negative half-axis. This was achieved by demonstrating that the value function, which maps the initial value to the optimal cost functional value, is quadratic. Based on this finding, a theory was developed that extends algebraic Riccati equations to system nodes. While penalizing the input is also possible in our setup (achieved through an artificial extension of the output, explained in detail in Section~\ref{sec:optcont}), the problem addressed in this work differs for two primary reasons. Firstly, we consider a finite time horizon, and secondly, we allow for input constraints. The latter aspect means that our problem may not necessarily yield a quadratic value function. Consequently, our problem, in general, cannot be addressed using the approaches outlined in \cite{OpmeerStaffans2014,OpmeerStaffans2019}. In this context of Riccati-based approaches, we also mention the earlier works~\cite{Lasiecka1991,Curtain1976} and the textbooks \cite{Lasiecka2000,Lasiecka2000a,Bensoussan2007} for a thorough treatment of infinite-horizon and unconstrained infinite-dimensional optimal control problems with well-posed dynamics. Further, for linear quadratic optimal control of e.g., the particular case of a non-well-posed heat equation with Dirichlet boundary control, we mention~\cite[Chapter 9]{LiYong12} and \cite[Section 9]{Lions71} where the smoothing properties of the heat semigroup are used to deduce existence of solutions and optimality conditions. The main difference to our approach is that we consider a more general setting, as we do not assume, e.g., analyticity of the semigroup.

\noindent Another line of research in optimal control is based on variational theory~\cite{Troeltzsch2010,Hinze08,Lions71}. This approach is complementary to semigroup theory and in particular enables variational discretization techniques such as finite elements. A central tool in these works is the derivation of a control-to-state map in a suitable functional analytic setting (note that in the well-posed semigroup setting, such a mapping can be directly obtained trough the variation of constants formula). In general, such a control-to-state mapping has to be deduced and analyzed in a case-by-case scenario: We refer the reader, e.g., to~\cite{Schiela2013} for a compact approach to linear parabolic equations, to~\cite{Kroener2011} for wave equations with control constraints, to~\cite{Braak2012} for fluid dynamical applications and to \cite{Bommer2016} for the treatment of Maxwell equations. The main difference to our approach is that we present our theory for the general class of system nodes. Therein, we define an unbounded control-to-output map in an abstract setting that allows for existence of solutions or optimality conditions. Consequently, for particular applications, the only remaining task is to verify that the problem under consideration can be formulated as a system node. To this end, however, we refer to \changed{\cite{PhilReis23,ReisSchal23,weiss2013maxwell}}, where system node formulations of dissipative heat, wave, Maxwell's and Oseen equations were presented.

Thus, the main novelty of this work is the abstract and operator-theoretic formulation of existence theory and optimality conditions with very minor assumptions on the structure. In this way, we include various classes of problems with unbounded input or observation, such as, e.g., a heat equation with Dirichlet boundary control and Neumann observation. 

\noindent\textbf{Notation.}
Let $X$ and $Y$ denote Hilbert spaces, consistently assumed to be complex throughout this work. 
The norm in $X$ will be written as $\|\cdot\|_{{X}}$ or simply $\|\cdot\|$, if clear from context. 
The identity mapping in $X$ is denoted as $\Id_{X}$ (or
just $\Id$, if context makes it clear).

The symbol $X^*$ stands for the {\em anti-dual} of $X$, that is, it consists of all continuous and conjugate-linear functionals. Correspondingly, $\<\cdot,\cdot\>_{X^*,X}$ stands for the corresponding duality product.
Further, note that
the Riesz map $R_X$, sending $x\in X$ to the functional $\<x,\cdot\>_{X}$ is a~\changed{unitary operator} from $X$ to $X^*$. If the spaces are clear from context, we may skip the subindices. Further, if not stated else, a~Hilbert space is canonically identified with its anti-dual. Note that, in this case, $R_X=\Id_X$.

The space of bounded linear operators from $X$ to $Y$ is denoted by $L(X,Y)$. As customary, we abbreviate $L(X):= L(X,X)$. The domain $\dom(A)$ of a potentially unbounded linear operator $A:X\supset\dom(A)\to Y$ is usually endowed with the graph norm, represented as $\|x\|_{\dom(A)}:=\big(\|x\|_{X}^2+\|Ax\|_{Y}^2\big)^{1/2}$.

The adjoint of a~densely defined linear operator $A:X\supset\dom(A)\to Y$ is $A^*:Y\supset\dom(A^*)\to X$ with domain
\[\dom(A^*)=\setdef{y\in Y}{\exists\, z\in X\text{ s.t.\ }\forall\,x\in\dom(A):\;\langle y,Ax\rangle_Y=\langle z,x\rangle_X }.\]
The vector $z\in X$ in the \changed{definition of $\dom(A^*)$} is uniquely determined by $y\in\dom(A^*)$, and we define $A^*y=z$.
A~self-adjoint operator $P:X\supset\dom(P)\to X$ is called {\em nonnegative}, if $\langle x,Px\rangle\geq0$ for all $x\in \dom(P)$. The {\em operator square root} of such a nonnegative operator, i.e., a~self-adjoint and
nonnegative (i.e., also self-adjoint) operator whose square is $P$, is denoted by $\sqrt{P}$.



We adopt the notation presented in the book by {\sc Adams} \cite{adams2003sobolev} for Lebesgue and Sobolev spaces. When referring to function spaces with values in a Hilbert space $X$, we indicate the additional notation "$;X$" following the specification of the domain. For instance, the Lebesgue space of $p$-integrable $X$-valued functions over the domain $\Omega$ is denoted as $L^p(\Omega;X)$.

For a finite time horizon $T>0$, the spaces
\begin{align}
H^2_{0l}([0,T];X) &:= \setdef{v\in H^2([0,T];X)}{v(0)= \tfrac{\mathrm{d}}{\mathrm{d}t}v(0) = 0},\label{eq:H2l}\\
H^2_{0r}([0,T];X) &:= \setdef{v\in H^2([0,T];X)}{v(T)= \tfrac{\mathrm{d}}{\mathrm{d}t}v(T) = 0}\label{eq:H2r}
\end{align}
play a~crucial role throughout this work. Note that due to boundedness of the operators defined by evaluation of the function and its derivative at zero, these spaces are closed and thus again are Hilbert spaces when endowed with the usual norm in $H^2([0,T];X)$.

Their dual spaces with respect to the pivot space $L^2([0,T];X)$ are denoted by
\begin{align}
H^{-2}_{0l}([0,T];X) &:= H^2_{0r}([0,T];X)^*,\label{eq:H-2l}\\
H^{-2}_{0r}([0,T];X) &:= H^2_{0l}([0,T];X)^*.\label{eq:H-2r}
\end{align}
The {\em second derivative} $\big(\tfrac{\mathrm{d}^2}{\mathrm{d}t^2}\big)_l:L^2([0,T];X)\to H^{-2}_{0l}([0,T];X)$
is defined via
\begin{align}\label{eq:l2nder}
\begin{split}
\langle \big(\tfrac{\mathrm{d}^2}{\mathrm{d}t^2}\big)_l v, w\rangle_{H^{-2}_{0l}([0,T];X),H^2_{0r}([0,T];X)}:=\langle v, &\tfrac{\mathrm{d}^2}{\mathrm{d}t^2}w\rangle_{L^2([0,T];X)} \quad \forall\,v\in L^2([0,T];X),\,w\in H^2_{0r}([0,T];X).
\end{split}
\end{align}
In an analogous manner, we can also contemplate \changed{another} \emph{second derivative} $\big(\tfrac{\mathrm{d}^2}{\mathrm{d}t^2}\big)_r: L^2([0,T];X) \to H^{-2}_{0r}([0,T];X)$ \changed{where the test functions are chosen to vanish on the left boundary}.\\
Moreover, we introduce the {\em time reflection operator}
\begin{equation}\label{eq:Refl}
\begin{aligned}
\Refl_T:&&L^2([0,T];X)&\to L^2([0,T];X),\\
&&v(\cdot)&\mapsto v(T-\cdot),
\end{aligned}
\end{equation}
which is a~self-adjoint \changed{and unitary operator}. Clearly, $\Refl_T$ also restricts to \changed{a~unitary operator} from the space $H^{2}_{0r}([0,T];X)$ to
$H^{2}_{0l}([0,T];X)$, and also from $H^{2}_{0l}([0,T];X)$ to
$H^{2}_{0r}([0,T];X)$. These are also denoted by $\Refl_T$ for sake of convenience. Moreover, by defining
\begin{multline*}
\langle \Refl_T v, w\rangle_{H^{-2}_{0l}([0,T];X),H^2_{0r}([0,T];X)}\\:=\langle v, \Refl_T w\rangle_{H^{-2}_{0r}([0,T];X),H^2_{0l}([0,T];X)} \ \ \ \forall\,v\in H^{-2}_{0r}([0,T];X),\,w\in H^2_{0r}([0,T];X).
\end{multline*}
we see that $\Refl_T$ extends to a~\changed{unitary operator} from $H^{-2}_{0r}([0,T];X)$ to $H^{-2}_{0l}([0,T];X)$. Its inverse is again an extension of the time reflection operator on $L^2([0,T];X)$, and therefore also denoted by $\Refl_T$.

\section{System nodes and solution operators}\label{sec:sysnode}

For Hilbert spaces $X$, $U$, and $Y$ and linear operators \(A\&B: \dom(A\&B) \subset {X}\times {U} \to {X}\), \changed{$C\&D: \dom(A\&B)\to {Y}$}, we introduce the necessary fundamentals for systems of the form \eqref{eq:ODEnode}.
Here, in the infinite-dimensional case, the operators $A\&B$ and $C\&D$ are not assumed to directly segregate into distinct components that correspond to the state and input and which act on $X$ and $U$, in contrast to the finite-dimensional scenario. This is primarily motivated by the application of boundary control in partial differential equations.
The autonomous dynamics (i.e, those with trivial input $u=0$) are determined by the so-called \emph{main operator}
$A\colon \dom(A)\subset X\to X$
with $\dom(A) \coloneqq \setdef{x\in X}{\spvek x0\in\dom(A\&B)}$ and $Ax \coloneqq A\&B\spvek x0$ for all $x\in\dom(A)$.

\begin{defn}[System node]\label{def:sysnode}
	A {\em system node} on the triple $({X},{U},{Y})$ of Hilbert spaces is a~linear operator $S =\sbvek{A\&B\\[-1mm]}{C\&D}$ with $A\&B:\dom(A\&B)\subset \mathcal{X}\times {U}\to {X}$, \changed{$C\&D:\dom(A\&B)\to {Y}$} satisfying the following conditions:
	\begin{enumerate}[(a)]
		\item\label{def:sysnodea} $A\&B$ is closed.
		\item\label{def:sysnodeb} $C\&D\in L(\dom (A\&B),{Y})$.
		\item\label{def:sysnodec} For all $u\in {U}$, there exists some $x\in {X}$ with $\spvek{x}{u}\in \dom(S)$.
		\item\label{def:sysnoded} The main operator $A$ generates a~strongly continuous semigroup $\mathfrak{A}(\cdot)\colon
		\R_{\ge 0}\to L({X})$ on ${X}$.
	\end{enumerate}
\end{defn}
Next, we define the term {\em solution} for the abstract evolution equation~\eqref{eq:ODEnode}.
\begin{defn}[Classical/generalized trajectories]\label{def:traj}
	Let $T>0$, and let $S = \sbvek{A\& B}{C\& D}$ be a~system node  on $(X,U,Y)$.\\
	A {\em classical trajectory} of \eqref{eq:ODEnode} on $[0,T]$ is a triple
	\[
	(x,u,y)\,\in\,  C^1([0,T];X)\times C([0,T];U)\times C([0,T];Y)
	\]
	which for all $t\in[0,T]$ satisfies \eqref{eq:ODEnode}.\\
	A {\em generalized trajectory} of \eqref{eq:ODEnode} on $[0,T]$ is a~limit of classical trajectories of \eqref{eq:ODEnode} on $[0,T]$ in the topology of $C([0,T];X)\times L^2([0,T];U)\times  L^2([0,T];Y)$.
\end{defn}
Any operator $A\&B$ with the properties \eqref{def:sysnodea}, \eqref{def:sysnodec} and \eqref{def:sysnoded} in Definition~\ref{def:sysnode} can be regarded as a~system node on $(X,U,\{0\})$. Consequently, we may further speak of classical and generalized trajectories $(x,u)$ of
\begin{align}\label{e:ODE}
\dot{x}=A\&B\spvek xu.
\end{align}
We rephrase a~solvability result from \cite{Staffans2005}. Here, besides demanding twice weak differentiability with absolutely integrable derivative of the input, it is required that the pair consisting of the initial state and initial input value lies in $\dom(A\&B)$. For boundary control systems, this means that the boundary value at $t=0$ is consistent with the corresponding boundary value of the prescribed initial state.

\begin{prop}[Existence of classical trajectories {\cite[Thm.~4.3.9]{Staffans2005}}]\label{prop:solex}
	Let $S$ be a system node on $(X,U,Y)$, let $T{>0}$, $x_0\in X$ and $u\in W^{2,1}([0,T];U)$ with $\spvek{x_0}{u(0)}\in \dom(S)$. Then there exists a unique
	classical trajectory $(x,u,y)$ of \eqref{eq:ODEnode} with $x(0)=x_0$.
\end{prop}
Whereas this result could be used to establish a control-to-state map for smooth controls, \changed{we require} the analogous map with square integrable controls, as the cost functional under consideration is merely coercive in the $L^2(0,T;U)$-norm, cf.~\eqref{eq:costfun2}. Thus, going towards a concept of a control to state map in suitable spaces, we will now present more specific results regarding the existence and regularity of the trajectories.

Initially, it is essential to highlight that the operator $A\&B$ can indeed be distinctly separated into components corresponding to the state and the input, aligning with the conventional framework used in numerous studies on infinite-dimensional systems, such as, for instance \cite{TuWe09}. However, to establish such a separation, it becomes necessary to consider extrapolation spaces. 

\begin{rem}[System nodes]\label{rem:nodes}
	Let $S = \sbvek{A\& B}{C\& D}$ be a~system node on $(X,U,Y)$.
	\begin{enumerate}[(a)]
		\item\label{rem:nodesa} For $k\in\N$, the operator $A$ extends to closed and densely defined operator $A_{-k} : X_{-k}\supset\dom(A_{-k}) = X_{-k+1}\to X_{-k}$, where the Hilbert space $X_{-k}$ is the completion of $X$ with respect to the norm $\|x\|_{X_{-k}}:= \|(\alpha \Id-A)^{-k}x\|$ for some $\alpha\in\R$ such that $\alpha \Id-A$ is bijective.
		The semigroup $\frakA(\cdot)$ generated by $A$ extends to a~semigroup $\frakA_{-k}(\cdot)$ on $X_{-k}$. The generator of this semigroup is $A_{-k}$. Moreover, for $X_k:=\dom(A^k)$, $A$ restricts to
		$A_{k} : X_{k}\supset\dom(A_{k}) = X_{k+1}\to X_{k}$, and the semigroup $\frakA_{k}(\cdot)$ generated by $A_k$ is
		the restriction of $\frakA(\cdot)$  to $X_{k}$ \cite[Prop.~2.10.3 \& 2.10.4]{TuWe09}.
		\item\label{rem:nodesb} There exists an operator $B\in L(U,X_{-1})$ such that $[A_{-1}\ B]\in L(X\times U,X_{-1})$ is an extension of $A\& B$. The domain of $A\&B$ (equally: the domain of $S$) satisfies
		\[
		\dom(A\&B)=\setdef{\spvek xu \in X\times U}{A_{-1}x+Bu\in X},
		\]
		see \cite[Def.~4.7.2 \& Lem.~4.7.3]{Staffans2005}.
		\item\label{rem:nodesc} For $k\in\mathbb{Z}$, we denote $X_{\mathrm{d},k}$ as the space constructed as in \eqref{rem:nodesa}, but now from $A^*$.
		Then \cite[Prop.~2.10.2]{TuWe09} yields $X_{\mathrm{d},k}=X_{-k}^*$, where the latter is the dual of $X_{-k}$ with respect to the pivot space~$X$.\\
		In particular, the adjoint of $B$ maps from $X_{\mathrm{d},1}=X_{-1}^*$ to $U$. Then, by using that $A^*$ generates the adjoint semigroup $\frakA^*$ \cite[Prop.~2.8.5]{TuWe09}, we obtain that, for $\frakA_{\mathrm{d},1}$ being the semigroup $\frakA^*$ restricted to $X_{\mathrm{d},1}$ (cf.\ \eqref{rem:nodesa}), it holds that for all $t\geq0$, $B^*\frakA_{\mathrm{d},1}(t)$ is a~bounded operator from $X_{\mathrm{d},1}$ to $U$.
	\end{enumerate}
\end{rem}
In the following we develop a~solution concept that is even more general than generalized trajectories in Definition~\ref{def:traj}. To this end, we first notice that a~generalized trajectory $(x,u)$ of \eqref{e:ODE} fulfills
\begin{equation}\label{eq:mildsol}
\forall\,t\in[0,T]:\quad x(t)=\frakA_{-1}(t)x(0)+\int_0^t \frakA_{-1}(t-\tau)Bu(\tau){\rm d}\tau,
\end{equation}
where the latter has to be interpreted as an integral in the space $X_{-1}$ with $B\in L(U,X_{-1})$ as in Remark \ref{rem:nodes}~\eqref{rem:nodesb}. The state satisfies
\begin{equation}x\in C([0,T];X_{-1}) \cap W^{1,1}([0,T];X_{-2})\label{eq:sol-2}
\end{equation}
as shown in \cite[Thm.~3.8.2]{Staffans2005}.
Consequently, however, the output evaluation $y(t)=C\&D \spvek {x(t)}{u(t)}$ is -- at a glance -- not necessarily well-defined for all $t\in[0,T]$. However, it is shown in \cite[Lem.~4.7.9]{Staffans2005} that the second integral of $\spvek {x}{u}$ is continuous as a~mapping from $[0,T]$ to $\dom(A\&B)$. Thus, the output may -- in the distributional sense -- be defined as the second derivative (as defined in \eqref{eq:l2nder}) of $C\&D$ applied to the second integral of $\spvek {x}{u}$. This can be used to show that $(x,u,y)$ is a~generalized trajectory of \eqref{eq:ODEnode} if, and only if, $(x,u)$ is a~generalized trajectory of \eqref{e:ODE} with $C\&D\int_0^\cdot(\cdot-\tau)\spvek {x(\tau)}{u(\tau)}{\rm d}\tau \in H^2_{0l}([0,T];Y)$ and
\begin{equation}
y=\big(\tfrac{{\rm d}^2}{{\rm d}t^2}\big)_l\,C\&D\int_0^\cdot(\cdot-\tau)\spvek {x(\tau)}{u(\tau)}{\rm d}\tau    \label{eq:IOmap} 
\end{equation}
cf.\ \cite[Eq.\ (4.7.6)]{Staffans2005}. 

Considering now the formula \eqref{eq:IOmap} in the distributional sense allows us to define a solution concept with outputs in the space $H^{-2}_{0l}([0,T];Y)$ defined in \eqref{eq:H-2l}. \changed{The naming is inspired by the concept of very weak solutions to partial differential equations, cf.~\cite{lewis11993very}.}
\begin{defn}[Very generalized trajectory]\label{def:verygen}
	Let $T>0$, and let $S = \sbvek{A\& B}{C\& D}$ be a~system node  on $(X,U,Y)$. Then
	\[
	(x,u,y)\,\in\,  C([0,T];X_{-1})\times L^2([0,T];U)\times H^{-2}_{0l}([0,T];Y)
	\]
	is a {\em very generalized trajectory} of \eqref{eq:ODEnode} on $[0,T]$, if
	\eqref{eq:mildsol} and \eqref{eq:IOmap} hold.
\end{defn}

The findings \changed{\eqref{eq:mildsol}--\eqref{eq:IOmap}} imply that $(x,u,y)$ is a generalized trajectory of \eqref{eq:ODEnode} on $[0,T]$, if, and only if, it is a~very generalized trajectory of \eqref{eq:ODEnode} on $[0,T]$
with $x\in C([0,T];X)$ and $y\in L^2([0,T];Y)$.

Now, we will introduce a series of operators associated with trajectories of \eqref{eq:ODEnode}, which are of essential importance for the addressed optimal control problem. To maintain clarity for both readers and, admittedly, the authors, these operators are systematically presented in a tabular format in Appendix~\ref{sec:ourops}. We recommend the reader to have this table at hand while reading the article. 

First, we define $C\in L(\dom(A),Y)$ by $Cx=C\&D\spvek{x}{0}$.
The introduction of very generalized trajectories gives rise to the introduction of the input-to-state map ${\frakB_T}$, the state-to-output map ${\frakC_T}$, and the input-to-output map ${\frakD_T}$. Namely for all {$x_0\in X$}, $u\in L^2([0,T];U)$, there exist unique $y\in H^{-2}_{0l}([0,T];Y)$ and $x\in C([0,T];X_{-1})$ with $x(0)=x_0$ (defined by \eqref{eq:mildsol} and \eqref{eq:IOmap}) such that $(x,u,y)$ is a~very generalized trajectory of \eqref{eq:ODEnode} on $[0,T]$. Thus, we may define the operators
\begin{align}
{\frakB_T}:&& L^2([0,T];U)&\to X_{-1},\label{eq:genISmap}\\
&&u&\mapsto\int_0^T \frakA_{-1}(T-\tau)Bu(\tau){\rm d}\tau,\hspace*{1.5cm} \text{(input-to-state map)}\nonumber\\[1mm]
{\frakC_T}:&& X_{-1}&\to H^{-2}_{0l}([0,T];Y),\label{eq:genSOmap}\\
&&x_0&\mapsto \big(\tfrac{{\rm d}^2}{{\rm d}t^2}\big)_l C\int_0^\cdot(\cdot-\tau) \changed{\frak{A}_{-1}}(\tau) x_0{\rm d}\tau, \hspace*{.75cm} \text{(state-to-output map)}\nonumber\\[1mm]
{\frakD_T}:&& L^2([0,T];U)&\to H^{-2}_{0l}([0,T];Y),\label{eq:genIOmap}\\
&&u&\mapsto \big(\tfrac{{\rm d}^2}{{\rm d}t^2}\big)_l\,C\&D\int_0^\cdot(\cdot-\tau)\spvek {\frakB_\tau u}{u(\tau)}{\rm d}\tau. \hspace*{.42cm} \text{(input-to-output map)}\nonumber
\end{align}
\changed{These operators are all bounded. For $\frakB_T$, this follows from the fact that $\frakA_{-1}$ is a~strongly continuous semigroup on $X_{-1}$. Further, boundedness of 
	$\frakC_T$ and $\frakD_T$ follows from a~combination of Proposition~\ref{prop:solex} with the closed-graph theorem \cite[Thm. 7.9]{Alt16}.} 
We note that, if $\im \left({\frakB_T}\right)\subset X$, then $B$ is called an~{\em admissible control operator for $\frakA$}. Further, $C$ is
called an~{\em admissible observation operator for $\frakA$}, if $\im \left({\frakC_T}\right)\subset L^2([0,T];Y)$.
Clearly, the latter two properties are fulfilled, if the system is well-posed.

As the optimal control problem introduced in the previous section involves a weighting of the terminal value $x(T)$, there arises a need for an operator which maps a~pair consisting of an initial value $x_0$ and input $u \in L^2([0,T];U)$ to $Fx(T)$, where $F \in L(X,Z)$ and a Hilbert space $Z$ are given. However, defining such an operator is not straightforward, especially if the control operator $B$ is not admissible for $\frakA$. In fact, we impose an additional condition on $F$ to facilitate a meaningful definition for at least a certain class of $F \in L(X,Z)$. 
The foundation for the subsequent definition is as follows, where we use Proposition~\ref{prop:solex} to verify that ${\frakB_T}u\in X$ for all $u\in H^2_{0l}([0,T];U)$.
\begin{prop}\label{prop:ISmap}
	Let $T>0$, let $S = \sbvek{A\& B}{C\& D}$ be a~system node  on $(X,U,Y)$. Let $Z$ be a~Hilbert space, and let $F\in L(X,Z)$, such that the set
	\begin{equation}
	\mathcal{B}_F:=\setdef{w\in Z}{
		B^*\int_0^\cdot(\cdot-\tau) \frakA(\tau)^* F^*w\, {\rm d}\tau\in H^2_{0l}([0,T];U)}\label{eq:dense1}
	\end{equation}
	is dense in $Z$. Then
	\begin{equation}\label{eq:brevB}   \begin{aligned}
	\breve{\frakT}_{F,T}:\quad \setdef{\spvek{x_0}u\in X_{-1}\times L^2([0,T];U)}{\changed{\frakA_{-1}(T)}x_0+{\frakB_T}u\in X}&\to Z\\
	\spvek{x_0}u&\mapsto F\changed{\frakA_{-1}(T)} x_0+F{\frakB_T}u
	\end{aligned}
	\end{equation}
	is closable with respect to $u$, that is, if $x_0\in X_{-1}$, and $(u_{1n}),(u_{2n})$ are sequences in $L^2([0,T];U)$ which are both converging to the same limit, and, moreover
	\begin{align*}
	\forall\,n\in\N,\, i\in\{1,2\}:&\quad \changed{\frakA_{-1}(T)}x_0+{\frakB_T}u_{in}\in X,\\
	\forall\,i\in\{1,2\}:&\quad \big(\breve{\frakT}_{F,T}\spvek{x_0}{u_{in}}\big)\to z_i\in Z,
	\end{align*}
	then $z_1=z_2$. Further, if $x_0\in X_{-1}$, such that there exists some $\hat{u}\in L^2([0,T];U)$, with  $\changed{\frakA_{-1}(T)}x_0+{\frakB_T}\hat{u}\in X$, then the set
	\begin{equation}
	\setdef{u\in L^2([0,T];U)}{\changed{\frakA_{-1}(T)}x_0+{\frakB_T}{u}\in X}\label{eq:initialset}
	\end{equation}
	is dense in $L^2([0,T];U)$.
\end{prop}
\begin{proof}
	The claimed density of \eqref{eq:initialset} follows, since
	$\changed{\frakA_{-1}(T)} x_0+{\frakB_T}(\hat{u}+u)\in X$ for all $u\in H^2_{0l}([0,T];U)$. 
	Thus, \changed{it remains} to show closability.
	
	Recall from Remark~\ref{rem:nodes} that, for all $t\geq0$, $B^*\frakA_{\mathrm{d},1}(t)$ is a~bounded operator from $X_{\mathrm{d},1}$ to $U$, where $\frakA_{\mathrm{d},1}$ is the restriction of $\frakA^*$ to $X_{\mathrm{d},1}$. As a~consequence, we have for all $x\in X$ that the second integral of $t\mapsto \frakA(t)^*x$ is continuous as an $X_{\mathrm{d},1}$-valued function.
	
	Assume that $x_0\in X_{-1}$, and $(u_{in})$, $i=1,2$ are sequences \changed{with properties as defined in the claim}. Then  $(u_n)=(u_{1n}-u_{2n})$ converges to zero, and
	we obtain  for all $w\in\mathcal{B}_F$,
	{\allowdisplaybreaks
		\begin{align*}
		\langle w,z_1-z_2\rangle_Z
		&=\lim_{n\to\infty}\langle w,\breve{\frakT}_{F,T} \spvek{0}{u_n}\rangle_Z\\
		&=\lim_{n\to\infty}\langle w,F\frakB_T u_n\rangle_Z\\
		&=\lim_{n\to\infty}\left\langle F^*w,\int_0^T\changed{\frakA_{-1}}(t)B u_n(T-t){\rm d}t\right\rangle_X\\
		&=\lim_{n\to\infty}\left\langle F^*w,\int_0^T\big(\tfrac{{\rm d}^2}{{\rm d}t^2}\big)_l\int_0^t(t-\tau)\changed{\frakA_{-1}}(\tau)B{\rm d}\tau u_n(T-t){\rm d}t\right\rangle_X\\
		& =\lim_{n\to\infty}\int_0^T\left\langle F^*w,\big(\tfrac{{\rm d}^2}{{\rm d}t^2}\big)_l\int_0^t(t-\tau)\changed{\frakA_{-1}}(\tau)B{\rm d}\tau u_n(T-t)\right\rangle_X {\rm d}t\\
		&=\lim_{n\to\infty}\int_0^T\left\langle \big(\tfrac{{\rm d}^2}{{\rm d}t^2}\big)_l\int_0^t(t-\tau)\frakA(\tau)^*F^*w{\rm d}\tau ,Bu_n(T-t)\right\rangle_{X_{\mathrm{d},1},X_{-1}} {\rm d}t\\
		&=\lim_{n\to\infty}\int_0^T\left\langle \big(\tfrac{{\rm d}^2}{{\rm d}t^2}\big)_lB^*\int_0^t(t-\tau)\frakA(\tau)^*F^*w{\rm d}\tau, u_n(T-t)\right\rangle_{U} {\rm d}t\\
		&=\lim_{n\to\infty}\left\langle
		\big(\tfrac{{\rm d}^2}{{\rm d}t^2}\big)_lB^*\int_0^\cdot(\cdot-\tau)\frakA(\tau)^*F^*w{\rm d}\tau,u_n(T-\cdot)\right\rangle_{L^2([0,T];U)}=0.
		\end{align*}}
	Density of $\mathcal{B}_F$ in $Z$ yields  $z_1=z_2$, and the statement is proven.
\end{proof}
\changed{Proposition~\ref{prop:ISmap}} gives rise to a~mapping arising from a~closure of $\breve{\frakT}_{F,T}$ with respect to $u$. In the following definition, \changed{we call this object the $F$-terminal value map.}
\begin{defn}[$F$-terminal value map]\label{def:FIS}
	Let $T>0$, let $S = \sbvek{A\& B}{C\& D}$ be a~system node  on $(X,U,Y)$. Let $Z$ be a~Hilbert space, and let $F\in L(X,Z)$, such that the set $\mathcal{B}_F$ as in \eqref{eq:dense1} is dense in $Z$.
	Then the {\em $F$-terminal value map}
	\[\frakT_{F,T}:X_{-1}\times L^2([0,T];U)\supset\dom(\frakT_{F,T})\to Z,\]
	is the closure of $\breve{\frakT}_{F,T}$ as in \eqref{eq:brevB} with respect to $u$. That is, $\spvek{x_0}u\in \dom (\frakT_{F,T})$ if, and only if, there exists sequence $(u_n)$ converging to $u$ in $L^2([0,T];U)$, with additionally, $\changed{\frakA_{-1}(T)}x_0+{\frakB_T}{u}_n\in X$ for all $n\in\N$, and $\big(\breve{\frakT}_{F,T}\spvek{x_0}{u_n}\big)$ converges to some $z\in Z$. In this case we set 
	\[{\frakT}_{F,T}\spvek{x_0}{u}:=z.\]
\end{defn}
By Proposition~\ref{prop:ISmap}, for fixed initial value,  the $F$-terminal value map is a densely defined mapping from $L^2([0,T];U)$. A detailed discussion of the central density assumption on $\mathcal{B}_F$ will be provided in Remark~\ref{rem:denseass}.
\begin{rem}[$F$-terminal value map]\
	\begin{enumerate}[(i)]
		\item The $F$-terminal value map is not necessarily closed as a~mapping on $X_{-1}\times L^2([0,T];U)$. Even the mapping $\breve{\frakT}_{F,T}$ is in general not closable under the assumptions made in Proposition~\ref{prop:ISmap}. To achieve this, it can be seen that, for the set $\mathcal{B}_F$ as in \eqref{eq:dense1}, a necessary and sufficient condition for closability of $\breve{\frakT}_{F,T}$ is that,
		\begin{equation}\setdef{w\in \mathcal{B}_F}{\frakA(T)^*F^*w\in \dom(A^*)}\label{eq:BFdenseadd}\end{equation}
		is dense in $Z$. Though this is true, if $\mathcal{B}_F$ is dense in $Z$ and, additionally, by using \cite[Chap.~II, Thm.~4.6]{EngeNage00}, if $\frakA$ is an~analytic semigroup, a~density claim on \eqref{eq:BFdenseadd} would however exclude a variety of interesting cases. We note that, as the first argument of the $F$-terminal value map stands for the initial value, it is fixed in our optimal control problem. Hence, there is actually no need to presume additional properties which guarantee the ``full closedness of $\frakT_{F,T}$''.
		\item For $x_0\in X_{-1}$ such that $\spvek{x_0}{u}$ for some $u\in L^2([0,T];U)$, it is required that there exists some $\hat{u}\in L^2([0,T];U)$ with
		$\frakA_{-1}x_0+{\frakB_T}\hat{u}\in X$. One can think about more general situations where $Fx(T)$ also makes sense even though some \changed{$\hat{u}\in L^2([0,T];U)$ with
			$\frakA_{-1}x_0+{\frakB_T}\hat{u}\in X$ does not exist} (such as, for instance, when we have that $F$ extends to a~bounded operator from $X_{-1}$ to $Z$). One can think about a~more general definition of the operator $\breve{\frakT}_{F,T}$ in Proposition~\ref{prop:ISmap} by presuming that there exists a~Hilbert space $\tilde{X}$ with dense embeddings $X\subset\tilde{X}\subset X_{-1}$, such that, in the definition of $\breve{\frakT}_{F,T}$ in \eqref{eq:brevB}, $X$ is replaced with $\tilde{X}$. To avoid further exaggerating the technical aspects of the whole approach, we refrain from making this generalization here.
		\item 
		If $B$ is an admissible control operator for $\frakA$, then $X\times L^2([0,T];U)$ is contained in the domain of the $F$-terminal \changed{value map}.
	\end{enumerate}
\end{rem}

To finalize the definitions necessary for the optimal control problem, summarized in Appendix~\ref{sec:ourops}, we need to further introduce $F$-input map and $G$-output map.
\begin{defn}[$F$-input map and $G$-output map]\label{def:FIO}
	Let $T>0$, let $S = \sbvek{A\& B}{C\& D}$ be a~system node  on $(X,U,Y)$. Let $Z$ be a~Hilbert space, and let $F\in L(X,Z)$, $G\in L(Z,X)$.
	\begin{enumerate}[(a)]
		\item Let $\mathcal{B}_F$ in \eqref{eq:dense1} be dense in $Z$.
		Then the {\em $F$-input map}
		\[{\frakI}_{F,T}:\;L^2([0,T];U)\supset \dom({\frakI}_{F,T})\to {Z}\times L^2([0,T];Y)\]
		is defined by
		\begin{align*}
		\dom({\frakI}_{F,T})&=\setdef{u\in L^2([0,T];U)}{\spvek{0}u \in \dom({\frakT}_{F,T})\,\wedge\,{\frakD_T}u\in L^2([0,T];Y)},\\
		\quad{\frakI}_{F,T}\changed{u}&=\pvek{{\frakT_{F,T}}\spvek{0}u}{{\frakD_T}u}.
		\end{align*}
		\item
		Further, the {\em $G$-output map}
		\[{\frakO}_{G,T}:\;{Z}\times L^2([0,T];U)\supset \dom({\frakO}_{G,T})\to  L^2([0,T];Y)\]
		is the mapping with
		\[\begin{aligned} \dom({\frakO}_{G,T})&=\setdef{\spvek{z}{u}\in Z\times L^2([0,T];U)}{{\frakC_T}Gz+{\frakD_T}u\in  L^2([0,T];Y)},\\
		{\frakO}_{G,T}\spvek{z}{u}&={\frakC_T}Gz+{\frakD_T}u. 
		\end{aligned}\]
	\end{enumerate}
\end{defn}


Now, we show under suitable assumptions, that the $F$-input map and the $G$-output map are closed and densely defined.

\begin{prop}\label{prop:closedop}
	Let $T>0$, let $S = \sbvek{A\& B}{C\& D}$ be a~system node  on $(X,U,Y)$. Let $Z$ be a~Hilbert space, and let $F\in L(X,Z)$, $G\in L(Z,X)$. Then the following holds:
	\begin{enumerate}[(a)]
		\item\label{prop:closedopa} If the set $\mathcal{B}_F$ as in \eqref{eq:dense1} is dense in $Z$, then the $F$-input map ${\frakI}_{F,T}$ is closed and densely defined.
		\item\label{prop:closedb} The $G$-output map ${\frakO}_{G,T}$ is closed. If, further, the set
		\begin{equation}
		\mathcal{C}_G:=\setdef{z\in Z}{{\frakC_T}Gz\in L^2([0,T];Y)}\label{eq:dense2}
		\end{equation}
		is dense in $Z$, then ${\frakO}_{G,T}$ is densely defined.
	\end{enumerate}
\end{prop}
\begin{proof}\
	\begin{enumerate}[(a)]
		\item Assume that $(u_n)\subset \dom({\frakI}_{F,T})$ converges in $L^2([0,T];U)$ to some $u$, and $({\frakI}_{F,T}u_n)$ converges in $Z\times L^2([0,T];Y)$ to $\spvek{z}{y}$. Then
		$({\frakT}_{F,T}\spvek{0}{u_n})$ converges in $Z$ to $z$, and closedness of the $F$-terminal value map \changed{(as defined in Definition~\ref{def:FIS})} with respect to $u$ yields
		\[{\frakT}_{F,T}\spvek{0}{u}=z.\]
		By further using that \changed{$L^2([0,T];Y)$ is continuously embedded in $H^{-2}_{0l}([0,T];Y)$}, we have that $({\frakD_T} u_n)$ converges to $y$ in $H^{-2}_{0l}([0,T];Y)$, and boundedness of
		${\frakD_T}$ leads to $y={\frakD_T}u$. Altogether, we have
		\[{\frakI}_{F,T}u=\spvek{z}{y},\]
		which shows that ${\frakI}_{F,T}$ is closed. Dense definition of ${\frakI}_{F,T}$ holds, since by Proposition~\ref{prop:solex},
		\[H^2_{0l}([0,T];U)\subset\dom({\frakI}_{F,T}).\]
		\item Assume that $\big(\spvek{z_n}{u_n}\big)$ is a~sequence in $\dom({\frakO}_{G,T})$ that converges in $Z\times L^2([0,T];U)$ to $\spvek{z}{u}\in Z\times L^2([0,T];U)$, and, further, $\big({\frakO}_{G,T}\spvek{z_n}{u_n}\big)$
		converges in $L^2([0,T];Y)$ to $y\in L^2([0,T];Y)$. Now boundedness of ${\frakC_T}$, $G$ and ${\frakD_T}$ together with the continuous embedding
		of $L^2([0,T];Y)$ in $H^{-2}_{0l}([0,T];Y)$ yields that
		\[{\frakC_T}Gz+{\frakD_T}u=y,\]
		and thus
		$\spvek{z}{u}\in \dom({\frakO}_{G,T})$ with $\dom({\frakO}_{G,T})\spvek{z}{u}=z$. This shows that ${\frakO}_{G,T}$ is a~closed operator.\\
		To show the claim on dense definition, we first observe that, by Proposition~\ref{prop:solex}, we have for all $u\in H^2_{0l}([0,T];U)$ that $\spvek0{u}\in \dom({\frakO}_{G,T})$. Then density of $\mathcal{C}_G$ immediately yields that ${\frakO}_{G,T}$ is densely defined.
	\end{enumerate}
\end{proof}

\begin{rem}\label{rem:denseass}\
	\begin{enumerate}[(i)]
		\item\label{rem:denseassi}     We briefly recall that the purpose of this article is the treatment of the optimal control problem as lined out in Section~\ref{sec:objective}. A~crucial ingredient in our approach to this problem will by the $F^*$-output map, where $F\in L(X,Z)$ is the operator that occurs in the weighting of the terminal value. This will be elaborated in the forthcoming section.
		\item\label{rem:denseassii} 
		The set $\mathcal{B}_F$ in \eqref{eq:dense1} is dense in $Z$, if one of the following conditions hold:
		\begin{enumerate}[(a)]
			\item\label{rem:denseassiia} $B^*$ is an admissible observation operator for $\frakA^*$. Equivalently, by \cite[Thm.~4.4.3]{TuWe09},
			$B$ is an admissible control operator for $\frakA$. In this case, the set specified in \eqref{eq:dense1} even coincides with~$Z$.
			\item\label{rem:denseassiib} The set
			\[\setdef{z\in Z}{F^*\changed{z} \in \dom(A^*)}\]
			is dense in $X$. This is for instance fulfilled in the case where $F$ is the zero operator, or $F$ is injective and has closed range.
		\end{enumerate}
		\item Likewise,
		the set $\mathcal{C}_G$ in \eqref{eq:dense2}
		is dense in $Z$, if one of the following conditions hold:
		\begin{enumerate}[(a)]
			\item $C$ is an admissible observation operator for $\frakA$. In this case, the set specified in \eqref{eq:dense2} even coincides with $Z$.
			\item The set
			\[\setdef{z\in Z}{G\changed{z} \in \dom(A)}\]
			is dense in $X$. This is for instance fulfilled in the case where $G$ is the zero operator or it is surjective.
		\end{enumerate}
		\item If the system is well-posed, we have that ${\frakB_T}\in L(L^2([0,T];U),X)$, ${\frakC_T}\in L(X,L^2([0,T];Y))$ and ${\frakD_T}\in L(L^2([0,T];U),L^2([0,T];Y))$. Consequently, for $F\in L(X,Z)$, $G\in L(Z,X)$, we have, in this case, that the $F$-input map and the $G$-output map \changed{(see Definition~\ref{def:FIO})} are as well bounded, and they moreover simplify to
		\[\frakI_{F,T}=\bvek{F{\frakB_T}}{{\frakD_T}},\quad 
		\frakO_{G,T}=\begin{bmatrix}{{\frakC_T}}G&{{\frakD_T}}\end{bmatrix}. 
		\]
		Thus, we observe that the incorporation of non-well-posed systems significantly complicates matters and demands a certain degree of technical finickiness.\end{enumerate}
\end{rem}
\noindent In the following let $S = \sbvek{A\& B}{C\& D}$ be a~system node  on $(X,U,Y)$. 
We now present a~result on the adjoint of the $F$-input map and the $G$-output map \changed{(as defined in Definition~\ref{def:FIO})} showing that this can be constructed from
the {\em adjoint system node} $S^*=\sbvek{A\&B}{C\&D}^*$. It indeed holds that the latter is a~system node, if $S$ itself is a~system node \cite[Lem.~6.2.14]{Staffans2005}. We denote
\begin{equation}\label{eq:Sstar}
S^*=\sbvek{{[A\&B]}^\mathrm{d}\\[-1mm]}{{[C\&D]}^\mathrm{d}},
\end{equation}
and we consider the {\em adjoint system}
\begin{equation}
\spvek{\dot{x}_\mathrm{d}(t)}{u_\mathrm{d}(t)}
= \sbvek{{[A\&B]}^\mathrm{d}\\[-1mm]}{{[C\&D]}^\mathrm{d}}\spvek{{x}_\mathrm{d}(t)}{y_\mathrm{d}(t)},\;\;x_\mathrm{d}(0)=x_{\mathrm{d}0}.\label{eq:ODEnodeadj}\end{equation}
It is moreover shown in \cite[Lem.~6.2.14]{Staffans2005} that the main operator of $S^*$ is given by $A^*$ \changed{which generates the semigroup $\frak{A}^*$ on $X$, see also Remark~\ref{rem:nodes}(c).} The adjoint system node will be employed to show that the adjoint $F$-input map can be constructed from the $F^*$-output map associated with the system defined by $S^*$ and additional ``time flips'', which are expressed by the time reflection operator as introduced in \eqref{eq:Refl} and the subsequent lines.

\changed{Further, for the adjoint system we have the following operators, analogously to \eqref{eq:genISmap}--\eqref{eq:genIOmap}: For $T>0$, define
	\begin{align}
	{\frakB_{\mathrm{d},T}}:&& L^2([0,T];Y)&\to X_{\mathrm{d},-1},\label{eq:genISmapdual}\\
	&&y_\mathrm{d}&\mapsto\int_0^T \frakA_{\mathrm{d},-1}(T-\tau)B_\mathrm{d}y_\mathrm{d}(\tau){\rm d}\tau,\nonumber\\[1mm]
	{\frakC_{\mathrm{d},T}}:&& X_{\mathrm{d},-1}&\to H^{-2}_{0l}([0,T];U),\label{eq:genSOmapdual}\\
	&&x_{\mathrm{d},0}&\mapsto \big(\tfrac{{\rm d}^2}{{\rm d}t^2}\big)_l C_\mathrm{d}\int_0^\cdot(\cdot-\tau) \frakA_{\mathrm{d},-1}(\tau) x_{\mathrm{d},0}{\rm d}\tau,\nonumber\\[1mm]
	{\frakD_{\mathrm{d},T}}:&& L^2([0,T];Y)&\to H^{-2}_{0l}([0,T];U),\label{eq:genIOmapdual}\\
	&&y_\mathrm{d}&\mapsto \big(\tfrac{{\rm d}^2}{{\rm d}t^2}\big)_l\,[C\&D]^d\int_0^\cdot(\cdot-\tau)\spvek {\frakB_{\mathrm{d},\tau} u_\mathrm{d}}{u_{\mathrm{d}}(\tau)}{\rm d}\tau.\nonumber
	\end{align}
	\changed{where \changed{$\frakA_{\mathrm{d},-1}$} is the extension of $\mathfrak{A}^*$ onto $X_{\mathrm{d},-1}$}, the output operator $C_{\mathrm{d}}\in L(\dom(A^*),U)$ is defined via $C_{\mathrm{d}}x_{\mathrm{d}}=[C\&D]^d\spvek{x_{\mathrm{d}}}{0}$ and $B_\mathrm{d}\in L(Y,X_{\mathrm{d},-1}$ is constructed by Remark~\ref{rem:nodes}(b) applied to the dual system node.}

Before the main result for the adjoints of $\frakI_{F,T}$ and $\frakO_{G,T}$ is presented, we advance an auxiliary result (which generalizes \cite[Lem.~6.2.16]{Staffans2005}) on an integration-by-parts like identity between the (very) generalized \changed{trajectories} of a~system node $S$ and its adjoint.
To this end, we recall from Remark~\ref{rem:nodes}\,\eqref{rem:nodesc} that $X_{\mathrm{d},1}=X_{-1}^*$ {and that $X_{\mathrm{d},-1} = X_{1}^*$.}

\begin{lem}\label{lem:adj}
	Let $T>0$, let $S = \sbvek{A\& B}{C\& D}$ be a~system node on $(X,U,Y)$, and let $S^*$ as in \eqref{eq:Sstar} be the adjoint system node.
	Assume that $(x,u,y)$ is a~classical trajectory of \eqref{eq:ODEnode} with \changed{$u\in H^2_{0}([0,T];U)$} and $x(0)=0$,
	and $(x_\mathrm{d},y_\mathrm{d},u_\mathrm{d})$ is a~very generalized trajectory of \eqref{eq:ODEnodeadj}. Then \changed{$x(T)\in X_1$ with}
	\begin{equation}
	\langle x(T),x_\mathrm{d}(0)\rangle_{X_1,X_{\mathrm{d},-1}}+\langle y,\Refl_T y_\mathrm{d}\rangle_{L^2([0,T];Y)}
	=
	\langle u,\Refl_T u_\mathrm{d}\rangle_{H^2_{0l}([0,T];Y),H^{-2}_{0r}([0,T];Y)}.
	\label{eq:vgenadj}    \end{equation}
\end{lem}
\begin{proof}
	As in the proof of \cite[Lem.~6.2.16]{Staffans2005}, it can be shown that, for
	any~classical trajectory $(x,y,u)$ of \eqref{eq:ODEnode},
	and any classical trajectory $(x_\mathrm{d},y_\mathrm{d},u_\mathrm{d})$ of \eqref{eq:ODEnodeadj}, it holds that
	\[\forall \,t\in[0,T]:\,
	\tfrac{{\rm d}}{{\rm d}t}\langle x(t),x_\mathrm{d}(T-t)\rangle_X+\langle y(t),y_\mathrm{d}(T-t)\rangle_Y=\langle u(t),u_\mathrm{d}(T-t)\rangle_U.
	\]
	Now an integration over $[0,T]$ yields that classical trajectories fulfill
	\begin{equation}
	\langle x(T),x_\mathrm{d}(0)\rangle_X+\langle y,\Refl_T y_\mathrm{d}\rangle_{L^2([0,T];Y)}
	=
	\langle x(0),x_\mathrm{d}(T)\rangle_X+\langle u,\Refl_T u_\mathrm{d}\rangle_{L^2([0,T];Y)}.\label{eq:genadj}
	\end{equation}
	Now assume that $(x,u,y)$ is a~classical trajectory of \eqref{eq:ODEnode} with \changed{$u\in H^2_{0}([0,T];U)$} and $x(0)=0$, and let $(x_\mathrm{d},y_\mathrm{d},u_\mathrm{d})$ be a~\changed{very generalized trajectory} for \eqref{eq:ODEnodeadj}. Then, by density of $\dom(S^*)$ in $X\times Y$ and density of $H^2([0,T];Y)$ in $L^2([0,T];Y)$, there exist sequences $(y_{\mathrm{d},n})$ in $H^2([0,T];Y)$, and $(x_{\mathrm{d},0,n})$ in $X_{\mathrm{d},1}$, such that
	\begin{itemize}
		\item[$\boldsymbol{\cdot}$] $(y_{\mathrm{d},n})$ converges in $L^2([0,T];Y)$ to $y_\mathrm{d}$,
		\item[$\boldsymbol{\cdot}$] $(x_{\mathrm{d},0,n})$ converges in $X_{\mathrm{d},-1}$ to $x_\mathrm{d}(0)$, and
		\item[$\boldsymbol{\cdot}$] $\spvek{x_{\mathrm{d},0,n}}{\changed{y_{\mathrm{d},n}(0)}}\in\dom(S^*)$ for all $n\in \N$.
	\end{itemize}
	Let ${\frakB_{\mathrm{d},T}}$, ${\frakC_{\mathrm{d},T}}$, ${\frakD_{\mathrm{d},T}}$ be the \changed{mappings defined in \eqref{eq:genISmapdual}-\eqref{eq:genIOmapdual} associated to the system node $S^*$}.
	Then, by a~combination of Proposition~\ref{prop:solex} with \eqref{eq:mildsol} and \eqref{eq:IOmap}, we obtain that, for $n\in\N$,
	$x_{\mathrm{d},n}:[0,T]\to X$ and $u_{\mathrm{d},n}\in L^2([0,T];U)$ with
	\begin{align*}
	x_{\mathrm{d},n}(t)&=\changed{\frakA_{\mathrm{d},-1}}(t)x_{\mathrm{d},0,n} + {\frakB_{\mathrm{d},t}}\changed{y}_{\mathrm{d},n},\quad t\in[0,T],\\
	u_{\mathrm{d},n}&={\frakC_{\mathrm{d},t}}x_{\mathrm{d},0,n}+{\frakD_{\mathrm{d},t}}y_{\mathrm{d},n},
	\end{align*}
	$(x_{\mathrm{d},n},y_{\mathrm{d},n},u_{\mathrm{d},n})$ is a~classical trajectory of \eqref{eq:ODEnodeadj}. Further, by boundedness of ${\frakB_{\mathrm{d},T}}$, ${\frakC_{\mathrm{d},T}}$, ${\frakC_{\mathrm{d},T}}$ and the fact that $\frakA^*$ extends to a~strongly continuous semigroup on $X_{\mathrm{d},-1}$ (see Remark~\ref{rem:nodes}\,\eqref{rem:nodesa}), we obtain that $x_{\mathrm{d},n}(T)$ converges in $X_{\mathrm{d},-1}$ to $x_\mathrm{d}(T)$, and  $(u_{\mathrm{d},n})$ converges in $H^{-2}_{0l}([0,T];Y)$ to $u_\mathrm{d}$. \changed{As $(x,u,y)$ is a~classical trajectory of \eqref{eq:ODEnode} with $u\in H^2_{0}([0,T];U)$, we have
		\[\spvek{x(T)}{u(T)}\in\dom(A\&B).\]
		As $u\in H^{2}_{0r}([0,T];Y)$ and $u(T)=0$, it holds that $x(T)\in \dom(A)=X_1$.}
	Now invoking that \eqref{eq:genadj} holds for classical trajectories, we obtain
	\begin{align*}
	\langle x(T),x_{\mathrm{d},n}(0)\rangle_{X_1,X_{\mathrm{d},-1}}+\langle y,\Refl_T y_{\mathrm{d},n}\rangle_{L^2([0,T];Y)}
	&=
	\langle x(T),x_{\mathrm{d},n}(0)\rangle_{X}+\langle y,\Refl_T y_{\mathrm{d},n}\rangle_{L^2([0,T];Y)}\\
	&=
	\langle u,\Refl_T u_{\mathrm{d},n}\rangle_{L^2([0,T];Y)}\\
	&=\langle u,\Refl_T u_{\mathrm{d},n}\rangle_{H^2_{0l}([0,T];Y),H^{-2}_{0r}([0,T];Y)}.
	\end{align*}
	Now the result follows by taking the limit $n\to\infty$.
\end{proof}

\begin{prop}[The adjoint of the $F$-input map]\label{prop:adjFinput}
	Let $T>0$, and let $S = \sbvek{A\& B}{C\& D}$ be a~system node on $(X,U,Y)$.
	Let $Z$ be a~Hilbert space, and let $F\in L(X,Z)$, 
	such that $\mathcal{B}_F$ as specified in \eqref{eq:dense1} is dense in $Z$.
	Further, let $\frakO_{\mathrm{d},F^*,T}$ be the $F^*$-output map corresponding to $S^*$. Then
	\begin{align}
	\frakI_{F,T}^*&
	= \Refl_T\frakO_{\mathrm{d},F^*,T}\sbmat{I}{0}{0}{\Refl_T}.\label{eq:Iadj}
	\end{align}
\end{prop}
\begin{proof}
	We show the equivalent statement
	\[    \left(\sbmat{I}{0}{0}{\Refl_T}\frakI_{F,T}\right)^*
	= \Refl_T\frakO_{\mathrm{d},F^*,T}.
	\]
	We first observe that, by Lemma~\ref{lem:adj}, for all $z\in Z$, $y_\mathrm{d}\in L^2([0,T];Y)$
	and \changed{$u\in H^2_{0}([0,T];U)$}, we have
	\begin{equation}
	\begin{aligned}
	\left\langle \sbmat{I}{0}{0}{\Refl_T}\frakI_{F,T}u,\spvek{z}{ y_\mathrm{d}}\right\rangle_{Z\times L^2([0,T];Y)}\!\!\!
	&=\left\langle \frakI_{F,T}u,\spvek{z}{\Refl_T y_\mathrm{d}}\right\rangle_{Z\times L^2([0,T];Y)}\\
	&=\left\langle \spvek{F{\frakB_{T}}u}{{\frakD_{T}}u},\spvek{z}{\Refl_T y_\mathrm{d}}\right\rangle_{Z\times L^2([0,T];Y)}\\
	&=\langle {F{\frakB_{T}}u},{z}\rangle_Z
	+\langle
	{{\frakD_{T}}u},\Refl_T y_\mathrm{d}
	\rangle_{L^2([0,T];Y)}
	\\
	&=\langle {{\frakB_{T}}u},{F^*z}\rangle_X
	+\langle
	{{\frakD_{T}}u},\Refl_T y_\mathrm{d}
	\rangle_{L^2([0,T];Y)}\label{eq:adjcomp}
	\\
	&=\langle {{\frakB_{T}}u},{F^*z}\rangle_{X_1,X_{\mathrm{d},-1}}
	+\langle
	{{\frakD_{T}}u},\Refl_T y_\mathrm{d}
	\rangle_{L^2([0,T];Y)}\\
	&=\left\langle u, \Refl_T\big({\frakC_{\mathrm{d},T}}F^*z+{\frakD_{\mathrm{d},T}}y_\mathrm{d}\big)\right\rangle_{H^{2}_{0l}([0,T];U),H^{-2}_{0r}([0,T];U)},
	\end{aligned}
	\end{equation}
	\changed{where in the last equality we used \eqref{eq:vgenadj} with $x(T) = \frak{B}_T u$, $y = \frak{D}_Tu$ and $x_{\mathrm{d}}(0)=F^*z$}.
	If $\spvek{z}{y_\mathrm{d}}\in\dom(\Refl_T\frakO_{\mathrm{d},F^*,T})=\dom(\frakO_{\mathrm{d},F^*,T})$, then
	\[\Refl_T\big({\frakC_{\mathrm{d},T}}F^*z+{\frakD_{\mathrm{d},T}}y_\mathrm{d}\big)\in L^2([0,T];U),\]
	and the last expression in the chain of equalities \eqref{eq:adjcomp} can be considered as inner product in $L^2([0,T];U)$.
	In this case, we further have
	\[\Refl_T\big({\frakC_{\mathrm{d},T}}F^*z+{\frakD_{\mathrm{d},T}}y_\mathrm{d}\big)=\Refl_T\frakO_{\mathrm{d},F^*,T}\spvek{z}{y_\mathrm{d}}.\]
	Since \eqref{eq:adjcomp} holds for all $u$
	in the dense subspace \changed{$H^2_{0}([0,T];U)$} of $L^2([0,T];U)$, we have shown
	that $\Refl_T\frakO_{\mathrm{d},F^*,T}$ is a~restriction of $\left(\sbmat{I}{0}{0}{\Refl_T}\frakI_{F,T}\right)^*$. Hence, to complete the proof, we have to prove that
	\[\dom\left(\left(\sbmat{I}{0}{0}{\Refl_T}\frakI_{F,T}\right)^*\right)\subset \dom\left(\Refl_T\frakO_{\mathrm{d},F^*,T}\right).\]
	Assume that $\spvek{z}{y_\mathrm{d}}\in \dom\left(\left(\sbmat{I}{0}{0}{\Refl_T}\frakI_{F,T}\right)^*\right)$. Then, by \eqref{eq:adjcomp}, we have for all \changed{$u\in H^2_{0}([0,T];U)$} that
	\[\begin{aligned}
	\left\langle u,\left(\sbmat{I}{0}{0}{\Refl_T}\frakI_{F,T}\right)^*\spvek{z}{y_\mathrm{d}}\right\rangle_{L^2([0,T];U)}
	&=
	\left\langle \sbmat{I}{0}{0}{\Refl_T}\frakI_{F,T}u,\spvek{z}{y_\mathrm{d}}\right\rangle_{L^2([0,T];\changed{Y})}\\
	&\!\!\!\stackrel{\eqref{eq:adjcomp}}{=}\left\langle u, \Refl_T\big({\frakC_{\mathrm{d},T}}F^*z+{\frakD_{\mathrm{d},T}}y_\mathrm{d}\big)\right\rangle_{H^{2}_{0l}([0,T];U),H^{-2}_{0r}([0,T];U)} .
	\end{aligned}
	\]
	This yields that $\Refl_T\big({\frakC_{\mathrm{d},T}}F^*z+{\frakD_{\mathrm{d},T}}y_\mathrm{d}\big)\in L^2([0,T];U)$, and thus
	\[\spvek{z}{y_\mathrm{d}}\in\dom\left(\Refl_T\frakO_{\mathrm{d},F^*,T}\right),\]
	which completes the proof.
\end{proof}

\section{The optimal control problem}\label{sec:optcont}
\noindent Having developed the operator-theoretic foundation, we are now prepared to analyze the optimal control problem
\begin{equation}
\begin{aligned}
\textrm{minimize}&\quad\frac{1}{2}\int_0^T\|y(t)-y_\mathrm{ref}(t) \|^2_{Y}{\rm d}t+ \frac{1}{2}\| Fx(T)-z_f\|_Z^2\\
\text{subject to}&\quad\spvek{\dot{x}(t)}{y(t)} = \sbvek{A\& B\\[-1mm] }{C\& D} \spvek{{x}(t)}{u(t)}, \quad
x(0) = x_0,\quad u \in \Uad.\nonumber
\end{aligned}\label{eq:OCP}\tag{OCP}
\end{equation}
Hereby, $\Uad\subset L^2([0,T];U)$ is the {\em set of admissible inputs} which, as will be illustrated later, allows to incorporate common control constraints of various types.

We first collect all the central assumptions on the problem~\eqref{eq:OCP}.
\begin{ass}\label{ass:ocp}\
	\begin{enumerate}[(a)]
		\item\label{ass:ocp1} $S=\sbvek{A\& B}{C\& D}$ is a~system node on the triple $(X,U,Y)$ of complex Hilbert spaces,
		\item\label{ass:ocp4} $F\in L(X,Z)$ for some Hilbert space $Z$, and, for the semigroup $\frakA:\R_{\ge0}\to L(X)$ generated by~$A$,
		the set $\mathcal{B}_F$ as in \eqref{eq:dense1}
		is dense in $Z$.
		\item\label{ass:ocp2} $x_0\in X_{-1}$, $z_f\in Z$, and $y_\mathrm{ref}\in L^{2}([0,T];Y)$.
		\item\label{ass:ocp3} $\Uad$ is a~closed and convex subset of $L^2([0,T];U)$.
		\item\label{ass:ocp5} At least one of the following two conditions hold:
		\begin{enumerate}[(i)]
			\item $\Uad$ is bounded.
			\item The cost functional is {\em coercive}. That is, there exists some $c>0$, such that all generalized \changed{trajectories} $(x,u,y)$ with $x(0)=0$ satisfy
			\[\|y\|_{L^2([0,T];U)}^2+\| Fx(T)\|_Z\geq c\,\|u\|_{L^2([0,T];U)}^2.\]
		\end{enumerate}
		\item\label{ass:ocp6} For the mappings introduced in \eqref{eq:genSOmap}, \eqref{eq:genIOmap} and Definition~\ref{def:FIS} (see also Table~\ref{tab:ourops}),
		there exists an admissible input $\hat{u}\in\Uad$ with 
		\[\spvek{x_0}{\hat{u}}\in \dom(\frakT_{F,T})\quad\text{ and }\quad
		{\frakC_{T}}{x_0}+{\frakD_T}{\hat{u}}\in L^2([0,T];Y).\]
	\end{enumerate}
\end{ass}
In the following, we discuss our assumptions.
\begin{rem}[Our assumptions]\label{rem:optass}\
	\begin{enumerate}[(a)]
		\item We indeed permit initial values $x_0$ in the extrapolation space $X_{-1}$. Of course, this includes initialization with $x_0\in X$.
		\item A~possible choice for set of admissible inputs is, for example, by means of (a.e.) pointwise control constraints, i.e., for some closed and convex set $U_{\rm ad}\subset U$,
		\[\Uad=\setdef{u\in L^2([0,T];U)}{u(t)\in U_{\rm ad}\text{ for almost all }t\in[0,T]}.\]
		This for instance allows to incorporate box constraints on the input. However, in \eqref{eq:OCP}, also time-varying constraints are possible, such as $u(t)\in U_{\rm ad}(t)$ for almost all $t\in[0,T]$, where $U_{\rm ad}(t)$ is closed and convex for all
		$t\in[0,T]$.
		\item Condition \eqref{ass:ocp4} in Assumptions~\ref{ass:ocp} guarantees, in view of Proposition~\ref{prop:closedop}, that the terminal cost
		\begin{equation}
		\frac12\|Fx(T)-z_f\|_Z^2 \label{eq:termcost}   \end{equation}
		is, in the weak sense, well-defined by $\|\frakT_{F,T}\spvek{x_0}{u}-z_f\|_Z^2$, if $\spvek{x_0}{\hat{u}}\in \dom(\frakT_{F,T})$ and $\infty$ otherwise.
		\item Condition \eqref{ass:ocp5} in Assumptions~\ref{ass:ocp} guarantees that any sequence of admissible inputs for which the cost functional tends to the infimal value, is bounded in $L^2([0,T];U)$ by means of the usual argumentation via closed and bounded sublevel sets. As is standard in optimal control theory, this gives rise to the existence of a~weakly convergent subsequence, whose limit will be shown to be the optimal control.\\
		Note that coercivity is, for instance, fulfilled, if the cost functional is of type \eqref{eq:costfun2} for some $c>0$.
		\item Condition \eqref{ass:ocp6} in Assumptions~\ref{ass:ocp} is a~crucial one, since it essentially expresses that there exists at least one control with finite cost.\\
		A~sufficient criterion for this is the existence of some generalized trajectory $(\hat{x},\hat{u},\hat{y})$ with $\hat{x}(0)=x_0$ and $\hat{u}\in\Uad$. This is, for instance, fulfilled if $\Uad$ is nonempty, 
		$y_{\rm ref}\in L^2([0,T];Y)$, $x_0,x_f\in X$,
		and the system is well-posed in the sense that there exists some $c>0$, such that all classical (and thus also generalized) trajectories fulfill \eqref{eq:wp}. 
	\end{enumerate}
\end{rem}




Now we define what we mean by an optimal control.
\begin{defn}[Optimal control]\label{def:cost}
	Consider an optimal control problem \eqref{eq:OCP} with Assumptions~\ref{ass:ocp}.
	Then, for the operators introduced in Section~\ref{sec:sysnode} (see also the table in Appendix~\ref{sec:ourops}), the cost of the input $u\in\Uad$ is given by
	\[\mathcal{J}(u)=\begin{cases}
	{\displaystyle\frac12}\left\|\pvek{\frakT_{F,T}\spvek{x_0}u-z_f}{{\frakC_T}{x_0}+{\frakD_T}{{u}}-y_{\rm ref}}\right\|^2_{Z\times L^2([0,T];Y)}\!\!\!\!\!&:\;\text{\parbox{63mm}{$\spvek{x_0}u\in\dom(\frakT_{F,T})\;\;\wedge\,$\\${\frakC_T}{x_0}+{\frakD_T}{{u}}\in L^2([0,T];Y),$}}\\
	\infty&:\;\text{\parbox{63mm}{$\spvek{x_0}u\notin\dom(\frakT_{F,T})\;\;\vee\,$\\${\frakC_T}{x_0}+{\frakD_T}{{u}}\notin L^2([0,T];Y)$.}}
	\end{cases}\]
	We call $u_{\rm opt}\in\Uad$ an {\em optimal control for \eqref{eq:OCP}}, if
	\[\mathcal{J}(u_{\rm opt})=\inf_{u\in\Uad}\mathcal{J}(u).\]
\end{defn}

\begin{rem}\label{rem:costfun}\
	\begin{enumerate}[(i)]
		\item If $x_0\in X$, then the cost functional can be rewritten to
		\[\mathcal{J}(u)=\begin{cases}
		{\displaystyle\frac12}\left\|\pvek{F\frakA(T)x_0+\frakT_{F,T}\spvek0u-z_f}{\frakO_{I,T}\spvek{x_0}{u}-y_\mathrm{ref}}\right\|^2_{Z\times L^2([0,T];Y)}\!\!\!\!\!\!\!&:\;\text{\parbox{40mm}{$\spvek{0}u\in\dom(\frakT_{F,T})\,\wedge\,$\\$\spvek{x_0}{u}\in \dom(\frakO_{I,T}),$}}\\
		\infty&:\;\text{\parbox{40mm}{$\spvek{0}u\notin\dom(\frakT_{F,T})\,\vee\,$\\$\spvek{x_0}{u}\notin \dom(\frakO_{I,T}).$}}
		\end{cases}\]
		\item Assume that $u,\delta u\in L^2([0,T];U)$, such that both
		$\mathcal{J}(u)<\infty$ and $\mathcal{J}(u+\delta u)<\infty$. Then the definition of the $F$-input map \changed{(Definition~\ref{def:FIO})} yields
		\begin{align}\label{eq:difference}
		\changed{\pvek{\frakT_{F,T}\spvek{x_0}{u+\delta u} -z_f}{{\frakC_T}{x_0}+{\frakD_T}{(u+\delta u)}-y_{\rm ref}} - \pvek{\frakT_{F,T}\spvek{x_0}u -z_f}{{\frakC_T}{x_0}+{\frakD_T}{{u}}-y_{\rm ref}}}
		=    \pvek{\frakT_{F,T}\spvek{0}{\delta u}}{{\frakD_{T}}\delta u}
		=\frakI_{F,T}\delta u.
		\end{align}
		In particular, for
		$u\in L^2([0,T];U)$ with
		$\mathcal{J}(u)<\infty$, we have
		\[\mathcal{J}(u+\delta u)<\infty\;\Leftrightarrow \delta u\in \dom(\frakI_{F,T}).\]
		By an expansion of \changed{the squared norm in} $\mathcal{J}(u+\delta u)$ \changed{and \eqref{eq:difference}}, we obtain that
		\begin{align*}
		&\phantom{=}\mathcal{J}(u+\delta u)-\mathcal{J}(u)\\
		&=\re\left\langle\pvek{{\frakT_{F,T}\spvek{x_0}u-z_f}}{{\frakC_T}{x_0}+{\frakD_T}u-y_{\rm ref}},\frakI_{F,T}\delta u\right\rangle_{Z\times L^2([0,T];Y)}+\frac12 \left\|\frakI_{F,T}\delta u\right\|^2_{X\times L^2([0,T];Y)}.
		\end{align*}
		This will be later on used for our optimality condition based on the derivative of the cost functional.
	\end{enumerate}
\end{rem}
Now we collect some important properties of the cost functional.
\begin{prop}\label{prop:cost}
	Consider an optimal control problem \eqref{eq:OCP} with Assumptions~\ref{ass:ocp}, then the cost functional
	\[\mathcal{J}:\Uad\to \R_{\ge0}\cup\{\infty\}\]
	as specified in Definition~\ref{def:cost} has the following properties:
	\begin{enumerate}[(a)]
		\item\label{prop:costa} $\mathcal{J}$ is {\em proper}. That is, there exists some $u\in\Uad$ with $\mathcal{J}(u)<\infty$.
		\item\label{prop:costb} $\mathcal{J}$ is {\em convex}. That is, for all $u_1,u_2\in\Uad$, $\lambda\in[0,1]$,
		\begin{equation}\mathcal{J}(\lambda u_1+(1-\lambda)u_2)\leq \lambda\mathcal{J}(u_1)+(1-\lambda)\mathcal{J}(u_2).\label{eq:convex}\end{equation}
		\item\label{prop:costc} $\mathcal{J}$ is {\em lower semicontinuous}. That is, for all $a\in\R_{\ge0}$, the {\em sublevel set}
		\begin{equation}
		\setdef{u\in\Uad}{\mathcal{J}(u)\leq a}\label{eq:sublevel}
		\end{equation}
		is closed subset of $L^2([0,T];U)$.
	\end{enumerate}
\end{prop}
\begin{proof}\
	Let $\hat{u}\in L^2([0,T];U)$ with $\mathcal{J}(\hat{u})<\infty$ (which exists by Assumptions~\ref{ass:ocp}\,\eqref{ass:ocp6}). Our proof is mainly based on the identity
	\begin{equation}
	\pvek{\frakT_{F,T}\spvek{x_0}u-z_f}{{\frakC_T}{x_0}+{\frakD_T}{u}-y_{\rm ref}}=\pvek{\frakT_{F,T}\spvek{x_0}{\hat{u}}-z_f}{{\frakC_T}{x_0}+{\frakD_T}\hat{u}-y_{\rm ref}}+\frakI_{F,T}(u-\hat{u}),\label{eq:costrel}
	\end{equation}
	which holds by Remark~\ref{rem:costfun}.
	\begin{enumerate}[(a)]
		\item This follows from $\mathcal{J}(\hat{u})<\infty$.
		\item Let $u_1,u_2\in \Uad$ and $\lambda \in [0,1]$. If at least one
		of $\mathcal{J}(u_1)$, $\mathcal{J}(u_2)$ is infinite, the inequality
		\eqref{eq:convex} is trivially fulfilled. Now assume that $\mathcal{J}(u_1), \mathcal{J}(u_2)\in\R_{\ge0}$, and let $\lambda\in[0,1]$. Then, by 
		linearity of the operators in Table~\ref{tab:ourops},
		\begin{multline*}
		\lambda\pvek{\frakT_{F,T}\spvek{x_0}{u_1}-z_f}{{\frakC_T}{x_0}+{\frakD_T}{u}_1-y_{\rm ref}}+(1-\lambda)\pvek{\frakT_{F,T}\spvek{x_0}{u_2}-z_f}{{\frakC_T}{x_0}+{\frakD_T}{u}_2-y_{\rm ref}}
	\\=    \pvek{\frakT_{F,T}\spvek{x_0}{\lambda u_1+(1-\lambda)u_2}-z_f}{{\frakC_T}{x_0}+{\frakD_T}(\lambda u_1+(1-\lambda)u_2)-y_{\rm ref}}.
		\end{multline*}
		Now using convexity of the mapping $\|\cdot\|_{X\times L^2([0,T];U)}^2$ (which is a~simple consequence of the Cauchy-Schwarz inequality), we can conclude that $\mathcal{J}$ is convex.
		\item Let $a\in \R_{\ge0}$, and let $(u_n)$ be a~sequence in $\Uad$ with $\mathcal{J}(u_n)<a$ for all $n\in\N$, which converges in $L^2([0,T];U)$ to $u\in L^2([0,T];U)$. Closedness of $\Uad$ leads to $u\in \Uad$. Then it follows from \eqref{eq:costrel} that
		\begin{equation}
		\big(\frakI_{F,T}(u_n-\hat{u})\big)\label{eq:imseq}
		\end{equation}
		is a~bounded sequence in $Z\times L^2([0,T];Y)$. Hence, by \cite[Thm.~8.10]{Alt16}, it has a~weakly convergent subsequence. It is therefore no loss of generality
		to assume that the sequence \eqref{eq:imseq} itself is weakly convergent. As a~consequence, there exist $z\in Z$, $y\in L^2([0,T];Y)$, such that weak convergence
		\begin{align}\label{eq:weakconv}
		\begin{split}&\left(\pvek{\hat{u}}{\spvek{\frakT_{F,T}\spvek{x_0}{\hat{u}}-z_f}{{\frakC_T}{x_0}+{\frakD_T}{\changed{\hat{u}}}-y_{\rm ref}}}+\pvek{u_n-\hat{u}}{\frakI_{F,T}(u_n-\hat{u})}\right)\\&\hspace*{3cm}=
		\left(\pvek{u_n}{\spvek{\frakT_{F,T}\spvek{x_0}{\hat{u}}-z_f}{{\frakC_T}{x_0}+{\frakD_T}{\changed{\hat{u}}}-y_{\rm ref}}+\frakI_{F,T}(u_n-\hat{u})}\right)
		\rightharpoonup\pvek{u}{\spvek{z}{y}}
		\end{split}
		\end{align}
		holds in $L^2([0,T];U)\times Z\times L^2([0,T];Y)$. The latter sequence evolves in the affine-linear space
		\begin{equation}
		\pvek{\hat{u}}{\spvek{\frakT_{F,T}\spvek{x_0}{\hat{u}}-z_f}{{\frakC_T}{x_0}+{\frakD_T}{\changed{\hat{u}}}-y_{\rm ref}}}+G(\frakI_{F,T})\label{eq:affspace}
		\end{equation}
		where $G(\frakI_{F,T})$ stands for the graph of $\frakI_{F,T}$. Since $G(\frakI_{F,T})$ is a~closed space by closedness of the operator $\frakI_{F,T}$ (see Proposition~\ref{prop:closedop})
		the affine-linear space \eqref{eq:affspace} is closed and convex, whence it is also weakly closed by the separation theorem \cite[Thm.~8.12]{Alt16}. This yields that
		\[\pvek{z}{y}=\spvek{\frakT_{F,T}\spvek{x_0}{\hat{u}}-z_f}{{\frakC_T}{x_0}+{\frakD_T}{u}-y_{\rm ref}}+\frakI_{F,T}(\changed{u}-\hat{u}).\]
		\changed{As the norm function is weakly lower semicontinuous, cf.~\cite[Remark 8.3]{Alt16}, we have}
		\begin{align*}
		\mathcal{J}(u)=\left\|\spvek{z}y\right\|_{Z\times L^2([0,T];Y)}^2&\leq\liminf_{n\to\infty}\left\|\spvek{\frakT_{F,T}\spvek{x_0}{\hat{u}}-z_f}{{\frakC_T}{x_0}+{\frakD_T}{u}-y_{\rm ref}}+\frakI_{F,T}(u_n-\hat{u})\right\|_{Z\times L^2([0,T];Y)}^2\\&=\liminf_{n\to\infty}\mathcal{J}(u_n)\leq a.
		\end{align*}
		As a~consequence,  $u$ belongs to the sublevel set \eqref{eq:sublevel}, which shows the desired result.
		%
	\end{enumerate}
\end{proof}
Next we show existence of optimal controls.
\begin{thm}[Existence of optimal controls]\label{thm:optcont}
	Consider an optimal control problem \eqref{eq:OCP} with Assumptions~\ref{ass:ocp}. Then there exists an optimal control. 
\end{thm}
\begin{proof}
	Here we follow the traditional argumentation for minimization of convex functions, such as, for instance, presented in \cite[Chap.~II]{EkelTema99}: Let $(u_n)$ be a~sequence in $\Uad$, such that
	\[\lim_{n\to\infty}\mathcal{J}(u_n)=V_{\rm opt}:=\inf_{u\in\Uad}\mathcal{J}(u).\]
	Since it is assumed that $\Uad$ is bounded or $\mathcal{J}$ is coercive, we can conclude from each of these two cases that $(u_n)$ is a~bounded sequence in $L^2([0,T];U)$.
	\changed{Consequently, $(u_n)$ has a~weakly convergent subsequence, whose weak limit is denoted by $u_{\rm opt}\in L^2([0,T];U)$.} Now we show that $u_{\rm opt}$ is indeed an optimal control. We first note that, by using that $\Uad$ is closed and convex, we can conclude from the separation theorem that it is weakly closed, i.e., $u_{\rm opt}\in\Uad$. Convexity of the cost functional (as shown in Proposition~\ref{prop:cost}\,\eqref{prop:costb}) implies that each sublevel set is convex. Combining this with lower semicontinuity of cost functional (which has been proven in Proposition~\ref{prop:cost}\,\eqref{prop:costc}), we can again use the separation theorem to see that each sublevel set is weakly closed. Consequently, we have
	\[\forall\,\varepsilon>0:\quad u_{\rm opt}\in \setdef{u\in\Uad}{\mathcal{J}(u)\leq V_{\rm opt}+\varepsilon},\]
	and thus $\mathcal{J}(u_{\rm opt})=V_{\rm opt}$. In other words, $u_{\rm opt}$ is an optimal control.
\end{proof}

\noindent The next result provides a~characterization of optimality for controls.
\begin{thm}\label{thm:optcond}
	Consider an optimal control problem \eqref{eq:OCP} with Assumptions~\ref{ass:ocp}. Then $u_{\rm opt}\in\Uad$ is an optimal control, if, and only if,
	\begin{multline}
	\forall\, u \in \Uad\ \mathrm{ with }\, u-u_\mathrm{opt}\in \dom(\frakI_{F,T}):\\
	\re\left\langle\pvek{\frakT_{F,T}\spvek{x_0}{{u}_{\rm opt}}-z_f}{{\frakC_T}{x_0}+{\frakD_T}{u}_{\rm opt}-y_{\rm ref}},\frakI_{F,T}(u-u_{\rm opt})\right\rangle_{X\times L^2([0,T];Y)} \geq 0.\label{eq:optder}
	\end{multline}
\end{thm}
\begin{proof}
	To prove the implication ``$\Rightarrow$'', let $u_\mathrm{opt}\in \dom(\frakI_{F,T})\cap\Uad$ be an optimal control. Further, let $u\in\Uad$ with $\delta u:=u-u_\mathrm{opt}\in \dom(\frakI_{F,T})$, and $\lambda\in[0,1]$. Convexity of $\Uad$ and $\dom(\frakI_{F,T})$ yields
	\[u_\mathrm{opt}+\lambda \delta u\in \dom(\frakI_{F,T})\cap\Uad,\]
	and, by
	invoking Remark~\ref{rem:costfun},
	and optimality of $u_\mathrm{opt}$, we have
	\begin{align*}
	0&\leq \tfrac1\lambda \big(\mathcal{J}(u_\mathrm{opt}+\lambda \delta u)-\mathcal{J}(u_\mathrm{opt})\big)\\
	&=\re\left\langle\pvek{\frakT_{F,T}\spvek{x_0}{{u}_{\rm opt}}-z_f}{{\frakC_T}{x_0}+{\frakD_T}{u}_{\rm opt}-y_{\rm ref}},\frakI_{F,T}\delta u\right\rangle_{Z\times L^2([0,T];Y)} +\frac{\lambda}2 \left\|\frakI_{Z,T}\delta u\right\|^2_{Z\times L^2([0,T];Y)}\\
	&\underset{\lambda\searrow0}{\rightarrow}\re\left\langle\pvek{\frakT_{F,T}\spvek{x_0}{{u}_{\rm opt}}-z_f}{{\frakC_T}{x_0}+{\frakD_T}{u}_{\rm opt}-y_{\rm ref}},\frakI_{F,T}\delta u\right\rangle_{Z\times L^2([0,T];Y)}.
	\end{align*}
	Now we show the converse implication: Assume that $u_\mathrm{opt}\in\Uad$ fulfills
	\eqref{eq:optder}. Let $u\in\Uad$ with $\mathcal{J}(u)<\infty$. Now using Remark~\ref{rem:costfun}, we obtain for $\delta u=u_{\rm opt}-u$ that $\delta\in \dom(\frakI_{F,T})$ with
	\begin{multline*}
	\mathcal{J}(u)-\mathcal{J}(u_{\rm opt})
	\\=\re\left\langle\pvek{\frakT_{F,T}\spvek{x_0}{{u}_{\rm opt}}-z_f}{{\frakC_T}{x_0}+{\frakD_T}{u}_{\rm opt}-y_{\rm ref}},\frakI_{F,T}\delta u\right\rangle_{Z\times L^2([0,T];Y)}+\frac12 \left\|\frakI_{F,T}\delta u\right\|^2_{Z\times L^2([0,T];Y)}\geq 0.
	\end{multline*}
	This shows that $\mathcal{J}(u)\ge \mathcal{J}(u_{\rm opt})$, i.e., the control $u_{\rm opt}$ is optimal.
	Hence, the claim is proven.
\end{proof}
Next we present a~sufficient criterion for uniqueness of optimal controls.
\begin{thm}[Uniqueness of optimal controls]\label{thm:optcont2}
	Consider an optimal control problem \eqref{eq:OCP} satisfying Assumptions~\ref{ass:ocp} and, additionally, let the $F$-input map $\frakI_{F,T}$ \changed{defined in Definition~\ref{def:FIO}} be injective.
	Then there exists exactly one optimal control.
\end{thm}
\begin{proof}
	Existence of at least one optimal control has already been proven in Theorem~\ref{thm:optcont}. Now, \changed{assume} that the $F$-input map is injective, and that $u_{\rm opt,1},u_{\rm opt,2}\in\Uad$  are optimal controls.
	Then, by Remark~\ref{rem:costfun}, we have $u_{\rm opt,1}-u_{\rm opt,2}\in\dom(\frakI_{F,T})$ and
	\begin{align*}
	0&{=}\mathcal{J}(u_{\rm opt,1})-\mathcal{J}(u_{\rm opt,2})\\
	&=\re\left\langle\pvek{\frakT_{F,T}\spvek{x_0}{{u}_{\rm opt,2}}-z_f}{{\frakC_T}{x_0}+{\frakD_T}u_{\rm opt,2}-y_{\rm ref}},\frakI_{F,T}(u_{\rm opt,1}-u_{\rm opt,2})\right\rangle_{X\times L^2([0,T];Y)}\\&\hspace{3cm}+\frac12 \left\|\frakI_{F,T}(u_{\rm opt,1}-u_{\rm opt,2})\right\|^2_{Z\times L^2([0,T];Y)}.
	\end{align*}
	The first summand on the right hand side is nonnegative by Theorem~\ref{thm:optcond}, 
	which means that both summands vanish.
	This gives
	\[\frakI_{F,T}(u_{\rm opt,1}-u_{\rm opt,2})=0,\]
	and injectivity of $\frakI_{F,T}$ now leads to $u_{\rm opt,1}=u_{\rm opt,2}$, and the result is shown.
	
\end{proof}
\noindent If $u-u_{\rm opt}\in H^2_{0l}([0,T];U)$, then $\frakI_{F,T}(u-u_{\rm opt})=\spvek{Fz(T)}{w}\in Z\times L^2([0,T];Y)$, where
\[\spvek{\dot{z}(t)}{w(t)}
= \sbvek{A\&B\\[-1mm]}{C\&D} \spvek{{z}(t)}{u(t)-u_{\rm opt}(t)},\quad z(0)=0.\]
This will be used in the following result, which is a~direct consequence of combining the Theorem~\ref{thm:optcond} with Proposition~\ref{prop:adjFinput}. The proof is therefore omitted.


\begin{cor}\label{cor:aufloesen}
	Consider an optimal control problem \eqref{eq:OCP} with Assumptions~\ref{ass:ocp}, and let $u_{\rm opt}\in\Uad$ be an optimal control (which exists by Theorem~\ref{thm:optcont}). 
	Then
	\begin{align*}
	\forall\, u \in \Uad\ \mathrm{ with }\, u-u_\mathrm{opt}\in H^2_{0l}([0,T];U):\quad 
	\re\left\langle u-u_{\rm opt}, u_\mathrm{d}\right\rangle_{H^2_{0l}([0,T];U),H^{-2}_{0r}([0,T];U)} \geq 0,
	\end{align*}
	where 
	\begin{equation}
	u_\mathrm{d} = \Refl_T [{\frakC_{\mathrm{d},T}},\,{\frakD_{\mathrm{d},T}}]  \pvek{F^*\frakT_{F,T}\spvek{x_0}{{u}_{\rm opt}}-z_f}{\Refl_T\left({\frakC_T}{x_0}+{\frakD_T}{u}_{\rm opt}-y_{\rm ref}\right)}.\label{eq:udopt}
	\end{equation}
\end{cor}
\noindent By defining the 
state and output of the system driven by the optimal control
\begin{equation}
\begin{aligned}
x_{\rm opt}(t)&=\frakA_{-1}(t)x_0+{\frakB_t}{u}_{\rm opt},\quad t\in[0,T],\\
y_{\rm opt}&={\frakC_T}{x_0}+{\frakD_T}{u}_{\rm opt},
\end{aligned}\label{eq:yopt}\end{equation}
and the {\em adjoint state}
\begin{equation}
\mu(t)=\frakA_{T-t}(t)^*F^*\frakB_{\bar{x},F,T}u_{\rm opt}+{\frakB_{\mathrm{d},T-t}}(y_{\rm opt}-y_{\rm ref}),\quad t\in[0,T],
\label{eq:adjstate}\end{equation}
we obtain that $u_\mathrm{d}$ in \eqref{eq:udopt} solves, in a~certain sense, 
the boundary value problem
\begin{equation}
\begin{aligned}
\spvek{\dot{x}_{\rm opt}(t)}{y_{\rm opt}(t)}&
= \sbvek{A\&B\\[-1mm]}{C\&D} \spvek{{x}_{\rm opt}(t)}{u_{\rm opt}(t)},& x_{\rm opt}(0)&=x_0,\\
\spvek{\dot{\mu}(t)}{u_\mathrm{d}(t)}
&= \sbvek{{-[A\&B]^\mathrm{d}}}{{\phantom{-}[C\&D]^\mathrm{d}}}
\spvek{{\mu}(t)}{y_{\rm opt}(t)-y_{\rm ref}(t)},& \mu(T)&=F^*(Fx_{\rm opt}(T)-z_f).
\end{aligned}\label{eq:udsys2}
\end{equation}

Let us now \changed{treat} the case where, loosely speaking, the constraints imposed by the set $\Uad$ of admissible controls are not active anywhere for the optimal control $u_{\rm opt}$. First observe that, if $u_{\rm opt}\pm\delta u\in \Uad$ for some $\delta u\in \dom(\frakI_{F,T})$, we obtain from Theorem~\ref{thm:optcond} that
\begin{equation}
\re\left\langle\pvek{\frakT_{F,T}\spvek{x_0}{{u}_{\rm opt}}-z_f}{{\frakC_T}{x_0}+{\frakD_T}{u}_{\rm opt}-y_{\rm ref}},\frakI_{F,T}\delta u\right\rangle_{Z\times L^2([0,T];Y)} = 0.\label{eq:reopt}
\end{equation}
This will be used in the following, where we give a~criterion for the variational inequality in Theorem~\ref{thm:optcond} becoming an equality.
\begin{prop}\label{prop:aufloesen}
	Consider an optimal control problem \eqref{eq:OCP} with Assumptions~\ref{ass:ocp}, and let $u_{\rm opt}\in\Uad$ be an optimal control (which exists by Theorem~\ref{thm:optcont}). Assume that the closure of the set
	\begin{align}
	\mathcal{T}:= \changed{\{\delta u\in \dom(\frakI_{F,T})\,|\, u_\mathrm{opt} \pm \delta u \in \Uad\}} \subset \dom(\frakI_{F,T})\label{eq:T}
	\end{align}
	has a~nonempty interior in $L^2([0,T];U)$, i.e.,
	\[\operatorname{int}\overline{\mathcal{T}}\neq\emptyset.\]
	Then
	\begin{align}\label{eq:aufloesen}
	\frakO_{\mathrm{d},F^*,T}\pvek{\frakT_{F,T}\spvek{x_0}{{u}_{\rm opt}}-z_f}{\Refl_T\big({\frakC_T}{x_0}+{\frakD_T}{u}_{\rm opt}-y_{\rm ref}\big)}=0.
	\end{align}
\end{prop}
\begin{proof}
	Assume that, for $\varepsilon>0$, $\delta u\in L^2([0,T];U)$, the $\varepsilon$-ball $U_\varepsilon(\delta u)$ centered at $\delta u$ is contained in $\operatorname{int}\overline{\mathcal{T}}$. Since,by definition, $\mathcal{T}=-\mathcal{T}$, convexity of $(\Uad-u_\mathrm{opt})\cap\dom(\frakI_{\changed{F},T})$, yields
	\[U_\varepsilon(0)\subset\operatorname{int}\overline{\mathcal{T}}.\]
	Consequently, the space $\operatorname{span}_\C\mathcal{T}$, \changed{the linear span with complex coefficiencs}, is dense in $L^2([0,T];U)$.\\
	The argumentation prior to this proposition that \eqref{eq:reopt} holds for all $\delta u\in\mathcal{T}$. Consequently, \eqref{eq:reopt} holds for all $\delta u$ in the dense subspace
	$\operatorname{span}_\C\mathcal{T}$. Hence, by multiplying with the imaginary unit, we further obtain that, for all $\delta u\in\operatorname{span}_\C\mathcal{T}$,\begin{multline*}  \IM\left\langle\pvek{\frakT_{F,T}\spvek{x_0}{{u}_{\rm opt}}-z_f}{{\frakC_T}{x_0}+{\frakD_T}{u}_{\rm opt}-y_{\rm ref}},\frakI_{F,T}\delta u\right\rangle_{Z\times L^2([0,T];Y)} \\
	=\imath \re\left\langle\pvek{\frakT_{F,T}\spvek{x_0}{{u}_{\rm opt}}-z_f}{{\frakC_T}{x_0}+{\frakD_T}{u}_{\rm opt}-y_{\rm ref}},\frakI_{F,T}(\imath\delta u)\right\rangle_{Z\times L^2([0,T];Y)}
	= 0.\end{multline*}
	Hence, for all $\delta u\in\operatorname{span}_\C\mathcal{T}$,
	\[\left\langle\pvek{\frakT_{F,T}\spvek{x_0}{{u}_{\rm opt}}-z_f}{{\frakC_T}{x_0}+{\frakD_T}{u}_{\rm opt}-y_{\rm ref}},\frakI_{F,T}\delta u\right\rangle_{Z\times L^2([0,T];Y)} = 0.\]
	This shows that
	\begin{align*}
	\pvek{\frakT_{F,T}\spvek{x_0}{{u}_{\rm opt}}-z_f}{{\frakC_T}{x_0}+{\frakD_T}{u}_{\rm opt}-y_{\rm ref}} &\in \dom(\frakI_{F,T}^*).
	\end{align*}
	with
	\[\frakI_{F,T}^*\pvek{\frakT_{F,T}\spvek{x_0}{{u}_{\rm opt}}-z_f}{{\frakC_T}{x_0}+{\frakD_T}{u}_{\rm opt}-y_{\rm ref}}=0.\]
	Then the result follows from the representation of $\frakI_{F,T}^*$ in Proposition~\ref{prop:adjFinput}.
\end{proof}
Under the assumptions Proposition~\ref{prop:aufloesen}, we have that $u_\mathrm{d}$ from Corollary~\ref{cor:aufloesen} vanishes.
That is, by considering the state $x_{\rm opt}$ and output $y_{\rm opt}$ of the system driven by the optimal control (see \eqref{eq:yopt}) and 
the adjoint state $\mu$ as in \eqref{eq:adjstate}, we are led to a~boundary value problem
\begin{equation}
\begin{aligned}
\spvek{\dot{x}_{\rm opt}(t)}{y_{\rm opt}(t)}&
= \sbvek{A\&B\\[-1mm]}{C\&D} \spvek{{x}_{\rm opt}(t)}{u_{\rm opt}(t)},& x_{\rm opt}(0)&=x_0,\\
\spvek{\dot{\mu}(t)}{0}
&= \sbvek{{-[A\&B]^\mathrm{d}}}{{\phantom{-}[C\&D]^\mathrm{d}}}
\spvek{{\mu}(t)}{y_{\rm opt}(t)-y_{\rm ref}(t)},& \mu(T)&=F^*(Fx_{\rm opt}(T)-z_f).
\end{aligned}\label{eq:udsys0}
\end{equation}
\begin{rem}[Regular systems and optimality Hamiltonians]
	\changed{Using the fact that $X$, $U$, and $Y$ are Hilbert spaces, it has been shown in \cite[Thm.~5.1.12]{Staffans2005} that the system nodes $\sbvek{A\&B}{C\&D}$ describing a well-posed system are "compatible" in the sense that, roughly speaking, $C\&D$ can be decomposed into two parts: one associated with the state and the other with the input. (Notably, well-posedness is not required for the proof of \cite[Thm.~5.1.12]{Staffans2005}.)} More precisely, for 
	\begin{equation}\label{eq:Vdef}
	V:=\setdef{x\in X}{\exists\ u\in U\text{ s.t.\ }\spvek{x}u\in\dom(A\&B)}
	\end{equation}
	which is a~Hilbert space endowed with the norm
	\[\|x\|_V:=\inf\setdef{\|\spvek{x}u\|_{\dom(A\&B)}}{\exists\ u\in U\text{ s.t.\ }\spvek{x}u\in\dom(A\&B)},\]
	we have that there exists some $\tilde{C}\in L(V,Y)$, $D\in L(U,Y)$, such that $C\&D\spvek{x}{u}=\tilde{C}x+Du$ for all $\spvek{x}{u}\in\dom(A\&B)$. Unfortunately, $\tilde{C}$ is neither uniquely determined by $C\&D$ (for instance, controlled boundary values can be put into $D$ and in $\tilde{C}$), nor do we have
	\begin{equation}\sbvek{{[A\&B]}^\mathrm{d}\\[-1mm]}{{[C\&D]}^\mathrm{d}}=
	\sbmat{A^*}{\tilde{C}^*\\[-1mm]}{B^*}{D^*},
	\label{eq:dualreg}\end{equation}
	in general. However, in case where the system is {\em regular} (that is, its transfer function has a~limit on the positive real axis, see \cite{SW04,Staffans2005}), then $\tilde{C}$ and $D$ can be chosen in a~way that \eqref{eq:dualreg} holds, see \cite[Thm.~3.5]{SW04}. In this case we obtain from \eqref{eq:udsys0} that 
	\begin{align}\label{eq:someequation}
	\begin{split}
	0&=B^*\mu(t)+D^*(\tilde{C}x_{\rm opt}(t)+Du_{\rm opt}(t))-D^*y_{\rm ref}(t)
	\\&=D^*\tilde{C}x_{\rm opt}(t)+B^*\mu(t)+D^*Du_{\rm opt}(t)-D^*y_{\rm ref}(t).
	\end{split}
	\end{align}
	If, additionally, $D$ has trivial nullspace and closed range, then $D^*D$ has a~bounded inverse, and thus we can solve \changed{\eqref{eq:someequation}} for $u_{\rm opt}(t)$, which gives
	\[u_{\rm opt}(t)=-(D^*D)^{-1}D^*\tilde{C}x_{\rm opt}(t)-(D^*D)^{-1}B^*\mu(t)+(D^*D)^{-1}D^*y_{\rm ref}(t).
	\]
	Now we can eliminate the control and output in \eqref{eq:udsys0}, to obtain, at least formally, the {\em abstract Hamiltonian system}
	\begin{multline*}
	\pvek{\dot{x}(t)}{\dot{\mu}(t)}=\begin{bmatrix}A-B(D^*D)^{-1}D^*\tilde{C}&-B(D^*D)^{-1}B^*\\-\tilde{C}^*D^*(D^*D)^{-1}D\tilde{C}&-A^*+\tilde{C}^*D(D^*D)^{-1}B^*\end{bmatrix}\pvek{{x}_{\rm opt}(t)}{{\mu}(t)}\\+\bvek{B(D^*D)^{-1}D^*}{\tilde{C}^*(I-D(D^*D)^{-1}D^*)}y_{\rm ref}(t),\qquad x_{\rm opt}(0)=x_0,\quad \mu(T)=F^*(Fx_{\rm opt}(T)-z_f).
	\end{multline*}
	Operators of the aforementioned type are analyzed in \cite{TreWy14} through the lens of spectral theory, with applications to algebraic Riccati equations that arise in linear-quadratic optimal control over an infinite time horizon.
\end{rem}
\begin{rem}[Real spaces and operators]\label{rem:realproblems}
	We assume that all spaces are complex. Real problems can be addressed by complexifying the operators and spaces involved. Importantly, coercivity of the cost functional is preserved under complexification.
\end{rem}

\section{The benefit of input penalization}\label{sec:inputpen}
Here, we examine optimal control problems with an added quadratic control penalization and elucidate the supplementary advantages.
That is, we consider the problem
\begin{equation*}
\begin{aligned}
\textrm{minimize}&\quad\frac{1}{2}\int_0^T\|y_1(t)-y_\mathrm{ref}(t) \|^2_{Y} + \|u(t)\|^2_U\,{\rm d}t+ \frac{1}{2}\|Fx(T)-z_f\|_Z^2\\
\text{subject to}&\quad\spvek{\dot{x}(t)}{y(t)} = \sbvek{A\& B\\[-1mm]}{C\& D} \spvek{{x}(t)}{u(t)}, \quad
x(0) = x_0,\quad u \in \Uad.\nonumber
\end{aligned}
\end{equation*}
As outlined in \eqref{eq:ODEnodeext}, we may equivalently reformulate this problem as
\begin{equation}
\begin{aligned}
\textrm{minimize}&\quad\frac{1}{2}\int_0^T\|y_{\rm ext}(t)-y_\mathrm{ref,ext}(t) \|^2_{Y_{\rm ext}}\,{\rm d}t+ \frac{1}{2}\|Fx(T)-z_f\|^2_Z\\
\text{subject to}&\quad\spvek{\dot{x}(t)}{y_{\rm ext}(t)} = \sbvek{A\&B}{\sbvek{C\&D}{0\ \ \Id}} \spvek{{x}(t)}{u(t)}, \quad
x(0) = x_0,\quad u \in \Uad.\nonumber
\end{aligned}
\end{equation}
with $Y_{\rm ext} = Y\times U$ and $y_\mathrm{ref,ext} = \spvek{y_{\mathrm{ref}}}{0}$. The coercivity condition of Assumption~\ref{ass:ocp}\eqref{ass:ocp5} is trivially fulfilled by
\begin{align*}
\|y_{\rm ext}\|_{L^2([0,T];Y_{\rm ext})}^2 = \|y\|_{L^2([0,T];Y)}^2 + \|u\|^2_{L^2([0,T];U)} \geq \|u\|^2_{L^2([0,T];U)}.
\end{align*}
As a~consequence, we do not have to impose the assumption of boundedness of $\Uad$ in $L^2([0,T];U)$.
Moreover, the adjoint system node of $S_{\rm ext}= \sbvek{A\&B}{\sbvek{C\&D}{0\ \ \Id}}$ on $(X,U,Y\times U)$ reads
\begin{equation}
S^*_{\rm ext}=\sbmat{[A\&B]^\mathrm{d}}{0}{[C\&D]^\mathrm{d} }{I},\label{eq:extadjsysnode}
\end{equation}
and it is again a system node, now on $(X,Y\times U,U)$. The corresponding maps \eqref{eq:genISmap}, \eqref{eq:genSOmap} of the system node $S=\sbvek{A\& B}{C\& D}$ will now be denoted as in the previous sections, whereas those corresponding to $S_{\rm ext}$ will be provided with the subscript $*_{\rm ext}$. We clearly have
\[{\frakB_{T,{\rm ext}}}={\frakB_{T}},\qquad \frakT_{F,T,{\rm ext}}=\frakT_{F,T}.\]
The state-to-output and input-to-output map of the system node $S$ split into
\begin{align*}
{\frakC_{T,{\rm ext}}} x_0 = \spvek{{\frakC_{T}}x_0}{0} \qquad {\frakD_{T,{\rm ext}}} u = \spvek{{\frakD_{T}}u}{u}.
\end{align*}
Hence, the input-to-state and state-to-output maps of the dual system nodes $S^*$ and $S_{\rm ext}^*$ are related by
\begin{align*}
{{\frakB_{\mathrm{d},F,T,{\rm ext}}}} \spvek{y_\mathrm{d}}{u} = {{\frakB_{\mathrm{d},F,T}}} y_\mathrm{d} \qquad {\frakD_{\mathrm{d},I,T,{\rm ext}}} \spvek{y_\mathrm{d}}{u} = {\frakD_{\mathrm{d},I,T}}y_\mathrm{d} + u.
\end{align*}
In particular, 
\begin{align*}
\dom({\frakB_{\mathrm{d},F,T,{\rm ext}}}) &= \dom({\frakB_{\mathrm{d},F,T}})\times L^2([0,T];U),\\
\dom({\frakD_{\mathrm{d},F,T,{\rm ext}}}) &= \dom({\frakD_{\mathrm{d},F,T}})\times L^2([0,T];U),
\end{align*}
and the condition \changed{\eqref{eq:aufloesen} of Proposition~\ref{prop:aufloesen}} may be rewritten in terms of
\begin{align*}
0 &= \frakO_{\mathrm{d},F^*,T,{\rm ext}} \pvek{\frakT_{F,T,{\rm ext}}\spvek{x_0}{{u}_{\rm opt}}-z_f}{\Refl_T\left({\frakC_{T,{\rm ext}}}{x_0}+{\frakD_{T,{\rm ext}}}{u}_{\rm opt}-y_{\rm ref,ext}\right)}\\
&= \frakO_{\mathrm{d},F^*,T,{\rm ext}} \pvek{\frakT_{F,T}\spvek{x_0}{{u}_{\rm opt}}-z_f}{\Refl_T\left(\spvek{{\frakC_{T}}{x_0}+{\frakD_{T}}{u}_{\rm opt}}{{u}_{\rm opt}}-\spvek{y_{\rm ref}}{0}\right)}\\
&= \frakO_{\mathrm{d},F^*,T} \pvek{\frakT_{F,T}\spvek{x_0}{{u}_{\rm opt}}-z_f}{\Refl_T\left({{\frakC_{T}}{x_0}+{\frakD_{T}}{u}_{\rm opt}}-y_\mathrm{ref}\right)}+\Refl_T u_{\rm opt}.
\end{align*}
Hence, we obtain a characterization of the control by means of the adjoint, i.e.,
\begin{align*}
u_\mathrm{opt} = -\Refl_T \frakO_{\mathrm{d},F^*,T} \spvek{\frakT_{F,T}\spvek{x_0}{{u}_{\rm opt}}-z_f}{\Refl_T({\frakC_{T}}{x_0}+{\frakD_{T}}{u}_{\rm opt}-y_\mathrm{ref})}.
\end{align*}
Note that the right-hand side corresponds to a backwards-in-time (due to the time flips) equation \changed{where the terminal value is given by the derivative of the terminal cost and the source term is the derivative of the integrand of the stage cost. More precisely, }
by considering the state $x_{\rm opt}$ and output $y_{\rm opt}$ of the system driven by the optimal control (see \eqref{eq:yopt}) and 
the adjoint state $\mu$ as in \eqref{eq:adjstate}, we obtain a~boundary value problem
\[\begin{aligned}
\spvek{\dot{x}_{\rm opt}(t)}{y_{\rm opt}(t)}&
= \sbvek{A\&B\\[-1mm]}{C\&D} \spvek{{x}_{\rm opt}(t)}{u_{\rm opt}(t)},& x_{\rm opt}(0)&=x_0,\\
\spvek{\dot{\mu}(t)}{u_{\rm opt}(t)}
&= \sbvek{{-[A\&B]^\mathrm{d}}}{{{-}[C\&D]^\mathrm{d}}}
\spvek{{\mu}(t)}{y_{\rm opt}(t)-y_{\rm ref}(t)},&\mu(T)&=F^*(Fx_{\rm opt}(T)-z_f).
\end{aligned}\]
\changed{Note that this boundary value problem is very similar to \eqref{eq:udsys2}. Here, however, besides the change of sign in the output equation, the output of the dual system node yields the optimal control $u_\mathrm{opt}$, while in \eqref{eq:udsys2} the output corresponds to the test functions in the optimality condition.}

\section{Constraints on the terminal state}\label{sec:termstate}

In this case we show how the incorporation of constraints on the terminal value is possible.
To this end, and under Assumptions~\ref{ass:ocp}, we consider the optimal control problem
\begin{align}\label{eq:optcons}
\begin{split}
\textrm{minimize}&\quad\frac{1}{2}\int_0^T\|y(t)-y_\mathrm{ref}(t) \|^2_{Y}{\rm d}t+ \frac{1}{2}\|Fx(T)-z_f\|^2_Z\\
\text{subject to}&\quad\spvek{\dot{x}(t)}{y(t)} = \sbvek{A\& B\\[-1mm] }{C\& D} \spvek{{x}(t)}{u(t)}, \quad
x(0) = x_0,\quad u \in \Uad,\quad F_cx(T)=z_c,
\end{split}
\end{align}
where $Z_c$ be a~Hilbert space,  $F_c\in L(X,Z_c)$ and $z_c\in Z$. 
In contrast to \eqref{eq:OCP}, the problem \eqref{eq:optcons} is provided with a~terminal constraint. Assuming that the set $\mathcal{B}_{F_c}$ as defined in \eqref{eq:dense1} is dense in $Z_c$, its mathematically precise interpretation is - analogously to the procedure for the terminal penalization - given by 
\begin{equation}\frakT_{F_c,T}\spvek{x_0}{u}=z_f.\label{eq:endconstr}
\end{equation}
where $\frakT_{F_c,T}$ is the unique extension of Definition~\ref{def:FIO}.
Clearly, to ensure non-emptyness of the set of admissible controls, we require the additional assumption that
there exists some admissible input $\hat{u}\in \Uad$ with the property as in Assumptions~\ref{ass:ocp}\,\eqref{ass:ocp6} and, additionally, $\frakT_{F_c,T}\spvek{x_0}{u}=z_f$. We further assume that the operator $u\mapsto\frakT_{F_c,T}\spvek{x_0}{u}$ (equivalently, $u\mapsto\frakT_{F_c,T}\spvek{0}{u}$) has closed range in $Z$, for reasons that will be explained later. Our supplementary assumptions for the optimal control problem with terminal state constraints are summarized below:
\begin{ass}\label{ass:ocpterm}
	$Z_c$ is a~Hilbert space, $F_c\in L(X,Z_c)$, $z_f\in Z_c$, and the following holds for the mappings in Table~\ref{tab:ourops}:
	\begin{enumerate}[(a)]
		\item\label{ass:ocpterma} \begin{equation}
		\mathcal{R}_{F_c}:=\setdef{\frakT_{F_c,T}\spvek{0}{u}}{u\in L^2([0,T];U)\text{ with }\spvek{0}{u}\in\dom(\frakT_{F_c,T})}
		\label{eq:RFc}\end{equation}is a~closed subspace of $Z_c$.
		\item\label{ass:ocptermb} For
		\begin{equation}{\widetilde\Uad}=\setdef{u\in\Uad}{\spvek{x_0}u\in\dom(\frakT_{F,T})\,\wedge\,\frakT_{F_c,T}\spvek{x_0}u=z_c}\label{eq:Uf}\end{equation}
		there exists an admissible input $\hat{u}\in\widetilde{\Uad}$ with 
		\[\spvek{x_0}{\hat{u}}\in \dom(\frakT_{F,T})\;\text{ and }\;
		{\frakC_{T}}{x_0}+{\frakD_T}{\hat{u}}-y_{\rm ref}\in L^2([0,T];Y).\]
	\end{enumerate}
\end{ass}
\changed{Loosely speaking, Assumption~\ref{ass:ocpterm}\eqref{ass:ocptermb} means that there exists a~state $x_c$ reachable from $x_0$ in time $T$ with, further, $Fx_c=z_f$. We now are in the position to formulate existence by the inclusion of the terminal constraint in the control constraint set $\widetilde{\Uad}$ as defined in \eqref{eq:Uf}.
	\begin{prop}
		Let Assumptions~\ref{ass:ocp} and~\ref{ass:ocpterm} hold. Then there is an optimal control
		$u_{\rm opt}\in\widetilde{\Uad}$ for \eqref{eq:optcons}
	\end{prop}
	\begin{proof}
		First, we note that \eqref{eq:optcons} is equivalent to $\inf_{u\in\widetilde{\Uad}}\mathcal{J}(u)$. Further, closedness and linearity of the $F$-input-to-state map yields that 
		the set
		\[\setdef{u\in L^2([0,T];U)}{\spvek{x_0}u\in\dom(\frakT_{F_c,T})\,\wedge\,\frakT_{F_c,T}\spvek{x_0}u=z_c}\]
		is closed and convex. Hence, the set $\widetilde{\Uad}$
		is the intersection of two sets which are both closed and convex and hence, it is closed and convex as well. Thus, Assumption~\ref{ass:ocp} is satisfied with $\widetilde{\Uad}$ instead of $\Uad$ such that Theorem~\ref{thm:optcont} is applicable to $\inf_{u\in\widetilde{\Uad}}\mathcal{J}(u)$.
\end{proof}}
\changed{We briefly comment on the extensions of the previous results in view of optimality conditions.}

To this end, we observe that we may directly apply Theorem~\ref{thm:optcond} to obtain that \eqref{eq:optder} holds for all $u \in \widetilde\Uad$ with $u-u_\mathrm{opt}\in \dom(\frakI_{F,T})$ and $\frakT_{F_c,T}\spvek{0}{u-u_\mathrm{opt}}=0$.
If, additionally $u-u_\mathrm{opt}\in H^2_{0l}([0,T];U)$, then Corollary~\ref{cor:aufloesen} can be applied. In particular, since $\frakT_{F_c,T}\spvek{0}{u-u_\mathrm{opt}}$ vanishes, we obtain, by using Lemma~\ref{lem:adj}, that
\[\forall\,\lambda\in Z:\langle {\frakC_{\mathrm{d},T}}F_c^*\lambda,u-u_\mathrm{opt}\rangle_{H_{0r}^{-2}([0,T];U),H_{0l}^{2}([0,T];U)}=0,\]
and the optimality condition in Corollary~\ref{cor:aufloesen} means that we additionally have
\begin{align}
\forall\, u \in \widetilde\Uad\ \mathrm{ with }\, u-u_\mathrm{opt}\in H^2_{0l}([0,T];U): 
\re\left\langle u-u_{\rm opt}, u_\mathrm{d}\right\rangle_{H^2_{0l}([0,T];U),H^{-2}_{0r}([0,T];U)} \geq 0\label{eq:optendpoint}
\end{align}
for all $u_\mathrm{d}\in H^{-2}_{0r}([0,T];U)$ with
\[
u_\mathrm{d} = \Refl_T[{\frakC_{\mathrm{d},T}},\,{\frakD_{\mathrm{d},T}}]  \pvek{F^*\frakT_{F,T}\spvek{x_0}{u_{\rm opt}}-z_f+F_c^*\lambda}{\Refl_T\left({\frakC_T}{x_0}+{\frakD_T}{u}_{\rm opt}-y_{\rm ref}\right)},\quad \lambda\in Z_c.
\]
In the case where the closure of $\mathcal{T}$ as defined in \eqref{eq:T} \changed{(with $\Uad$ replaced by $\widetilde{\Uad}$)} possesses a nonempty relative interior in the closed space
\[\mathcal{N}_{F_c}:=\setdef{u\in L^2([0,T];U)}{\frakT_{F_c,T}\spvek{0}{u}=0},\]
then, \changed{in analogy to Proposition~\ref{prop:aufloesen}}, the inequality \eqref{eq:optendpoint} becomes an equation. 

\changed{We will now deduce the corresponding optimality condition in terms of an optimality boundary value problem including an adjoint state as an analogon to \eqref{eq:udsys0}}. To this end, note that for all $\delta u\in \dom(\frakI_{F,T})\cap \mathcal{N}_{F_c}$, we have that\[\left\langle\pvek{\frakT_{F,T}\spvek{x_0}{u_{\rm opt}}-z_f}{{\frakC_T}{x_0}+{\frakD_T}{u}_{\rm opt}-y_{\rm ref}},\frakI_{F,T}\delta u\right\rangle_{X\times L^2([0,T];Y)} = 0.\]
This shows that for all $\delta u\in \mathcal{N}_{F_c}\cap H^2_{0l}([0,T];U)$, 
\[\left\langle[{\frakC_{\mathrm{d},T}}F^*,\,{\frakD_{\mathrm{d},T}}]\pvek{\frakT_{F,T}\spvek{x_0}{u_{\rm opt}}-z_f}{{\frakC_T}{x_0}+{\frakD_T}{u}_{\rm opt}-y_{\rm ref}},\delta u\right\rangle_{H^{-2}_{0r}([0,T];U),H^2_{0l}([0,T];U)} = 0.\]
Now incorporating that $\mathcal{R}_{F_c}$ as in \eqref{eq:RFc} is a~closed subspace of $Z_c$, this holds for the adjoint of the closed and densely defined mapping $u\mapsto\frakT_{F_c,T}^*$ as well. The latter is, by applying Proposition~\ref{prop:adjFinput} with $Y=\{0\}$, given by the mapping $\lambda\mapsto \Refl_T{\frakC_{\mathrm{d},T}}F_c^*\lambda$ defined on the space of all $\lambda\in Z$  and we obtain that there exists some $\lambda\in Z_c$ with ${\frakC_{\mathrm{d},T}}F_c^*\lambda\in L^2([0,T];U)$. Consequently, there exists some $\lambda\in Z$, such that
\[\left[{\frakC_{\mathrm{d},T}}F^*,\,{\frakD_{\mathrm{d},T}}\right]\pvek{\frakT_{F,T}\spvek{x_0}{u_{\rm opt}}-z_f}{\Refl_T\big({\frakC_T}{x_0}+{\frakD_T}{u}_{\rm opt}-y_{\rm ref}\big)}+{\frakC_{\mathrm{d},T}}F_c^*\lambda=0.\]
By invoking that, for $G=\diag(F^*,F_c^*)$
\[[{\frakC_{\mathrm{d},T}}F^*,\,{\frakC_{\mathrm{d},T}}F_c^*]={\frakC_{\mathrm{d},T}}G,\]
we obtain that there exists some $\lambda\in Z$ such that 
\[\left(\begin{smallmatrix}
{\frakT_{F,T}\spvek{x_0}{u_{\rm opt}}-z_f}\\\Refl_T\big({\frakC_T}{x_0}+{\frakD_T}{u}_{\rm opt}-y_{\rm ref}\big)
\end{smallmatrix}\right)\in \dom\frakO_{\mathrm{d},G,T}\]
with 
\[\frakO_{\mathrm{d},G,T}\left(\begin{smallmatrix}\lambda\\
{\frakT_{F,T}\spvek{x_0}{u_{\rm opt}}}\\\Refl_T\big({\frakC_T}{x_0}+{\frakD_T}{u}_{\rm opt}-y_{\rm ref}\big)
\end{smallmatrix}\right)=0.\]
This corresponds to the solution of the system
\[
\begin{aligned}
\spvek{\dot{x}_{\rm opt}(t)}{y(t)}&
= \sbvek{A\&B\\[-1mm]}{C\&D} \spvek{{x}_{\rm opt}(t)}{u_{\rm opt}(t)},& x_{\rm opt}(0)&=x_0,\quad F_cx_{\rm opt}(T)=z_c,\\
\spvek{\dot{\mu}(t)}{0}
&= \sbvek{{-[A\&B]^\mathrm{d}}}{{[C\&D]^\mathrm{d}}}
\spvek{{\mu}(t)}{y_{\rm opt}(t)-y_{\rm ref}(t)},& \mu(T)&=F^*(Fx_{\rm opt}(T)-z_f)+F_c^*\lambda
\end{aligned}\]
for the unknowns $x_{\rm opt},u_{\rm opt},y,\mu$ and $\lambda$.

Let us finally give some remarks on the case where the full terminal state is prescribed.
\begin{rem}\label{rem:fullterm}
	Consider an optimal control problem with full end point constraint, i.e.,
	\begin{align}\label{eq:ocpremark}
	\begin{split}
	\textrm{minimize}&\quad\frac{1}{2}\int_0^T\|y(t)-y_\mathrm{ref}(t) \|^2_{Y}\, {\rm d}t\\
	\text{subject to}&\quad\spvek{\dot{x}(t)}{y(t)} = \sbvek{A\& B\\[-1mm] }{C\& D} \spvek{{x}(t)}{u(t)}, \quad
	x(0) = x_0,\quad u \in \Uad,\quad x(T)=x_f\in X.
	\end{split}
	\end{align}
	It is clear that the penalization of the terminal state is obsolete in this case. We adopt Assumptions~\ref{ass:ocp} and, furthermore, Assumptions~\ref{ass:ocpterm}~\eqref{ass:ocpterma}.
	To ensure the validity of Assumption~\ref{ass:ocpterm}\eqref{ass:ocpterma}, we must further assume that the space $\im\left(\frakB_{T}\right)\cap X$ \changed{(which is \eqref{eq:RFc} for $F_c=I$)}
	is a closed subspace of $X$. This indeed imposes a restriction, as it precludes the consideration of various (in particular parabolic) problems. To extend the applicability of our theory to optimal control problems with hard terminal state constraints and non-closed reachability space, additional efforts, such as selecting alternative norms in the state space, are required. These endeavors are beyond the scope of this article.
	
	However, if the reachability space is closed, and, additionally, the closure of $\mathcal{T}$ as defined in \eqref{eq:T} possesses a nonempty relative interior in $\ker\frakB_{I,T}$, the optimal control problem \changed{\eqref{eq:ocpremark}} leads to the solution of the~system
	\[
	\begin{aligned}
	\spvek{\dot{x}_{\rm opt}(t)}{y_{\rm opt}(t)}&
	= \sbvek{A\&B\\[-1mm]}{C\&D} \spvek{{x}_{\rm opt}(t)}{u_{\rm opt}(t)},&x_{\rm opt}(0)=x_0,\ x_{\rm opt}(T)=x_f,\\
	\spvek{\dot{\mu}(t)}{0}
	&= \sbvek{{-[A\&B]^\mathrm{d}}}{{\phantom{-}[C\&D]^\mathrm{d}}}
	\spvek{{\mu}(t)}{y_{\rm opt}(t)-y_{\rm ref}(t)},
	\end{aligned}\]


\end{rem}

\section{Control of port-Hamiltonian systems with minimal energy supply}

As a specific application of the theory presented thus far, we examine state transition of port-Hamiltonian systems while supplying the system with the minimal amount of (physical) energy. This optimal control problem has been suggested in \cite{SPFWM23,FMPSW22} for finite-dimensional port-Hamiltonian systems. In \cite{PSTFMW21}, an extension to infinite-dimensional systems was given, where, however, the considered system class includes a~bounded input operator with finite-dimensional control space and a governing main operator which can be split into a dissipative and a skew-adjoint part. 

We commence with a brief introduction to port-Hamiltonian system nodes, as developed in \cite{PhilReis23}. For sake of brevity, our setup is slightly simpler than that presented in \cite{PhilReis23}. The class nevertheless covers a~wide range of physical examples, such as Maxwell's equations, advection-diffusion equations and linear hyperbolic systems in one spatial variable, a~wave equation in more spatial variables \cite{Farkas2023}, and Oseen's equations \cite{ReisSchal23}.

\begin{defn}[Port-Hamiltonian system]\label{def:pHsys}
	Let $X$, $U$ be Hilbert spaces.
	Let
	\[M=\sbvek{\tF\&\tG}{\tK\&\tL}:X^*\times U\supset\dom(M)\to X\times U^*\]
	be a~{\em dissipation node} on $(X,U)$, that is, 
	\changed{$M$ is dissipative, and}
	\begin{enumerate}[(a)]
		\item\label{def:dissnode2} $\tF\&\tG:X^*\times U\supset\dom(\tF\&\tG)\to X$ with $\dom(\tF\&\tG)=\dom(M)$ is closed;
		\item\label{def:dissnode1} $\tK\&\tL\in L(\dom(\tF\&\tG),U^*)$;
		\item\label{def:dissnode3} for all $u\in U$, there exists some $x'\in X^*$ with $\spvek{x'}{u}\in \dom(M)$;
		\item\label{def:dissnode4} for the {\em main operator} $\tF : X\supset\dom(\tF)\to X$ defined by \[\dom(\tF) := \setdef{x'\in X^*}{\spvek {x'}0\in\dom(M)}\] and $Fx' := \tF\&\tG\spvek {x'}0$, there exists some $\lambda>0$ such that $\lambda R_X^{-1}-\tF$ has dense range, where $R_X$ is the Riesz isomorphism on $X$.
	\end{enumerate}
	Further, let $H\in L(X,X^*)$ be positive, self-dual and bijective.
	Then we call
	\begin{equation}
	\spvek{\dot{x}(t)}{y(t)}
	= \sbvek{\tF\&\tG\\[-1mm]}{-\tK\&\tL} \spvek{{H}{x}(t)}{u(t)}.\label{eq:pHODEnode}\end{equation}
	a~{\em port-Hamiltonian system on $(X,U)$}.
\end{defn}

It follows from \cite[Prop.~3.6]{PhilReis23} that $\sbvek{\tF\&\tG}{-\tK\&\tL}\sbmat{H}00{\Id_{U}}$ defines a~system node, whence we call it a~{\em port-Hamiltonian system node}. Since $H$ is assumed to be bounded, self-dual, positive and boundedly invertible, the mapping 
\begin{equation}
x\mapsto \|x\|_{H}:=\langle x(t),Hx(t)\rangle_{X,X^*}^{1/2}\label{eq:Hnorm}
\end{equation} is equivalent to the original norm $\|\cdot\|_X$ on $X$. It can be moreover seen that $\tF H$ is a~maximally dissipative (see \cite{EngeNage00} for a~definition) on $X$ endowed with the norm $\|\cdot\|_{H}$ (and thus with the scalar product~$\langle \cdot, H\cdot\rangle$), whence, by \cite[Chap.~2, Sec.~3b]{EngeNage00},
it generates a~contractive semigroup on $X$ equipped with $\|\cdot\|_{H}$.
Further, by \cite[Prop.~3.11]{PhilReis23}, for all $T>0$, the generalized trajectories $(y,x,u)\,\in\,L^2([0,T];U^*)\times C([0,T];X)\times  L^2([0,T];U)$ (and thus also the classical \changed{trajectories}) fulfill the {\em dissipation inequality}
\begin{multline}
\tfrac12 \langle x(T),Hx(T)\rangle_{X,X^*} -\tfrac12 \langle x(0),Hx(0)\rangle_{X,X^*}\\=
\int_0^T \re\left\langle  M\spvek{{H}x(t)}{u(t)},\spvek {{H}x(t)}{u(t)}\right\rangle_{X\times U^*,X\times U}{\rm d}t+\int_0^T \re\langle {y(t)},{u(t)}\rangle_{U^*,U}{\rm d}t\\
\leq\int_0^t \re\langle {y(t)},{u(t)}\rangle_{U^*,U}{\rm d}t.
\label{eq:enbalinf}
\end{multline}
The physical interpretation of the latter is that $\tfrac12 \langle x(t),Hx(t)\rangle_{X,X^*}$ represents the stored energy at time $t\in[0,T]$, while $\re\langle {y(t)},{u(t)}\rangle_{U^*,U}$ signifies the power supplied to the system. Consequently, the term
\[-\re\left\langle  M\spvek{{H}x(t)}{u(t)},\spvek {{H}x(t)}{u(t)}\right\rangle_{X\times U^*,X\times U} \ge 0\]
stands for the power dissipated by the system at time $t$. Let us first consider the problem of minimizing the the supplied energy while controlling the system from a~given initial value $x_0\in X$ to a~prescribed terminal state $z_f\in X$, that is,
\begin{align}\label{eq:energymin}
\begin{split}
\textrm{minimize}&\quad \int_0^T\re\langle {y(t)},{u(t)}\rangle_{U^*,U}\, {\rm d}t\\
\text{subject to}&\quad\spvek{\dot{x}(t)}{y(t)}
= \sbvek{\phantom{-}\tF\&\tG\\[-1mm]}{-\tK\&\tL} \spvek{{H}{x}(t)}{u(t)}, \quad
x(0) = x_0\in X,\quad u \in \Uad,\quad x(T)=x_c\in X.
\end{split}
\end{align}
The optimal control problem discussed here falls - at hand - outside the scope of the class addressed in the previous section. Nevertheless, drawing inspiration from a technique highlighted in the finite-dimensional scenario in \cite{SPFWM23}, we propose a reformulation of the aforementioned optimal control problem, which enables a direct application of the theory presented earlier. Specifically, we create an artificial output, the squared norm of which serves as a representation of dissipated power. The dissipation inequality leads to the observation that minimizing the norm of this artificial output yields the solution to the energy minimization problem~\eqref{eq:energymin}. The construction of this artificial output is subject of the following result.
\begin{prop}\label{prop:disseq}
	Let
	\[M=\sbvek{\tF\&\tG\\[-1mm] }{\tK\&\tL}:X^*\times U\supset\dom(M)\to X\times U^*\]
	be a~dissipation node on $(X,U)$. Then there exists a~Hilbert space $W$ and $\tR\&\tS\in L(\dom(\tF\&\tG),W)$, such that
	\begin{equation}
	\forall\,\spvek{x'}{u}\in\dom(\tF\&\tG):\quad
	-\re\left\langle  M\spvek{x'}{u},\spvek {x'}{u}\right\rangle_{X\times U^*,X\times U}=2\left\|\tR\&\tS\spvek {x'}{u}\right\|_W^2.
	\label{eq:disseq}\end{equation}
\end{prop}
\begin{proof}
	Consider the mapping
	\[\begin{aligned}
	\mathcal{Q}:\quad \dom(\tF\&\tG)&\to\R,\\
	\spvek{x'}{u}&\mapsto-2\re\left\langle  M\spvek{x'}{u},\spvek {x'}{u}\right\rangle_{X\times U^*,X\times U}.
	\end{aligned}
	\]
	Then $\mathcal{Q}$ is a~bounded symmetric nonnegative quadratic form in the sense of \cite[Chap.~6]{Kato80}.
	By using
	Kato's First Representation Theorem \cite[Sec.\ VI.2, Thm.\ 2.1]{Kato80}, there exists some self-adjoint and nonnegative  $Q\in L(\dom(\tF\&\tG))$, such that
	\[\forall\,\spvek{x'}{u}\in\dom(\tF\&\tG):\quad\mathcal{Q}\big(\spvek{x'}{u}\big)=\left\langle \spvek{x'}{u},Q\spvek{x'}{u}\right\rangle_{\dom(\tF\&\tG)}.\]
	Now we set
	\[W:=(\ker Q)^\bot,\quad \tR\&\tS\spvek{x'}{u}:=\sqrt{Q}\spvek{x'}{u}\quad \forall \, \spvek{x'}{u}\in\dom(\tF\&\tG),\]
	where the orthogonal complement has to be understood with respect to the inner product in $\dom(\tF\&\tG)$.
	This gives $\tR\&\tS\in L(\dom(\tF\&\tG),W)$, and \eqref{eq:disseq} holds.
\end{proof}
It can be seen that $W$ and $\tR\&\tS$ can be chosen in a~way that $\im \left(\tR\&\tS\right)$ is dense in $W$. In this case, $W$ is uniquely determined up to \changed{isometric isomorphy}, and $\tR\&\tS$ is uniquely determined up to the multiplication from the left with an isometric isomorphism.\\
One can readily observe that, for $\tR\&\tS$ as in Proposition~\ref{prop:disseq},
\[
\left[\begin{smallmatrix}\phantom{-}\tF\&\tG\\[-1mm] \\-\tK\&\tL\\[-1mm] \\\phantom{-}\tR\&\tS\end{smallmatrix}\right]
\sbmat{H}00{\Id_{U}}\]
is a~system node on $(X,U,U^*\times Z)$. Further, the dissipation inequality \eqref{eq:enbalinf} gives rise to the property that the generalized trajectories $(x,u,y)$ of the port-Hamiltonian system \eqref{eq:pHODEnode} on $[0,T]$ fulfil $\tR\&\tS\spvek{Hx}{u}\in L^2([0,T];W)$, and the 
classical (and thus also the generalized) trajectories $(x,u,\spvek{y}{w})$ of the system
\begin{equation}\left(\begin{smallmatrix}\dot{x}(t)\\y(t)\\w(t)\end{smallmatrix}\right)=\left[\begin{smallmatrix}\phantom{-}\tF\&\tG\\[-1mm]\\-\tK\&\tL\\[-1mm]\\\phantom{-}\tR\&\tS\end{smallmatrix}\right]\spvek{Hx(t)}{u(t)}\label{eq:extsys}
\end{equation}
on $[0,T]$ satisfy
\begin{align}
\int_0^T \re\langle {y(t)},{u(t)}\rangle_{U^*,U}{\rm d}t
=
\frac12 \langle x(T),Hx(T)\rangle_{X,X^*} -\frac12 \langle x(0),Hx(0)\rangle_{X,X^*}+\changed{2} 	\int_0^T \|w(t)\|^2_{W}{\rm d}t,
\label{eq:enbalinfZ}
\end{align}
which, in passing, shows that $\tR H$ is an admissible output operator for the semigroup generated by $\tF H$. We will call
\begin{equation}w(t)=\tR\&\tS\spvek{Hx(t)}{u(t)}\label{eq:dissout}\end{equation}
a~{\em dissipation output} for the port-Hamiltonian system \eqref{eq:pHODEnode}.\\
Since the initial and terminal state are prescribed, the minimization of the supplied energy, by using \eqref{eq:enbalinfZ}, corresponds to the minimization of the dissipated energy. Consequently, we may reformulate the \changed{optimal control problem~\eqref{eq:energymin}} to
\[\begin{aligned}
\textrm{minimize}&\quad \frac12 \int_0^T\|w(t)\|_{W}^2\, {\rm d}t\\
\text{subject to}&\quad\spvek{\dot{x}(t)}{w(t)}
= \sbvek{\tF\&\tG\\[-1mm]}{\tR\&\tS} \spvek{{H}{x}(t)}{u(t)}, \quad
x(0) = x_0\in X,\quad u \in \Uad,\quad x(T)=x_c\in X,\nonumber
\end{aligned}\]
which is now truly belonging to the class treated in Section~\ref{sec:termstate}. Recall from Remark~\ref{rem:fullterm} that the reachability space \changed{as defined in \eqref{eq:RFc} has to be closed in view of Assumption~\ref{ass:ocpterm} to deduce optimality conditions as done in Section~\ref{sec:termstate}}.

The construction of the dissipation output 
also allows for appropriate reformulations of more general energy-efficient (optimal) control problems, leading to ones discussed in earlier sections.
\begin{rem}\label{eq:modpHcont}
	Let a~port-Hamiltonian system \eqref{eq:pHODEnode} be given, and let $\tR\&\tS\in L(\dom(\tF\&\tG),W)$ be with the properties as in Proposition~\ref{prop:disseq}, i.e., $w$ with \eqref{eq:dissout} is a~dissipation output for \eqref{eq:pHODEnode}.  Rather than seeking for optimal controls that ensure a~perfect landing at a predetermined terminal state, we now only aim for a partial terminal condition on the state (or even none at all, which can be achieved by choosing $Z=\{0\}$). Our objective is to penalize both energy wastage and a weighted deviation from a desired terminal state $x_f$. That is, we consider the optimal control problem
	
	\begin{align}\label{eq:energymin_terminalweight}
	\begin{split}
	\textrm{minimize}&\quad
	\int_0^T\re\langle {y(t)},{u(t)}\rangle_{U^*,U}\, {\rm d}t+\frac12\|Fx(T)-z_f\|^2_Z\\
	\text{subject to}&\quad\spvek{\dot{x}(t)}{y(t)} = \sbvek{\phantom{-}\tF\&\tG\\[-1mm] }{-\tK\&\tL} \spvek{H{x}(t)}{u(t)}, \quad
	x(0) = x_0,\quad u \in \Uad,\quad F_cx(T)=z_c.
	\end{split}
	\end{align}
	
	Now using \eqref{eq:enbalinfZ}, we obtain that this leads to the minimization of 
	\[    \frac12 \langle x(T),Hx(T)\rangle_{X,X^*} +	\frac12\int_0^T \|w(t)\|^2_{W}{\rm d}t+\frac12\|Fx(T)-z_f\|^2_Z.
	\]
	To lead this back to a~problem treated in earlier sections, we 
	first observe that, for the Riesz isomorphism $R_X:X\to X^*$, we have that $R_XH$ is self-adjoint and positive, and thus it possesses an operator square root $\sqrt{R_XH}:X\to X$. 
	
	Now the \changed{optimal control problem~\eqref{eq:energymin_terminalweight}} can be reformulated to
	\[
	\begin{aligned}
	\textrm{minimize}&\quad
	\frac12\int_0^T \|w(t)\|_W^2\, {\rm d}t+\frac12\left\|\sbvek{F\\[-1mm]}{\sqrt{R_XH}}x(T)-\spvek{z_f\\[-1mm]}0\right\|_{Z\times X}^2\\
	\text{subject to}&\quad\spvek{\dot{x}(t)}{w(t)} = \sbvek{\tF\&\tG\\[-1mm] }{\tR\& \tS} \spvek{H{x}(t)}{u(t)}, \quad
	x(0) = x_0,\quad u \in \Uad,\quad F_cx(T)=z_c,\nonumber
	\end{aligned}
	\]
	Last, we examine yet another category of energy-optimal control problems. For this purpose, let us consider a system 
	\[\left(\begin{smallmatrix}
	{\dot{x}(t)}\\{y(t)}\\v(t)
	\end{smallmatrix}\right)
	= \left[\begin{smallmatrix}
	\phantom{-}\tF\&\tG\\[-1mm] \\{-\tK\&\tL}\\[-1mm] \\\phantom{-}\tC\&\tD
	\end{smallmatrix}\right] \spvek{{H}{x}(t)}{u(t)},\]
	such that $\sbvek{\tF\&\tG\\[-1mm]}{-\tK\&\tL}\sbmat{H}00{\Id_{U}}$ is a~port-Hamiltonian system node, and $\tC\&\tD\in L(\dom(\tF\&\tG),V)$ for some Hilbert space $V$ and let $v_{\rm ref}\in L^2([0,T],V)$. Now, let us delve into the optimal control problem that seeks for a~compromise between energy conservation and tracking of $v_{\rm ref}$ by $v$, i.e.,
	\[
	\begin{aligned}
	\textrm{minimize}&\quad
	\int_0^T\re\langle {y(t)},{u(t)}\rangle_{U^*,U}\,+ \frac12\|v(t)-v_\mathrm{ref}(t) \|^2_{V}{\rm d}t+\frac12\|Fx(T)-z_f\|^2_Z\\
	\text{subject to}&\quad\left(\begin{smallmatrix}\dot{x}(t)\\y(t)\\v(t)\end{smallmatrix}\right)=\left[\begin{smallmatrix}\phantom{-}\tF\&\tG\\[-1mm]\\-\tK\&\tL\\[-1mm]\\\phantom{-}\tC\&\tD\end{smallmatrix}\right]\spvek{Hx(t)}{u(t)}, \quad
	x(0) = x_0,\quad u \in \Uad,\quad F_cx(T)=z_c.\nonumber
	\end{aligned}
	\]
	This leads, \changed{using the energy balance~\eqref{eq:enbalinfZ}}, to an optimal control problem
	\begin{align}\label{eq:ocpterminalenergy}
	\begin{split}
	\textrm{minimize}&\quad
	\frac12\int_0^T\left\|\spvek{w(t)}{v(t)}-\spvek0{v_\mathrm{ref}(t)}\right\|_{W\times V}^2\,{\rm d}t+\frac12\left\|\sbvek{F\\[-1mm]}{\sqrt{R_XH}}x(T)-\spvek{z_f\\[-1mm]}0\right\|_{Z\times X}^2\\
	\text{subject to}&\quad\left(\begin{smallmatrix}\dot{x}(t)\\w(t)\\v(t)\end{smallmatrix}\right)=\left[\begin{smallmatrix}\tF\&\tG\\[-1mm]\\\tR\&\tS\\[-1mm]\\\tC\&\tD\end{smallmatrix}\right]\spvek{Hx(t)}{u(t)}, \quad
	x(0) = x_0,\quad u \in \Uad,\quad F_cx(T)=z_c.
	\end{split}
	\end{align}
	
	Note that, by providing $X$ with the norm $\|\cdot\|_H$ as in \eqref{eq:Hnorm} and identifying the Hilbert spaces $W$, $Z$ with their respective anti-duals, the adjoint of the system node in the \changed{optimal control problem \eqref{eq:ocpterminalenergy}} fulfils
	\[    \left(\left[\begin{smallmatrix}\tF\&\tG\\[-1mm] \\\tR\&\tS\\[-1mm] \\\tC\&\tD\end{smallmatrix}\right]
	\sbmat{H}00{\Id_{U}}\right)^*=
	\left[\begin{smallmatrix}I_X&0\\0&R_{U}^{-1}\end{smallmatrix}\right]\left[\begin{smallmatrix}\tF\&\tG\\[-1mm] \\\tR\&\tS\\[-1mm] \\\tC\&\tD\end{smallmatrix}\right]^*
	\left[\begin{smallmatrix}{H}&0&0\\0&{\Id_{W}}&0\\0&0&I_V\end{smallmatrix}\right],\]
	where the second factor on the right hand side stands for the dual of 
	\[\left[\begin{smallmatrix}\tF\&\tG\\[-1mm] \\\tR\&\tS\\[-1mm] \\\tC\&\tD\end{smallmatrix}\right],\]
	and $R_U:U\to U^*$ is the Riesz map.
	In particular, its domain is a~dense subspace of $X^*\times W\times V$, and it maps to $X\times U^*$. We would like to draw the reader's attention, however, to the fact that the adjoints of $F$ and $F_c$ with respect to the norm $\|\cdot\|_H$ must be considered in this context.
\end{rem}

\section{Applications}
In this part, we illustrate the results through two examples involving boundary control. The first example considers a advection-diffusion-reaction equation with Dirichlet control and Neumann observation, while the second one involves a~wave equation on an L-shaped domain.

Indeed, both of the addressed problems are real, not complex. While our general theory is primarily formulated for complex problems, it is equally applicable to real scenarios, as highlighted in Remark~\ref{rem:realproblems}. Specifically, it is important to emphasize that all the spaces involved are now considered to be real.

\subsection{Advection-diffusion-reaction equation with Dirichlet control and Neumann observation}
Consider the open unit interval $[0,1]\subset \R$, a heat diffusion coefficient $a:[0,1]\to \R$, a convection field $b:[0,1]\to \R$ and a reaction field $c:[0,1]\to \R$. We inspect the advection-diffusion-reaction equation 
\begin{subequations}\label{eq:heateq}
	\begin{align}
	\tfrac{\partial}{\partial t} x(t,\xi)&= \tfrac{\partial}{\partial \xi} \left(a(\xi) \tfrac{\partial}{\partial \xi}x(t,\xi)\right) + b(\xi)\tfrac{\partial}{\partial \xi} x(t,\xi) + c(\xi) x(t,\xi),\qquad t\geq 0,\ \xi \in [0,1],
	\intertext{with single input and output given by}
	u(t) &= x(t,1), \qquad 
	y(t) = a(0) \tfrac{\partial}{\partial \xi} x(t,0),\qquad
	t\geq 0.
	\end{align}
\end{subequations}

We first gather our manageable prerequisites on the coefficients.
\begin{ass}\label{ass:heat}\
	\begin{enumerate}[(a)]
		\item $a\in L^\infty([0,1])$ is positive almost everywhere, and $a^{-1}\in L^\infty([0,1])$;
		\item $b,c\in L^\infty([0,1])$.
	\end{enumerate}
\end{ass}
We note that the \changed{system~\eqref{eq:heateq}} is not well-posed, cf.\ \cite{Schw20}, such that in optimal control, this setup is a delicate problem. We refer to \cite[Section 9]{Lions71}, where a Dirichlet optimal control problem in higher dimension with square integrable controls and distributional output (in fact $H^{-1}$-valued Neumann observation) was analysed. The stationary problem was thoroughly analyzed in \cite{KuniVexl2007}.

We choose the state space $X = L^2([0,1])$, and, since our system is single-input-single-output, we have $U=Y=\R$. 

By denoting the spatial derivative by a~prime, and further setting, in coherence to the spaces introduced in the first section,
\[H^1_{0l}([0,1])=
\setdef{v\in H^1([0,T])}{v(0)= 0}\]
which due to its closedness is a~Hilbert space equipped when with the standard inner product in $H^1([0,1])$ \changed{ and the induced norm, that is,
	\begin{align*}
	\langle y,v\rangle_{H^1([0,1])} = \langle y',v'\rangle_{L^2([0,1])} + \langle y,v\rangle_{L^2([0,1])}
	\end{align*}
	for all $y,v\in H^1([0,1])$ and $\|y \|^2_{H^1([0,1])} = \langle y,y\rangle_{H^1([0,1])}$ for all $y\in H^1([0,1])$}. We note that the evaluation of $x$ at one represents a~bounded operator due to the continuous embedding $H^1([0,1])\hookrightarrow C([0,1])$, cf.~\cite{adams2003sobolev}. The system node corresponding to the advection-diffusion-reaction equation is defined by $S=\sbvek{A\&B}{C\&D}$ with
\begin{subequations}\label{eq:advdiff_dissnode}
	\begin{align}
	\dom{A\&B}:=\setdef{\spvek{x}{u}\in H^1_{0l}([0,1])\times \R}{\,(a x')'\in
		L^2([0,1])\;\wedge\;x(1)=u}
	\end{align}
	and
	\begin{align}
	\forall \, \spvek{x}{u}\in\dom(A\&B):\quad A\&B\spvek{x}{u}&= \big(a x'\big)'+b x' + cx,\;\; C\&D\spvek{x}{u}= {(ax')(0)}.
	\end{align}
\end{subequations}
The evaluation of $a x'$ is well-defined since $a x'\in H^{1}([0,1])$.

We first verify that this indeed defines a system node in the sense of Definition~\ref{def:sysnode}. In \cite[Sec.~4.1]{PhilReis23}, a comparable equation was explored within a port-Hamiltonian framework. However, in that scenario, the advection field was divergence-free, the reaction field vanished, and the complete Dirichlet and Neumann traces were selected as the input and output, respectively. Consequently, both the main operator and the system node itself were dissipative operators. Though dissipativity streamlined certain steps, we are now able to present a~shorter proof \changed{of \eqref{eq:advdiff_dissnode} defining a system node} for a~more general scenario.

To this end, we consider the ``bilinear form associated to $A\&B$''
\begin{equation}\begin{aligned}
q:\quad \changed{H^1_{0}([0,1])\times H^1_{0}([0,1])}&\to \R\\
(v,w) &\mapsto -\langle av',w'\rangle_{L^2([0,1])} + \langle bv'+cv,w\rangle_{L^2([0,1])}.
\end{aligned}\label{eq:qform}\end{equation}
We see that $q$ is continuous in the sense that
\begin{equation}\exists\,c>0:\quad |q(v,w)|\leq c \|v\|_{H^1([0,1])}\|w\|_{H^1([0,1])}\label{eq:qcont}
\end{equation}
\\
A crucial component for proving that the aforementioned operators constitute a system node is that the \changed{negative of the form} $q$ ensures a specific, albeit weaker, form of coercivity.
\begin{lem}\label{lem:weakcoer}
	Under Assumptions~\ref{ass:heat}, there exists some $\mu\in \R$ and some $\alpha > 0$ such that the form $q$ as in \eqref{eq:qform} fulfils
	\begin{align*}
	\forall\,v\in \changed{H^1_{0}([0,1])}:\qquad       q(v,v)   \leq -\alpha \|v\|^2_{H^1([0,1])}+\mu \|v\|^2_{L^2([0,1])}.
	\end{align*}
\end{lem}
\begin{proof}
	By using that $a$ is bounded from below by $\underline{a}:=\|a^{-1}\|_{L^\infty([0,1])}^{-1}>0$, we obtain that for all \changed{$v\in H^1_{0}([0,1])$},
	\begin{align*}
	q(v,v) &= -\langle av',v'\rangle_{L^2([0,1])} + \langle bv'+cv,v\rangle_{L^2([0,1])}\\
	&\leq -\underline{a}\|v'\|_{L^2([0,1])}^2 + \|b\|_{L^\infty(0,1)} \|v'\|_{L^2([0,1])}\|v\|_{L^2([0,1])} + \|c\|_{L^\infty([0,1])}\|v\|^2_{L^2([0,1])},
	\end{align*}
	and    Young's inequality implies $$\|b\|_{L^\infty(0,1)} \|v'\|_{L^2([0,1])}\|v\|_{L^2([0,1])} \leq \tfrac12\left(\underline{a}^{-1}{\|b\|^2_{L^\infty([0,1])}\|v\|^2_{L^2([0,1])}} + \underline{a}\|v'\|_{L^2([0,1])}^2\right)$$
	such that 
	\begin{align*}
	q(v,v) \leq -\tfrac{\underline{a}}{2}\,\|v'\|_{L^2([0,1])}^2 + \left(\|c\|_{L^\infty([0,1])} +  \tfrac{1}{2\underline{a}}\,\|b\|^2_{L^\infty([0,1])}\right)\|v\|^2_{L^2([0,1])}.
	\end{align*}
	Hence, choosing, e.g.,
	\begin{align*}
	\mu = \tfrac{\underline{a}}{2} + \left(\|c\|_{L^\infty([0,1])} +  \tfrac{1}{2\underline{a}}\, \|b\|^2_{L^\infty(0,1)}\right)
	\end{align*}
	yields the claim for $\alpha = \frac{\underline{a}}{2}$.
\end{proof}


\changed{Lemma~\ref{lem:weakcoer}} is the main ingredient for showing the system node properties of the operator composed by $A\&B$ and $C\&D$ \changed{as defined in \eqref{eq:advdiff_dissnode}}. We additionally \changed{show} that $S$ is bijective which will be used later on.
\begin{prop}\label{prop:heatsysnode}
	Under Assumptions~\ref{ass:heat}, the operator $S=\sbvek{A\&B}{C\&D}$ with $A\&B$ and $C\&D$ as in \eqref{eq:advdiff_dissnode} is a~system node. Further, $S:\dom(S)\to L^2([0,1])\times\R$ is bijective. 
\end{prop}
\begin{proof}
	We have to successively verify that $S$ has the properties as indicated in Definition~\ref{def:sysnode}.\\
	\noindent \textbf{\ref{def:sysnoded}):}
	We show that $A$ is the generator of a~strongly continuous semigroup: By definition of the weak derivative, we have 
	$Ax=z$ for $x\in \dom(A)$, $z\in L^2([0,1])$ if, and only if, $x\in H^1_0([0,1])$, and $q$ as in \eqref{eq:qform} fulfills
	\[\forall\, \varphi\in H^1_0([0,1]):\qquad\langle z,\varphi\rangle_{L^2([0,1])}=q(x,\varphi).\]
	Now using that $-q$ fulfills the weaker kind of coercivity in the sense of Lemma~\ref{lem:weakcoer},
	$A$ generates a strongly continuous (in fact, even analytic) semigroup due to \cite[Theorem 4.2]{ArenElst12}. 
	
	\noindent\textbf{\ref{def:sysnodec}):} We show that for all $u\in \R$, there exists some $x\in L^2([0,1])$ such that $\spvek{x}{u}\in \dom(S) = \dom(A\&B)$.
	
	First, with $\alpha$ and $\mu$ as in Lemma~\ref{lem:weakcoer}, we see that the form
	\changed{\[
		q_\mu:\quad H^1_{0}([0,1])\times H^1_{0}([0,1])\to \R,\quad
		(v,w) \mapsto q(v,w)- \mu\langle v,w\rangle_{L^2([0,1])}
		\]}
	fulfills \changed{$q_\mu(v,v)\leq -\alpha \|v\|^2_{H^1([0,1])}$}
	for all $v\in H^1_{0}([0,1])$.
	Further, 
	let $x_D \in H^1_{0l}([0,1])$, such that $x_D(1)=1$ (which exists by a~simple linear interpolation). By the Lax-Milgram lemma, there exists some $x_0\in H^1_0([0,1])$ such that
	\begin{equation*}
	\forall \varphi \in H^1_0([0,1]): \changed{q_\mu}(x_0, \varphi)  = u\langle \varphi',ax_D'\rangle_{L^2([0,1])} - u \langle bx_D'+(c-\mu)x_D,\varphi\rangle_{L^2([0,1])}
	\end{equation*}
	Clearly, $x = x_0+u x_D\in H^1_{0l}([0,1])$ fulfills $x(1) = u$ by construction and
	\begin{align*}
	\forall \varphi \in H^1_0([0,1]): -\langle \varphi',ax'\rangle_{L^2([0,1])}=
	-  \langle \varphi, bx'+(c-\mu)x\rangle_{L^2([0,1])} 
	\end{align*}
	This however means that $-bx'-cx+\mu x\in L^2([0,1])$ is the weak derivative of $ax'$, which implies that 
	$\spvek{x}{u}\in \dom(A\&B)$.
	

	\noindent\textbf{~\ref{def:sysnodea}):}
	Assume that $\spvek{x_n}{u_n}$ is a~sequence such that $\spvek{x_n}{u_n} \to \spvek{x}{u}\in L^2([0,1]) \times \R$ and $A\&B \spvek{x_n}{u_n}\to {z} \in L^2([0,1])$. By the already proven statement \textbf{~\ref{def:sysnodec})}, there exists some $x_{D}\in H^1_{0l}([0,1])$, such that $\spvek{x_D}1\in\dom(A\&B)$. Then $\spvek{x_n-u_nx_D}0\in\dom(A\&B)$ for all $n\in\N$, $\spvek{x-ux_D}0\in\dom(A\&B)$, 
	the sequence $(x_n-u_nx_D)$ converges in $L^2([0,1])$ to $x-ux_D$, and $A\&B\spvek{x_n-u_nx_D}0$ converges in $L^2([0,1])$ to $z-A\&B\spvek{ux_D}u$. The latter means that 
	$A({x_n-u_nx_D})$ converges in $L^2([0,1])$ to $z-A\&B\spvek{ux_D}u$.
	Since, by the already proven statement \textbf{\ref{def:sysnoded})}, $A$ generates a~strongly continuous semigroup, it is closed by \cite[Chap.~I, Thm.~1.4]{EngeNage00}. This yields 
	\[x-ux_D\in\dom(A)\text{ with }A(x-ux_D)=z-A\&B\spvek{ux_D}u,\] and thus
	\begin{align*}
	\spvek{x}u&=\spvek{x-ux_D}0+u\spvek{x_D}1\in\dom(A\&B)\\
	\text{with}\qquad    
	A\&B\spvek{x}u&=A\&B\spvek{x-ux_D}0+A\&B\spvek{ux_D}u=\big(z-A\&B\spvek{ux_D}u\big)+uA\&B\spvek{x_D}1=z,
	\end{align*}
	which shows that $A\&B$ is closed.\\
	\noindent\textbf{~\ref{def:sysnodeb}):} 
	Assume that the sequence $\spvek{x_n}{u_n}$ is bounded in $\dom(A\&B)$, that is,
	\begin{eqnarray}\label{eq:boundedness1}
	\spvek{x_n}{u_n}&  &\mathrm{is\ bounded\ in\ }L^2([0,1]) \times \R,\\
	\label{eq:boundedness2}
	(A\&B \spvek{x_n}{u_n})&  &\mathrm{is\ bounded\ in\ } L^2([0,1]).
	\end{eqnarray}
	We first show that $(x_n)$ is bounded in $H^1([0,1])$. Again, by
	the already proven statement \textbf{~\ref{def:sysnodec})}, we can choose some $x_{D}\in H^1_{0l}([0,1])$ with $\spvek{x_D}1\in\dom(A\&B)$. Then both sequences $A(x_n-u_nx_D)$ and $(x_n-u_nx_D)$ are bounded in $L^2([0,1])$, whence the scalar sequence \[\big(\langle x_n-u_nx_D,A(x_n-u_nx_D)\rangle_{L^2([0,1])}\big)\] is bounded. Now using that, by the definition of the weak derivative,
	\[\langle x_n-u_nx_D,A(x_n-u_nx_D)\rangle_{L^2([0,1])}=q(x_n-u_nx_D,x_n-u_nx_D),\]
	whence $(q(x_n-u_nx_D,x_n-u_nx_D))$ is a~bounded sequence. 
	By combining this with boundedness of $(x_n-u_nx_D)$ in $L^2([0,1])$ and Lemma~\ref{lem:weakcoer}, we deduce that $(x_n-u_nx_D)$ (and consequently $(x_n)$) forms a bounded sequence in $H^1([0,1])$. As a result, the sequence $(ax_n')$ is bounded in $L^2([0,1])$. Further combining this information with \eqref{eq:boundedness2}, we conclude that $(ax_n')$ is bounded in $H^1([0,1])$. Since point evaluation represents a bounded operator in $H^1([0,1])$, we can then affirm that the sequence
	\[(C\&D\spvek{x_n}{u_n})=((ax_n')(0))\]
	is also bounded.\\
	Now it remains to show bijectivity of $S:\dom(S)\to L^2([0,1])\times\R$. To this end, 
	let
	$v\in L^2([0,1])$, $w\in\R$. Then, by
	\cite[Thm.~5.1]{Hale80}, there exists a~unique solution of the initial value problem
	\[{\spvek{x\\[-1mm]}{x_1}}'=\sbmat{0}{a^{-1}}{-c}{-ba^{-1}}{\spvek{x\\[-1mm]}{x_1}}+\spvek{0\\[-1mm]}v,\quad \spvek{x(0)}{x_1(0)}=\spvek{0\\[-1mm]}w\]
	in the sense of Carath\'eodory, which furthermore exists on $[0,1]$ by \cite[Cor.~6.4]{Hale80}. This gives $x_1',x'\in L^2([0,1])$, whence also $x_1\in H^1([0,1])$, $x\in H^1_{0l}([0,1])$.
	As a~consequence, for $x_1=ax'$, we have
	\begin{equation}(ax')'+bx'+cx=v,\quad x\in H^1_{0l}([0,1]),\; (ax')(0)=w,\label{eq:invode}\end{equation}
	which means that, for $u=x(1)$, we have $\spvek{x}u\in\dom(S)$ with $\spvek{v}w=S\spvek{x}u$.\\
	On the other hand, to show injectivity of $S$, let $\spvek{x}u\in\ker S$. Then, for $x_1=ax'$,  \eqref{eq:invode} holds with $v=0$ and $w=0$. Since, by \cite[Thm.~5.3]{Hale80}, the solution of \eqref{eq:invode} is unique, we have $x=0$, and thus also $u=x(1)=0$, and the result is proven.

\end{proof}
We further identify the adjoint system node. To accomplish this, we introduce an operator-theoretic auxiliary result.
\begin{lem}\label{lem:Adj}
	Let $V,W$ be Hilbert spaces, and let $F:V\supset\dom(A)\to W$, 
	$G:W\supset\dom(B)\to V$, are closed, bijective, and
	\begin{equation}
	\forall\,v\in\dom(F),w\in\dom(G):\quad \langle Fv,w\rangle_W=\langle v,Gw\rangle_V.\label{eq:FGadj}    \end{equation}
	Then $F^*=G$.
\end{lem}
\begin{proof}
	\eqref{eq:FGadj} means that $\dom(G)\subset\dom(F^*)$ with $Gw=F^*w$ for all $w\in\dom(G)$. Hence, it remains to show that
	$\dom(F^*)\subset\dom(G)$. Let $w\in\dom(F^*)$. Then there exists some $x\in V$ with
	\[\forall\, v\in\dom(F):\quad 
	\langle Fv,w\rangle_W=\langle v,x\rangle_V.\]
	By bijectivity of $G$, there exists some $\hat{w}\in\dom(G)$ with $x=G\hat{w}$, and thus,
	\begin{equation*}
	\forall\, v\in\dom(F):\quad    \langle Fv,w\rangle_W
	=\langle v,x\rangle_V
	=\langle v,G\hat{w}\rangle_V
	=\langle Fv,\hat{w}\rangle_W
	\end{equation*}
	Surjectivity of $F$ now gives $w=\hat{w}\in\dom(G)$.
\end{proof}
Now we present our result on the adjoint of the system node \eqref{eq:advdiff_dissnode}. For sake of convenience, we impose the additional property that $b$ is Lipschitz continuous. 
\begin{prop}
	Under Assumptions~\ref{ass:heat} and, additionally, $b\in W^{1,\infty}([0,1])$, the adjoint of the operator $S=\sbvek{A\&B}{C\&D}$, with $A\&B$ and $C\&D$ as in \eqref{eq:advdiff_dissnode}, is given by
	$S^*=\sbvek{{[A\&B]^\mathrm{d}}}{{[C\&D]^\mathrm{d}}}$ with
	\begin{subequations}
		\label{eq:heat:Sd}    
		\begin{equation}
		\dom({[A\&B]^\mathrm{d}}):=\setdef{\spvek{x_\mathrm{d}}{y_\mathrm{d}}\in H^1([0,1])\times \R}{\quad\parbox{4.5cm}{$(a x_\mathrm{d}')'-(bx_\mathrm{d})'\in
				L^2([0,1])$\\$\;\wedge\;x_\mathrm{d}(0)=y_\mathrm{d}\;\wedge\;x_\mathrm{d}(1)=0$}}
		\end{equation}
		and
		\begin{equation}
		\begin{aligned}
		\forall \, \spvek{x_\mathrm{d}}{y_\mathrm{d}}\in\dom([A\&B]^\mathrm{d}):\quad [A\&B]^\mathrm{d}\spvek{x_\mathrm{d}}{y_\mathrm{d}}&= \big(a x_\mathrm{d}'\big)'-(b x_\mathrm{d})' + cx_\mathrm{d},\\ [C\&D]^\mathrm{d}\spvek{x_\mathrm{d}}{y_\mathrm{d}}&= {-(ax_\mathrm{d}')(1)}.
		\end{aligned}
		\end{equation}\end{subequations}
\end{prop}
\begin{proof}
	Consider the operator $\widetilde{S}=\sbvek{{[A\&B]^\mathrm{d}}}{{[C\&D]^\mathrm{d}}}$ with $\dom(\widetilde{S})=\dom({[A\&B]^\mathrm{d}})$ and ${[A\&B]^\mathrm{d}}$, ${[C\&D]^\mathrm{d}}$ as in \eqref{eq:heat:Sd}. Then, for all $\spvek{x_\mathrm{d}}{y_\mathrm{d}}\in\dom(\widetilde{S})$,
	$\spvek{x}{y}\in\dom({S})$, we have
	\begin{align*}
	&\phantom{=}   \langle\spvek{x_\mathrm{d}}{y_\mathrm{d}},S\spvek{x}{u}\rangle_{L^2([0,1])\times \R}\\ &= \langle x_\mathrm{d},(ax')'+bx'+cx\rangle_{L^2([0,1])} + y_\mathrm{d} (ax')(0)\\
	&= -\langle x_\mathrm{d}',ax'\rangle_{L^2([0,1])} + \langle bx_\mathrm{d},x'\rangle_{L^2([0,1])} + \langle cx_\mathrm{d},x\rangle_{L^2([0,1])} + x_d(0)(ax')(0) \\&\hspace*{4cm}+ \underbrace{x_d(1)(ax')(1)}_{=0} - x_d(0)(ax')(0)\\
	&= -\langle ax_\mathrm{d}',x'\rangle_{L^2([0,1])} - \langle (bx_\mathrm{d})',x\rangle_{L^2([0,1])} + \langle cx_\mathrm{d},x\rangle_{L^2([0,1])}  + \underbrace{bx_\mathrm{d}(1)x(1)}_{=0} - \underbrace{bx_\mathrm{d}(0)x(0)}_{=0}\\
	&= \langle (ax_\mathrm{d}')',x\rangle_{L^2([0,1])} - \langle (bx_\mathrm{d})',x\rangle_{L^2([0,1])} + \langle cx_\mathrm{d},x\rangle_{L^2([0,1])} - (ax_\mathrm{d}')(1)\underbrace{x(1)}_{=u} + \underbrace{(ax_\mathrm{d}')(0)x(0)}_{=0}\\
	&= \langle \widetilde{S}\spvek{x_\mathrm{d}}{y_\mathrm{d}},\spvek{x}{u}\rangle_{L^2([0,1])\times \R}.
	\end{align*}
	Moreover, $S$ is bijective by Proposition~\ref{prop:heatsysnode}. The product rule for the weak derivative \cite[Thm.~4.25]{Alt16} yields $(bx_\mathrm{d})' = b'x_\mathrm{d} + b x_\mathrm{d}'\in L^2([0,1])$, whence
	\[\forall\,\spvek{x_\mathrm{d}}{y_\mathrm{d}}\in\dom(\widetilde{S}):\quad [A\&B]^\mathrm{d}\spvek{x_\mathrm{d}}{y_\mathrm{d}}=\big(a x_\mathrm{d}'\big)'-b x_\mathrm{d}' + (c-b')x_\mathrm{d},\]
	whereas the assumption $b\in W^{1,\infty}([0,1])$ gives rise to $c-b'\in L^\infty([0,1])$. As a~consequence,
	\[\widetilde{\widetilde{S}}:=\sbmat{\Refl_1}00I\widetilde{S}\sbmat{\Refl_1}00I\]
	is a~system node of a~type that is subject of Proposition~\ref{prop:heatsysnode}. In particular, we can conclude from Proposition~\ref{prop:heatsysnode} that $\widetilde{\widetilde{S}}$ (and thus also $\widetilde{S}$) is bijective. Altogether, we have that 
	$S$ and $\widetilde{S}$ fulfill the requirements of Lemma~\ref{lem:Adj}, and we can conclude that $S^*=\widetilde{S}$.

\end{proof}

\noindent \textbf{The optimal control problem.} Having provided a formulation of the advection-diffusion-reaction equation via system node, we will now consider the corresponding optimal control problem.

Let $y_\mathrm{ref} \in L^2([0,T])$ and $x_f\in L^2([0,1])$. We further consider a full terminal weight $F=I_{L^2([0,1])}$, and
we assume that the initial state fulfills $x_0\in H^1_{0l}([0,1])$.

For $\alpha > 0$ and the system node as in \eqref{eq:advdiff_dissnode}, we consider the optimal control problem 
\begin{equation}\begin{aligned}
\textrm{minimize}&\quad\frac{1}{2}\int_0^T(y(t)-y_\mathrm{ref}(t))^2+ \alpha u(t)^2\,{\rm d}t +  \frac{1}{2}\|x(T)-x_f\|^2_{L^2([0,1])}\\
\text{subject to}&\quad\spvek{\dot{x}(t)}{y(t)} = \sbvek{A\& B\\[-1mm] }{C\& D} \spvek{{x}(t)}{u(t)}, \quad
x(0) = x_0,\quad u \in \Uad := L^2([0,T]).
\end{aligned}\label{eq:OCP_parab}
\end{equation}
This problem falls into the class treated in Section~\ref{sec:inputpen}. To conclude existence of optimal controls and optimality conditions, we verify the Assumptions~\ref{ass:ocp}. 
\begin{itemize}
	\item[\eqref{ass:ocp1}] We have shown in Proposition~\ref{prop:heatsysnode} that $\sbvek{A\& B}{C\& D}$ as in \eqref{eq:advdiff_dissnode} is a system node. 
	\item[\eqref{ass:ocp4}] As \changed{$F=I$}, the density assumption is fulfilled in view of Remark~\ref{rem:denseass}\,\eqref{rem:denseassiib}.
	\item[\eqref{ass:ocp2}--\eqref{ass:ocp3}] These are satisfied by assumption.
	\item[\eqref{ass:ocp5}] This is satisfied by construction, as the extended output is coercive with respect to the input. 
	\item[\eqref{ass:ocp6}] We choose $\hat{u} \equiv {x}_0(1)$. 
	Let $x_D\in L^2([0,1])$, such that $\spvek{x_D}1\in\dom(A\&B)$.
	For \[f=\hat{u}\,A\&B\spvek{x_D}1,\quad z_0=x_0-\hat{u}x_D,\] consider the solution $z:[0,T]\to L^2([0,1])$ of
	\begin{equation}\dot{z}(t)=Az(t)+f,\quad z(0)=z_0.\label{eq:zsol}
	\end{equation}
	Since $x_0\in H^1_{0l}([0,1])$, the construction of $x_D$ yields $z_0\in H^1_{0}([0,1])$. By further using that, for $\mu$, $\alpha$ as in Lemma~\ref{lem:weakcoer}, the definition of the weak derivative gives
	\begin{multline}
	\forall\, v\in \dom(A):\ \langle v,(-A+\mu I)v\rangle=
	\langle av',v'\rangle_{L^2([0,1])} - 
	\langle bv',v\rangle_{L^2([0,1])}
	- \langle (c-\mu)v,v\rangle_{L^1([0,1]),L^\infty([0,1])},
	\end{multline}
	we obtain from \eqref{eq:qcont} and Lemma~\ref{lem:weakcoer} that
	\[\forall\, v\in \dom(A):\quad
	\alpha \|v\|^2_{H^1([0,1])}\leq
	\langle v,(-A+\mu I)v\rangle\leq (c+\mu) \|v\|^2_{H^1([0,1])}
	\]
	with $c>0$ as in \eqref{eq:qcont}.
	A~consequence is that the domain of the fractional power $(-A+\mu I)^{1/2}$ is the closure of $\dom(A)$ with respect to the norm in $H^1([0,1])$. This gives rise to 
	\[z_0\in H^1_0([0,1])=\dom(-A+\mu I)^{1/2}.\]
	Then \cite[Part II-1,~Thm.~3.1,~p.143]{Bensoussan2007}
	yield that $z\in L^2([0,T];\dom(A))$. Consequently, $\hat{x} := z + \hat{u}x_D$ satisfies
	\begin{equation}
	\spvek{\hat{x}}{\hat{u}}=\spvek{z}{0}+\spvek{\hat{u}x_D}{\hat{u}}\in L^2([0,T];\dom(A\&B)).\label{eq:xuindom}\end{equation}     
	We have
	$x(0)=z_0+\hat{u}x_D=x_0$
	and, by invoking that, by \eqref{eq:sol-2}, $z\in W^{1,1}([0,T];X_{-2})$, we have
	\begin{align*}
	\dot{\hat{x}}&=\dot{z}=Az+f=A(\hat{x}-\hat{u}x_D)+u\,A\&B\spvek{x_D}1\\&=A\&B\spvek{\hat{x}-\hat{u}x_D}{0}+\hat{u}=A\&B\spvek{x_D}1=
	A\&B\spvek{\hat{x}}{\hat{u}},
	\end{align*}
	we obtain that 
	$\hat{u}\in \Uad$ with $\spvek{x_0}{\hat{u}}\in \dom(\frakT_{I,T})$.
	Further, by using \eqref{eq:xuindom} together with $C\&D\in L(\dom(A\&B),\R)$, we also have 
	\begin{align*}
	y=   \frakA_{-1}x_0+{\frakD_T}\hat{u} = 
	C\&D\spvek{\hat{x}}{\hat{u}} \in L^2([0,T]).
	\end{align*}
\end{itemize}
We now may conclude existence of an optimal control using Theorem~\ref{thm:optcont} and its uniqueness, which $\alpha>0$ by using the findings in Section~\ref{sec:inputpen}.  
\begin{cor}
	Let $\alpha>0$, $y_\mathrm{ref} \in L^2([0,T])$, $x_0\in H^1_{0l}([0,1])$ and $x_f\in L^2([0,1])$. 
	Then, under Assumptions~\ref{ass:heat}, $X=L^2([0,1])$, $U=\R$, and $A\&B$, $C\&D$ as in \eqref{eq:advdiff_dissnode},
	the optimal control problem~\eqref{eq:OCP_parab} has solution. That is, there exists a unique optimal control $u_{\rm opt}\in L^2([0,T];\R^2)$ in the sense of Definition~\ref{def:cost}. \end{cor}
Further, optimality conditions may be expressed using the variational inequalities of Theorem~\ref{thm:optcond} and Corollary~\ref{cor:aufloesen}. In particular, in view of the results in Section~\ref{sec:inputpen}, the optimal control satisfies
\[\begin{aligned}
\spvek{\dot{x}_{\rm opt}(t)}{y_{\rm opt}(t)}&
= \sbvek{A\&B\\[-1mm]}{C\&D} \spvek{{x}_{\rm opt}(t)}{u_{\rm opt}(t)},& x_{\rm opt}(0)&=x_0,\\
\spvek{\dot{\mu}(t)}{u_{\rm opt}(t)}
&= \sbvek{{-\phantom{\alpha^{-1}}[A\&B]^\mathrm{d}}}{{{-}\alpha^{-1} [C\&D]^\mathrm{d}}}
\spvek{{\mu}(t)}{y_{\rm opt}(t)-y_{\rm ref}(t)},& \mu(T)&=x_{\rm opt}(T)-x_f.
\end{aligned}\]

\begin{figure}[htb]
	\centering
	\vspace*{-.1cm}
	\includegraphics[width=0.4\linewidth]{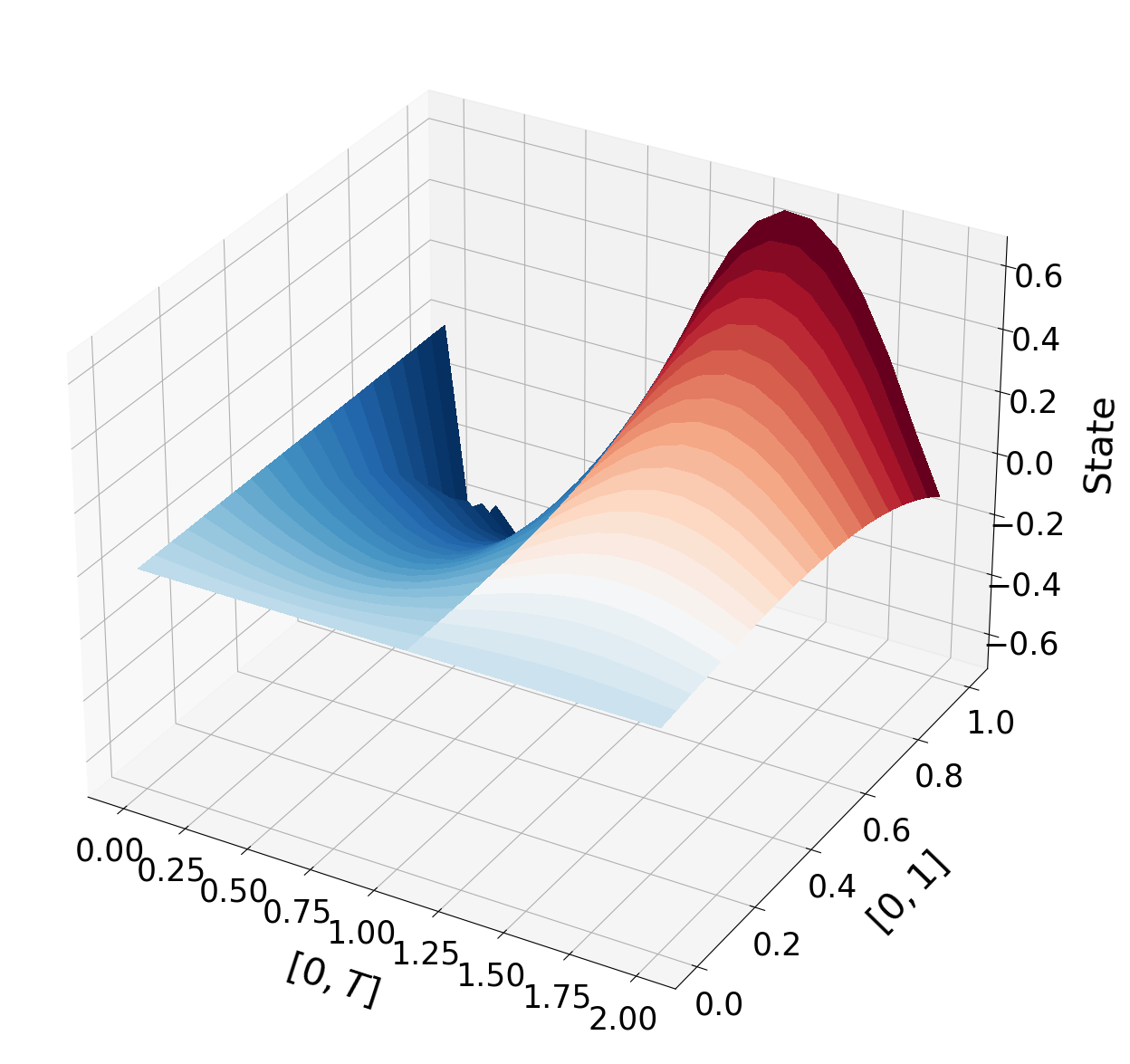}\hspace*{.5cm}
	\includegraphics[width=0.4\linewidth]{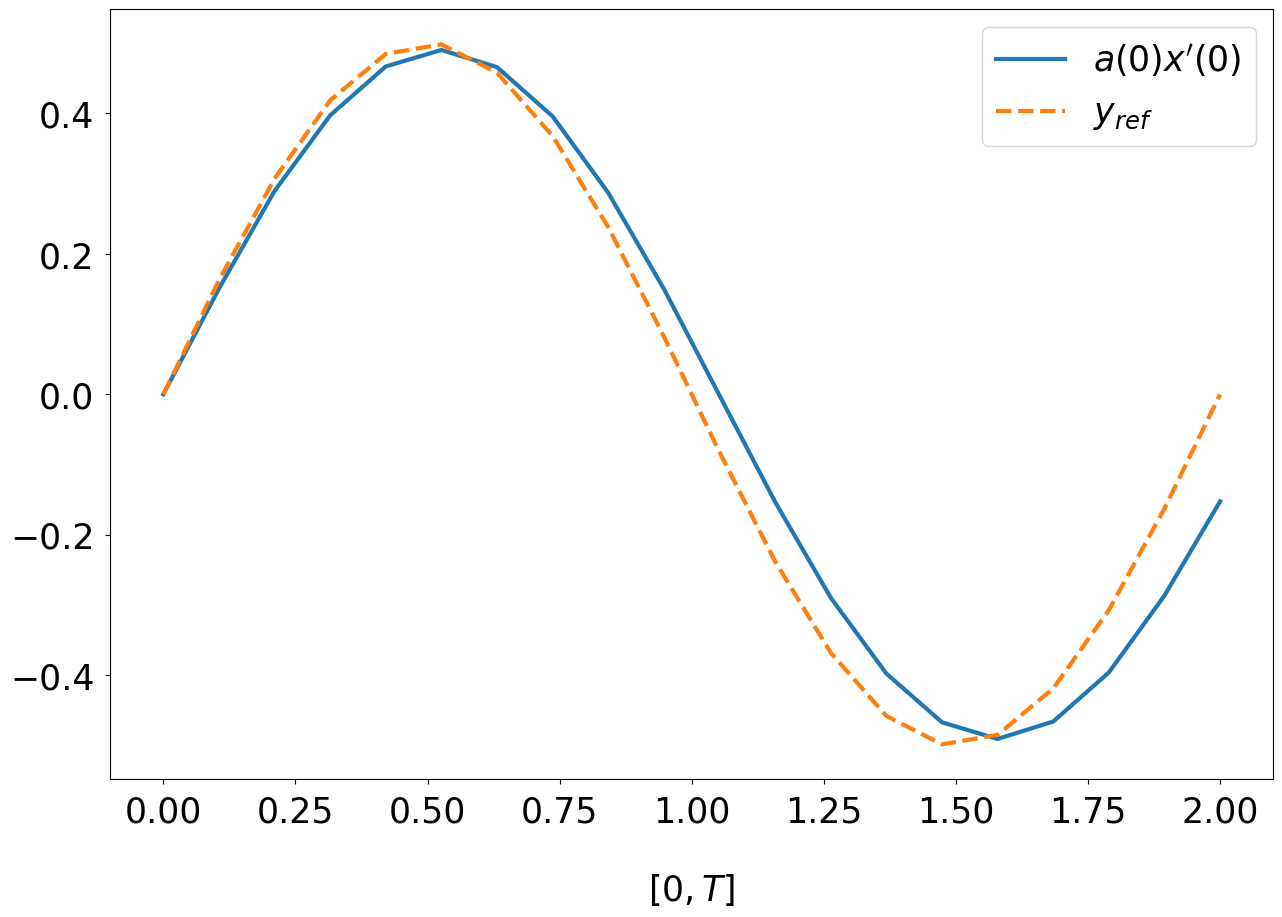}
	\caption{Reaction coefficient $c\equiv 1$. Left: Optimal state over time $[0,T]$ and space $[0,L]$. Right: Optimal output over time.}
	\label{fig:plots_heat}
\end{figure}

\begin{figure}[htb]
	\centering
	\includegraphics[width=0.4\linewidth]{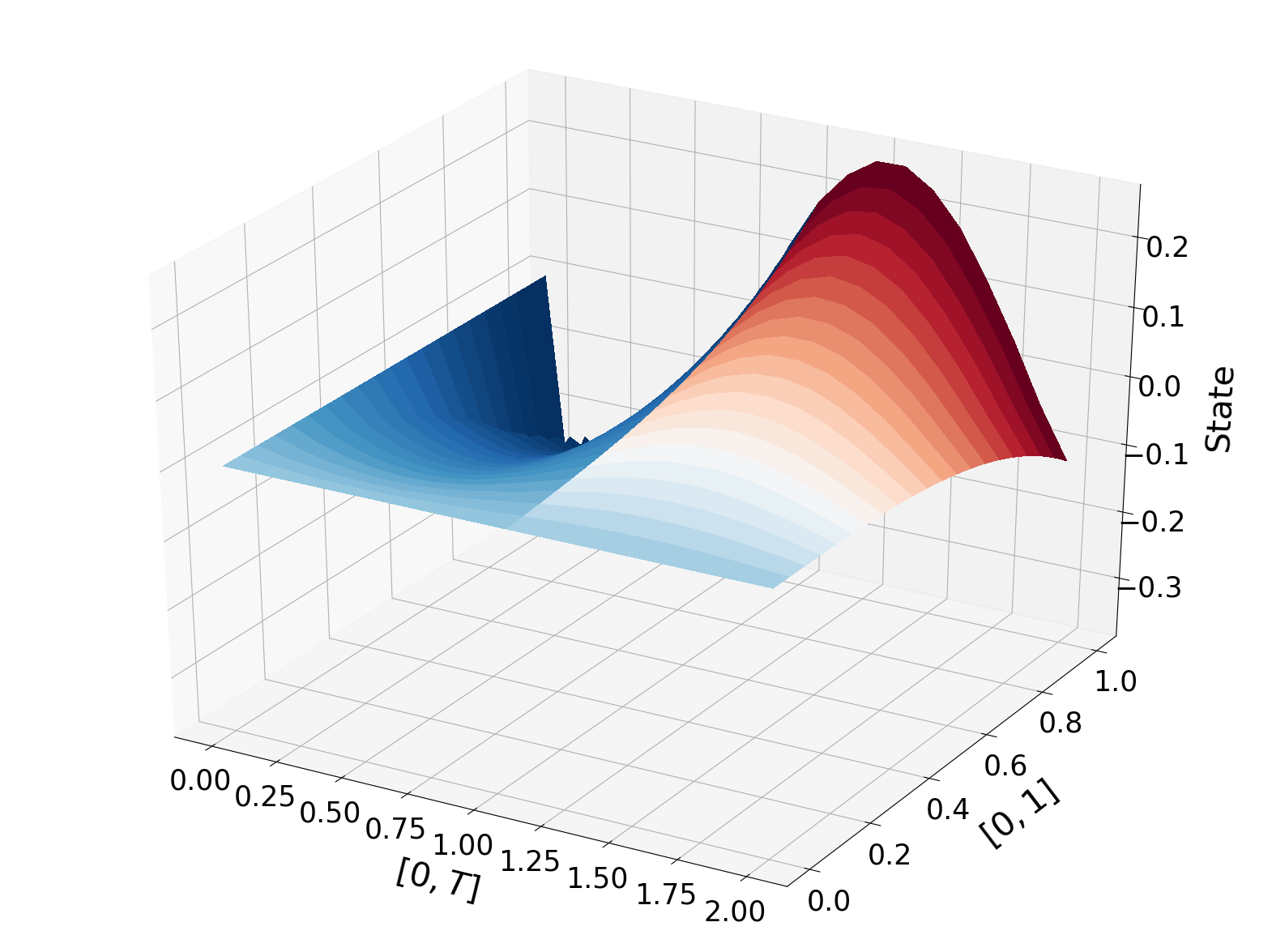}\hspace*{.5cm}
	\includegraphics[width=0.4\linewidth]{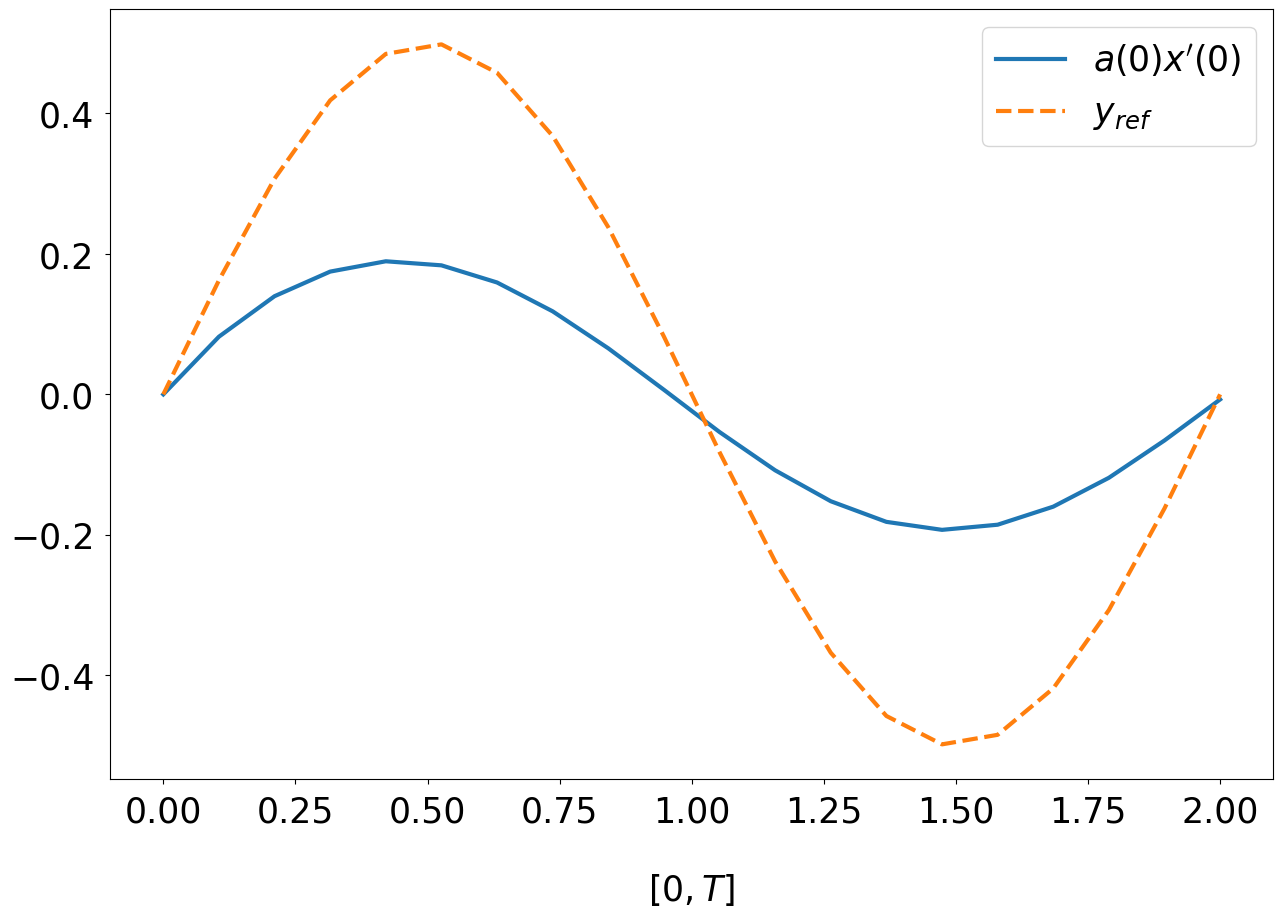}
	\caption{Reaction coefficient $c\equiv 5$. Left: Optimal state over time $[0,T]$ and space $[0,L]$. Right: Optimal output over time.}
	\label{fig:plots_heat2}
\end{figure}

We provide now a brief numerical example illustrating the solution of this problem. We solve the corresponding problem with an explicit Euler discretization and linear finite elements using {FEniCS}~\cite{alnaes2015fenics} and dolfin-adjoint~\cite{mitusch2019dolfin}. The optimal control problems aims at tracking $y_\mathrm{ref}$ as a Neumann trace on the left side of the domain while actuating the Dirichlet trace on the right side. For the PDE, we choose $a(x)= 1,\ b(x) = -x, \ x_0 = 0$ and either \changed{$c(x)=1$} (Figure~\ref{fig:plots_heat}) or $c(x) = 5$ (Figure~\ref{fig:plots_heat2}). For the cost functional, we set the reference trajectory $y_\mathrm{ref}(t) = \sin(\pi t)$, the  
parameter $\alpha = 0.1$ and the time horizon $T=2$.

The result for the smaller reaction coefficient is shown in Figure~\ref{fig:plots_heat}. To reach the desired Neumann trace at zero, the state at the respective other side of the domain is chosen negative by means of the Dirichlet actuation, leading to a positive Neumann trace at the other side. Then, after time $t=1$, the Dirichlet boundary control is chosen positive such that the Neumann trace on the opposite side becomes negative. 

In Figure~\ref{fig:plots_heat2}, we show the same quantities for a higher reaction coefficient. As this renders the uncontrolled PDE unstable, the corresponding optimal state is smaller due to the penalization of the control. Correspondingly, the optimal output has a higher disparity with the reference signal.

\FloatBarrier\subsection{Wave equation on an L-shaped domain}

We consider a~boundary controlled wave equation as in \cite{KuZw15,Farkas2023}, which is given by
\begin{equation}\begin{aligned}
\rho(\xi) \tfrac{\partial^2}{\partial t^2}\mathbf{w}(\xi,t) &= \Div(\mathcal{T}(\xi) \grad \mathbf{w}(\xi,t)) -  d(\xi) \tfrac{\partial}{\partial t}\mathbf{w}(\xi,t), &&  \xi \in \Omega, t \ge 0, \\
0 &= \mathbf{w} (\xi ,t ) && \xi\in \Gamma_0, t\ge0,\\
u(\xi,t) &= \nu \cdot (\mathcal{T}(\xi)\grad \mathbf{w}(\xi,t)) && \xi\in\Gamma_1, t\ge0,\\
y(\xi,t) &= \tfrac{\partial}{\partial t}\mathbf{w}(\xi,t), && \xi\in \Gamma_0, t\ge0, \\
\mathbf{w}(\xi,0)&= \mathbf{w}_0(\xi), \quad \mathbf{w}_t(\xi,0) = \mathbf{v}_0 (\xi) && \xi\in \Omega,
\end{aligned}\label{eq:waveeq}
\end{equation}
where $\nu\colon \partial\Omega\to\R^2$ is the unit outward normal vector of the L-shaped domain $\Omega\subset\R^2$. Its boundary is decomposed into two parts, $\Gamma_0$ and $\Gamma_1$, as depicted below.

\begin{minipage}[t]{.45\linewidth}
	\begin{align*}
	\Omega&={\mathrm{int}}(\overline{\Omega_1}\cup\overline{\Omega_2}),\\
	\Omega_1&=(0,1)\times(0,2),\\
	\Omega_2&=(1,2)\times(0,1),\\
	\Gamma_1&=(0,1)\times\{2\}\cup \{2\}\times (0,1)\subset\partial\Omega,\\
	\Gamma_0&=\partial\Omega\setminus\overline{\Gamma_1}.
	\end{align*}
\end{minipage}\hspace{.05\linewidth}
\begin{minipage}[t]{.45\linewidth}
	\vspace*{-.8cm}     
	\begin{tikzpicture}[scale=2]
	\draw[dashed, thick] (0,0)--(2,0)--(2,1)--(1,1)--(1,2)--(0,2)--cycle;
	\draw[thick] (2,0)--(2,1);
	\draw[thick] (0,2)--(1,2);    
	\node at (2.2,0.5) {$\Gamma_1$};
	\node at (0.5,2.2) {$\Gamma_1$};
	\node at (1.3,1.3) {$\Gamma_0$};
	\node at (1.,-.2) {$\Gamma_0$};
	\node at (-.2,1.) {$\Gamma_0$};
	\node at (0.6,0.8) {\LARGE$\Omega$};
	\end{tikzpicture}
\end{minipage}

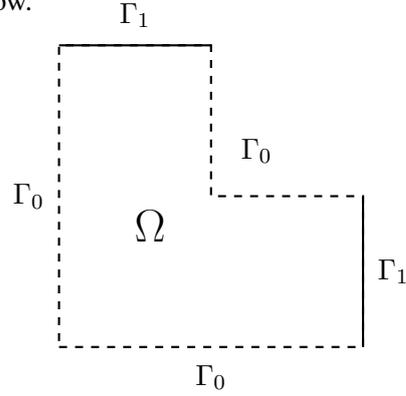
\captionof{figure}{\label{fig:recwave}L-shaped domain for the wave equation.}
~\\
As usual, $\mathbf{w}(\xi,t)$ denotes the displacement of the wave at point $\xi \in \Omega$ and time $t\ge 0$, $u$ is the input given by a force on the boundary part $\Gamma_1$ and $y$ is the output measured as the velocity at $\Gamma_1$. The physical parameters are included via the Young's modulus $\mathcal{T}(\cdot)$ and the mass density $\rho(\cdot)$, which are both assumed to be measurable, positive, and they have bounded inverses. The term with  $d$ can be interpreted as an internal damping which is assumed to be a bounded nonnegative and measurable function on $\Omega$.  Our assumptions are summarized below.\\
\begin{ass}\label{ass:wave}\
	\begin{enumerate}[(a)]
		\item $\Omega$, $\Gamma_0$, $\Gamma_1$ are as in Fig.~\ref{fig:recwave}.
		\item $\mathcal{T},\rho\in L^\infty(\Omega)$ are positive almost everywhere, with $\rho^{-1},\mathcal{T}^{-1}\in L^\infty(\Omega)$;
		\item $d\in L^\infty(\Omega)$ is nonnegative almost everywhere.
	\end{enumerate}
\end{ass}

Now we show that the system can be formulated as one which is port-Hamiltonian in the sense of Definition~\ref{def:pHsys}. 
The state will consist of the kinetic and potential energy variables, i.e.,
\begin{equation}
x = \spvek{\mathbf{p}\\[-1mm]}{\mathbf{q}} = \spvek{\rho \tfrac{\partial}{\partial t}\mathbf{w}}{\mathcal{T}^{-1}\grad \mathbf{w}}. \label{eq:wavestate}
\end{equation}
That is, $\mathbf{p}$ is the infinitesimal momentum, whereas $\mathbf{q}$ reflects the stress.
As state space we choose
\begin{equation}X := {L^2(\Omega)}\times{L^2(\Omega;\R^2)},\label{eq:Xwave}\end{equation}
which is further identified with its anti-dual in the canonical manner. Further, let $H:X\to X$ with 
\begin{equation}H\spvek{\mathbf{p}\\[-1mm]}{\mathbf{q}}=\spvek{\rho^{-1}\mathbf{p}}{\mathcal{T}\mathbf{q}},\label{eq:Hwave}\end{equation}
and we have, by $\rho>0$, $\mathcal{T}>0$, $\rho,\rho^{-1},\mathcal{T},\mathcal{T}^{-1}\in L^\infty(\Omega)$ that $H$ is self-adjoint with $H,H^{-1}\in L(X)$.
To introduce the corresponding dissipation node, we first mention that, 
by $H(\operatorname{div};\Omega)$ we mean the space of square integrable functions with weak divergence in $L^2(\Omega)$. 
The trace operator is denoted by $\gamma:H^1(\Omega)\to H^{1/2}(\partial \Omega)$, and we have that $\gamma$ is bounded and surjective due to the trace theorem~\cite[Thm.\ 1.5.1.3]{Gris85}. For a relatively open subset of the boundary $\Gamma \subset \partial \Omega$, $\gamma_\Gamma:H^1(\Omega)\to H^{1/2}(\Gamma)$ maps $f\in H^1(\Omega)$ to the restriction of $\gamma f$ to $\Gamma$ and we set
\[H^1_\Gamma(\Omega) = \{f\in H^1(\Omega)\,|\,\gamma_{\Gamma} = 0\}.\]
Then we obtain a~bounded and surjective operator
\[\gamma_{\Gamma}: H^1_{\partial \Omega \setminus \Gamma}(\Omega)\to H^{1/2}_0(\Gamma),\]
where $H^{1/2}_0(\Gamma)$ contains the elements of $H^{1/2}(\Gamma)$ which can be extended by zero to an element of $H^{1/2}(\partial \Omega)$.
Further, we consider $H^{-1/2}(\Gamma) := H^{1/2}_0(\Gamma)^*$, and we may define the normal trace on $\Gamma \subset \partial \Omega$ of $x\in H(\operatorname{div},\Omega)$ by $w = \gamma_{N,\Gamma} x\in H^{1/2}_{0}(\Gamma)^*$, where
\begin{align*}
\forall z\in H^1(\Omega): \langle w,\gamma z\rangle_{H^{-1/2}(\Gamma),H^{1/2}(\Gamma)} = \langle \operatorname{div} x,z\rangle_{L^2(\Omega)} + \langle x,\grad z\rangle_{L^2(\Omega;\R^2)}.
\end{align*}
Now assume that
$\Omega$, $\Gamma_0$, $\Gamma_1$ are as in Fig.\ \ref{fig:recwave}.
As input space, we choose
\begin{equation}U={H^{-1/2}(\Gamma_1)},\label{eq:Uwave}\end{equation}
whence, consequently, $U^*={H_0^{1/2}}(\Gamma_1)$.
Then, for $M=\sbvek{\tF\&\tG}{\tK\&\tL}$ with
\begin{subequations}\label{eq:wave_node}
	\begin{align}\dom(M)=\dom(\tF\&\tG):=\setdef{\spvek{\spvek{\mathbf{v}}{\mathbf{F}}}{u}\in {X}\times U}{\,\mathbf{v} \in H^1_{\Gamma_0}(\Omega)\,,\mathbf{F}\in H(\operatorname{div},\Omega)\;\wedge\;\gamma_{N,\Gamma_1}\mathbf{F}=u}
	\end{align}
	and
	\begin{align}
	\forall \, \spvek{\spvek{\mathbf{v}}{\mathbf{F}}}{u}\in\dom(M):\quad \tF\&\tG\spvek{\spvek{\mathbf{v}}{\mathbf{F}}}{u}&= \spvek{\operatorname{div}\mathbf{F}-d\mathbf{v}}{\grad(\mathbf{v})},\;\; \tK\&\tL\spvek{\spvek{\mathbf{v}}{\mathbf{F}}}{u}= \gamma_{\Gamma_1}(\mathbf{v}).
	\end{align}
\end{subequations}
the system introduced at the beginning of this section is represented by a~port-Hamiltonian system \eqref{eq:pHODEnode} with $H$ as in \eqref{eq:Hwave}. Here, we have replaced the clamping condition $\mathbf{w}(t)|_{\Gamma_0}=0$ with $\rho^{-1}\mathbf{p}(t)|{\Gamma_0}=0$ for $t\ge0$, achieved by taking the derivative with respect to time. To the best of the authors' knowledge, no criteria for the well-posedness of boundary-controlled wave equations on two-dimensional spatial domains are currently available (in contrast to the one-dimensional spatial case, as discussed in \cite[Chap.~13]{Jacob2012}). \changed{In this context we would like to mention the counterexample on a two-dimensional half-space in \cite[Section 6.3]{LasiTrig03}, where it is shown that the input-to-output map corresponding to Neumann boundary control and Dirichlet observation is not bounded as an operator from $L^2(0,T;L^2(\partial \Omega))$ to itself.}

It has been shown in \cite{Farkas2023}, that $M$ is indeed a~dissipative system node, and we can thus conclude that it is a~dissipation node. Further, it follows from \cite[Thm.~3.2]{KuZw15}, that 
\begin{equation}
\begin{aligned}
&    M_l:{X}\times U\supset\dom(M_l)=\dom(M)\to {X}\times U^*,\\
&\qquad M_l\spvek{\spvek{\mathbf{v}}{\mathbf{F}}}{u}=
\left(\begin{smallmatrix}
\operatorname{div}\mathbf{F}\\
{\grad(\mathbf{v})}\\
-\gamma_{\Gamma_1}(\mathbf{v})
\end{smallmatrix}\right)
\end{aligned}\label{eq:waveadj}\end{equation}
is skew-adjoint. In particular, we have a~dissipation output \eqref{eq:dissout} with $\tR\&\tS\in L(\dom(A\&B), W)$, for $W=X=L^2(\Omega)$ and
\begin{equation}
\forall\, \spvek{\spvek{\mathbf{v}}{\mathbf{F}}}{u}\in\dom(M):\quad 
\tR\&\tS \spvek{\spvek{\mathbf{v}}{\mathbf{F}}}{u}={\textstyle\sqrt{\frac{d}{{2}}}}\,\mathbf{v}.\label{eq:disswave}
\end{equation}
That is, in fact, $\tR\&\tS$ extends to a~bounded operator on whole $X\times U$. This means that, actually, 
$\tR\&\tS=[\tR\,\tS]$ with $\tS=0$, and  $\tR\in L(X, W)$ with $\spvek{\mathbf{v}}{\mathbf{F}}\mapsto{\textstyle\sqrt{\frac{d}{{2}}}}\,\mathbf{v}$.\medskip

\noindent \textbf{Reconstructing the displacement} 

We revisit the system \eqref{eq:waveeq} whose indeterminate is the displacement $\mathbf{w}$.
The choice of the state $x$ as in \eqref{eq:wavestate}, however, results in a partial loss of information about $\mathbf{w}$ - at first glance. \changed{To still be able to formulate cost functionals and terminal conditions involving the displacement}, we briefly discuss how to reconstruct $\mathbf{w}(t)\in H^1_{\Gamma_0}(\Omega)$ from the stress $\mathbf{q}(t)\in L^2(\Omega)$ at time $t\ge0$. 

\begin{prop}\label{prop:disprec}
	Under Assumptions~\ref{ass:wave}, let $\mathbf{q}\in L^2(\Omega)$. Then the following are equivalent:
	\begin{enumerate}[(i)]
		\item $\mathbf{q}\in\mathcal{T}^{-1}\grad H^1_{\Gamma_0}(\Omega)$. That is, there exists some $\mathbf{w}\in H^1_{\Gamma_0}(\Omega)$, such that
		\begin{equation}
		\mathbf{q}=\mathcal{T}^{-1}\grad\mathbf{w}\label{eq:gradfield}
		\end{equation}
		\item For all $\psi\in H(\operatorname{div},\Omega)$ with
		$\operatorname{div}\psi=0$ and $\gamma_{N,\Gamma_1} \psi=0$, it holds that
		\[\langle \psi,\mathcal{T}\mathbf{q}\rangle_{L^2(\Omega;\R^2)}=0.\]
	\end{enumerate}
	Moreover, in the case where \changed{(i) or (ii) (and hence both) are valid}, then $\mathbf{w}$ with the properties \changed{as stated in (i)} is the unique solution of the variational problem
	\begin{equation}\forall\,\varphi\in H^1_{\Gamma_0}(\Omega):\quad \langle\grad\varphi,\mathcal{T}^{-1}\grad \mathbf{w}\rangle_{L^2(\Omega;\R^2)}=\langle\grad\varphi,\mathbf{q}\rangle_{L^2(\Omega;\R^2)}.\label{eq:wLax}
	\end{equation}
\end{prop}
\begin{proof}
	To show that (i) implies (ii), assume that $\mathbf{q}=\mathcal{T}^{-1}\grad\mathbf{w}$
	for some $\mathbf{w}\in H^1_{\Gamma_0}(\Omega)$. By the definition of the normal trace, any $\psi\in  H(\operatorname{div},\Omega)$ and
	$\operatorname{div}\psi=0$, $\gamma_{N,\Gamma_1} \psi=0$ fulfills
	\begin{align*}
	\langle \psi,\mathcal{T}\mathbf{q}\rangle_{L^2(\Omega)} \changed{=
		\langle \psi,\grad\mathbf{w}\rangle_{L^2(\Omega;\R^2)}=
		-\underbrace{\langle \operatorname{div}\psi, \mathbf{w}\rangle_{L^2(\Omega)}}_{=0}+\langle \underbrace{\gamma_N\psi}_{=0},\gamma \mathbf{w}\rangle_{H^{-1/2}(\Gamma_1),H^{1/2}_0(\Gamma_1)}=0.    }
	\end{align*}
	Next, let us assume that \changed{$\mathbf{q}$ satisfies (ii)}. \changed{First, we observe that due to Assumption~\ref{ass:wave}, $\Gamma_0$ has positive measure. Thus, the Poincaré inequality \cite[Theorem 13.6.9]{TuWe09} implies that  the bilinear form on the left of \eqref{eq:wLax} is coercive. As $\mathcal{T}^{-1}\in L^\infty(\Omega)$, it is also continuous, such that the} Lax-Milgram lemma implies the existence of a unique $\mathbf{w}\in H^1_{\Gamma_0}(\Omega)$ satisfying \eqref{eq:wLax}. Since the solution $\mathbf{w}\in H^1_{\Gamma_0}(\Omega)$ to \eqref{eq:gradfield} is unique (provided it exists), our remaining objective is to prove that $\mathbf{w}\in H^1_{\Gamma_0}(\Omega)$ with \eqref{eq:wLax} indeed satisfies \eqref{eq:gradfield}, a task which will be accomplished in the following.
	
	First, we note that, by \eqref{eq:wLax}, 
	\[\forall\, \varphi\in \changed{H^1_{\Gamma_0}(\Omega)}:\quad \langle \grad\varphi,\mathbf{q}-\changed{\mathcal{T}^{-1}\grad \mathbf{w}}\rangle_{L^2(\Omega;\R^2)}{=}0,\]
	and the definition of the weak divergence yields \changed{$\mathbf{q}-\mathcal{T}^{-1}\grad\mathbf{w}\in H(\operatorname{div},\Omega)$} with $\operatorname{div}(\mathbf{q}-\changed{\mathcal{T}^{-1}\grad \mathbf{w}})=0$.
	Once again using \eqref{eq:wLax}, we obtain, for all $\varphi\in H^1_{\Gamma_0}(\Omega)$ that
	\begin{multline*}
	\langle {\gamma \varphi},\gamma_N(\mathbf{q}-\changed{\mathcal{T}^{-1}\grad \mathbf{w}})\rangle_{H^{1/2}_0(\Gamma_1),H^{-1/2}(\Gamma_1)}\\= \langle\varphi,\underbrace{\operatorname{div}(\mathbf{q}-\changed{\mathcal{T}^{-1}\grad \mathbf{w}})}_{=0}\rangle_{L^2(\Omega)}
	+\underbrace{\langle\grad\varphi,\mathbf{q}-\changed{\mathcal{T}^{-1}\grad \mathbf{w}}\rangle_{L^2(\Omega;\R^2)}}_{\stackrel{\eqref{eq:wLax}}{=}0}=0,
	\end{multline*}
	whence we can choose $\psi:=\mathbf{q}-\changed{\mathcal{T}^{-1}\grad \mathbf{w}}$ in (ii) to obtain that
	\begin{align*}
	\langle\mathbf{q}-\mathcal{T}^{-1}\grad\mathbf{w},\mathcal{T}\big(\mathbf{q}-\mathcal{T}^{-1}\grad\mathbf{w}\big)\rangle_{L^2(\Omega;\R^2)}&=\underbrace{\langle\psi,\mathcal{T}\mathbf{q}\rangle_{L^2(\Omega;\R^2)}}_{\stackrel{\text{(ii)}}{=}0}-\underbrace{\langle\psi, \mathcal{T}\big(\mathcal{T}^{-1}\grad\mathbf{w}\big)\rangle_{L^2(\Omega;\R^2)}}_{\stackrel{\text{(i)$\Rightarrow$(ii)}}{=}0}\\&=0.
	\end{align*}
	This indeed gives $\mathbf{q}=\mathcal{T}^{-1}\grad\mathbf{w}$.
\end{proof}
Coercivity of the form in \eqref{eq:wLax} together with the Lax-Milgram lemma yields the existence of a~bounded operator
\begin{equation}
\begin{aligned}
F_{\mathrm{disp}}: \quad L^2(\Omega;\R^2)&\to H^1_{\Gamma_0}(\Omega),\\
\mathbf{q}&\mapsto \mathbf{w}\text{ with \eqref{eq:wLax}}.
\end{aligned}\label{eq:Tdisp}
\end{equation}
Clearly, $F_{\mathrm{disp}}$ is also bounded as an~operator from $L^2(\Omega;\R^2)$ to $L^2(\Omega)$.
Proposition~\ref{prop:disprec} yields that we always have that $\mathbf{w}=F_{\mathrm{disp}}\mathbf{q}$ fulfills $\mathbf{q}=\mathcal{T}^{-1}\grad\mathbf{w}$ whenever $\mathbf{q}\in \mathcal{T}^{-1}\grad H^1_{\Gamma_0}(\Omega)$.  We give a brief result on the range of the adjoint of $F_{\mathrm{disp}}$ regarded as an operator from $L^2(\Omega;\R^2)$ to $L^2(\Omega)$. 
\begin{prop}\label{prop:disprec_range}
	Under Assumptions~\ref{ass:wave}, consider the operator $F_{\mathrm{disp}}\in L(L^2(\Omega;\R^2),L^2(\Omega))$ as in \eqref{eq:Tdisp}. Then 
	\[\im \left(F_{\mathrm{disp}}^*\right)\subset \setdef{\mathbf{z}\in H(\operatorname{div},\Omega)}{\gamma_{N,\Gamma_1}\mathbf{z}=0}.\]
\end{prop}
\begin{proof}
	Let $\varphi\in L^2(\Omega)$.
	By Proposition~\ref{prop:disprec}, we have 
	\[\forall\,\mathbf{w}\in  H^1_{\Gamma_0}(\Omega):\quad\mathbf{w}=F_{\mathrm{disp}}\changed{\mathcal{T}^{-1}}\grad\mathbf{w}.\]
	Thus, for all $\mathbf{w}\in  H^1_{0}(\Omega)$,
	we have
	\[\langle\varphi,\mathbf{w}\rangle_{L^2(\Omega)}
	=\langle\varphi,F_{\mathrm{disp}}\changed{\mathcal{T}^{-1}}\grad\mathbf{w}\rangle_{L^2(\Omega)}
	=\langle F_{\mathrm{disp}}^*\varphi,\changed{\mathcal{T}^{-1}}\grad\mathbf{w}\rangle_{L^2(\Omega)}.
	\]
	The definition of the weak divergence yields $\changed{\mathcal{T}^{-1}}F_{\mathrm{disp}}^*\varphi \in H(\operatorname{div},\Omega)$ with $-\operatorname{div}\changed{\mathcal{T}^{-1}}F_{\mathrm{disp}}^*\varphi=\varphi$. Further, the definition of the normal trace gives rise to the fact that, for all $\mathbf{w}\in  H^1_{\Gamma_0}(\Omega)$
	\begin{align*}
	\langle \gamma_N \changed{\mathcal{T}^{-1}} F_{\mathrm{disp}}^*\varphi,\gamma \mathbf{w}\rangle_{H^{-1/2}(\Gamma_1),H^{1/2}_0(\Gamma_1)}&=\langle \operatorname{div}\changed{\mathcal{T}^{-1}}F_{\mathrm{disp}}^*\varphi,\mathbf{w}\rangle_{L^2(\Omega)}+
	\langle F_{\mathrm{disp}}^*\varphi,\changed{\mathcal{T}^{-1}}\grad\mathbf{w}\rangle_{L^2(\Omega;\R^2)}\\
	&=-\langle \varphi,\mathbf{w}\rangle_{L^2(\Omega)}+\changed{\langle \varphi,F_{\mathrm{disp}}\changed{\mathcal{T}^{-1}}\grad\mathbf{w}\rangle_{L^2(\Omega)}}=0,
	\end{align*}
	which shows the result.
\end{proof}
Now we show that, if the stress is initialized with an element of $\mathbf{q}\in \mathcal{T}^{-1}\grad H^1_{\Gamma_0}(\Omega)$ then the stress remains in that space.
\begin{prop}\label{prop:Fdisp}
	Under Assumptions~\ref{ass:wave}, let a~port-Hamiltonian system 
	with operators as in \eqref{eq:Hwave} and 
	\eqref{eq:wave_node} be given, and let $\frakA$ be the corresponding semigroup on $X={L^2(\Omega)}\times{L^2(\Omega;\R^2)}$. Then the following holds:
	\begin{enumerate}[(a)]
		\item\label{prop:Fdispa} For all $T>0$, the space $L^2(\Omega)\times \mathcal{T}^{-1}\grad H^1_{\Gamma_0}(\Omega)\subset X$ is $\frakA(T)$-invariant.
		\item\label{prop:Fdispb} The space of reachable states is contained in $L^2(\Omega)\times \mathcal{T}^{-1}\grad H^1_{\Gamma_0}(\Omega)$, i.e.
		\[\forall\, u\in L^2([0,T];U)\text{ s.t.\ }\spvek{0}u\in\dom(\frakT_{I,T}):\quad \frakT_{I,T}\spvek{0}u\in L^2(\Omega)\times \mathcal{T}^{-1}\grad H^1_{\Gamma_0}(\Omega).\]
		\item\label{prop:Fdispc} For $F_{\mathrm{disp,ext}}:X\to X$ with 
		\[F_{\mathrm{disp,ext}}\spvek{\mathbf{p}\\[-1mm]}{\mathbf{q}}=\spvek{\mathbf{p}\\[-1mm]}{{F}_{\mathrm{disp}}\mathbf{q}},\]
		it holds that 
		the set $\mathcal{B}_{F_{\mathrm{disp,ext}}}$ (as defined in \eqref{eq:dense1}) is dense in $X$. In particular, the
		$F_{\mathrm{disp,ext}}$-terminal value map \changed{(as defined in Definition~\ref{def:FIS})}, the $F_{\mathrm{disp,ext}}$-input map, and the $F_{\mathrm{disp,ext}}^*$-output map (as defined in Definition~\ref{def:FIO}) are well-defined, cf.\ Table~\ref{tab:ourops} for an overview of these operators.
	\end{enumerate}
\end{prop}
\begin{proof}\
	\begin{itemize}
		\item[(a)] Assume that $x_0=\spvek{\mathbf{p}_0}{\mathbf{q}_0}\in 
		L^2(\Omega)\times \mathcal{T}^{-1}\grad H^1_{\Gamma_0}(\Omega)$, and denote $\spvek{\mathbf{p}(t)}{\mathbf{q}(t)}=\frakA(t)x_0$, $t\in[0,T]$. By \cite[Chap.~II, Lem.~1.3]{EngeNage00}, we have $\int_0^T\frakA(t)x_0{\rm d}t\in\dom(\tF H)=H^{-1}\dom(\tF)$ with
		\[\frakA(T)x_0-x_0=\tF H\int_0^T\frakA(t)x_0{\rm d}t.\]
		In particular,
		\[\mathbf{q}\changed{(T)}=\mathbf{q}_0+\mathcal{T}^{-1}\grad\int_0^T\rho^{-1}\mathbf{p}(t){\rm d}t.
		\]    
		Then the claim follows, since $\int_0^T\frakA(t)x_0{\rm d}t\in\dom(\tF H)$ implies $\rho^{-1}\mathbf{p}(t)\in H^1_{\Gamma_0}(\Omega)$. 
		
		
		\item[(b)]
		Let $u\in L^2([0,T];U)$ such that $\spvek{0}u\in\dom(\frakT_{I,T})$. The definition of the $F$-terminal value map gives rise to the existence of a~sequence $(u_n)$ in $H^2_{0l}([0,T];U)$, such that $(u_n)$ converges in 
		$L^2([0,T];U)$ to $u$, and $(\frakT_{I,T}\spvek{0}{u_n})=({\frakB}_Tu_n)$ converges in $X$ to \changed{$\frakB_{T}u$}. Denote 
		\[\spvek{\mathbf{p}_n(t)}{\mathbf{q}_n(t)}:={\frakB_T}(t)u_n,\quad n\in\N.\] 
		Then, by using that $u_n\in H^2_{0l}([0,T];U)$, we have 
		\[\forall\,t\in[0,T]:\quad \spvek{\spvek{\rho^{-1}\mathbf{p}_n(t)}{\mathcal{T}\mathbf{q}_n(t)}}{u_n(t)}\in\dom(\tF\&\tG)\]
		which particularly implies that $\rho^{-1}\mathbf{p}_n(t)\in H^1_{\Gamma_0}(\Omega)$ \changed{for all $t\in [0,T]$}.
		Further,
		\[{\frakB_T}u_n\spvek{\mathbf{p}_n(T)}{\mathbf{q}_n(T)}=\int_0^T\tF\&\tG\spvek{\spvek{\rho^{-1}\mathbf{p}_n(t)}{\mathcal{T}\mathbf{q}_n(t)}}{u_n(t)}{\rm d}t,\]
		which again leads to 
		\[
		\mathbf{q}_n\changed{(T)}=\mathcal{T}^{-1}\grad\int_0^T\rho^{-1}\mathbf{p}_n(t){\rm d}t\in \mathcal{T}^{-1}\grad H^1_{\Gamma_0}(\Omega).
		\]
		Now Proposition~\ref{prop:disprec} gives rise to $\langle \psi,\mathcal{T}\mathbf{q}_n(T)\rangle_{L^2(\Omega;\R^2)}=0$ for all
		$\psi\in H(\operatorname{div},\Omega)$ with
		$\operatorname{div}\psi=0$ and $\gamma_{N,\Gamma_1} \psi=0$. Now taking the limit $n\to\infty$, we obtain $\langle \psi,\mathcal{T}\mathbf{q}(T)\rangle_{L^2(\Omega;\R^2)}=0$, 
		and another use of Proposition~\ref{prop:disprec} yields that
		$\mathbf{q}(T)\in \mathcal{T}^{-1}\grad H^1_{\Gamma_0}(\Omega)$.
		\item[(c)] By using Proposition~\ref{prop:disprec_range}, we obtain that 
		\begin{align*}
		&\forall \spvek{\mathbf{v}}{\mathbf{w}}\in H^1_{\Gamma_0}(\Omega)\times L^2(\Omega):
		\\&\quad F_{\mathrm{disp,ext}}^*\spvek{\mathbf{v}}{\mathbf{w}}=\spvek{\mathbf{v}}{F_{\mathrm{disp}}^*\mathbf{w}}\in \setdef{\spvek{\mathbf{v}}{\mathbf{F}}\in {X}}{\,\mathbf{v} \in H^1_{\Gamma_0}(\Omega)\,,\mathbf{F}\in H(\operatorname{div},\Omega)\;\wedge\;\gamma_{N,\Gamma_1}\mathbf{F}=0}.
		\end{align*}
		In particular, the set 
		\[\setdef{\spvek{\mathbf{v}}{\mathbf{F}}\in L^2(\Omega)\times L^2(\Omega)}{\changed{F_{\mathrm{disp,ext}}^*\spvek{\mathbf{v}}{\mathbf{F}}} \in \dom(\tF)}\]
		is dense in $X$.
		\changed{We have $\dom(\tF)=\dom(\tF^*)=\dom((\tF H)^*)$ (now with respect to the standard inner product in $L^2$). This follows, since the operator in \eqref{eq:waveadj} is skew-adjoint, and $d$ represents a~bounded operator.}  Thus, we find ourselves precisely in the situation outlined in Remark~\ref{rem:denseass} \eqref{rem:denseassiib}. Leveraging the argumentation presented therein, we can verify that the $F_{\mathrm{disp,ext}}$-input-to-state map is well-defined.\end{itemize}   
\end{proof}

\noindent \textbf{The optimal control problem} 

We are considering an optimal control problem \changed{of the form \eqref{eq:energymin_terminalweight} where no constraints are imposed on the terminal state, that is, $F_c = 0$. Our objective is to achieve an energy-efficient control that (approximately) realizes a given displacement profile $\mathbf{w}_f\in L^2(\Omega)$ at time $T>0$ in a resting state, meaning that the displacement velocity at time $T>0$ is small.  The latter is represented by a terminal weight of the form $\frac12\| F_{\mathrm{disp,ext}}x-z_f\|_{L^2(\Omega)\times L^2(\Omega)}$, with $F_{\mathrm{disp,ext}}$ as in Proposition~\ref{prop:Fdisp}.} The corresponding optimal control problem is, under Assumptions~\ref{ass:wave}, and with $X$ as in \eqref{eq:Xwave}, $H$ as in \eqref{eq:Hwave}, $U$ as in \eqref{eq:Uwave}, $\tF\&\tG$, $\tK\&\tL$ as in \eqref{eq:wave_node},
\begin{align}\label{eq:waveopt}
\begin{split}
\textrm{minimize}&\quad
\int_0^T\langle {y(t)},{u(t)}\rangle_{U^*,U}{\rm d}t+\frac12\| F_{\mathrm{disp,ext}}x\changed{(T)}-\spvek{0}{\mathbf{w}_f}\|^2_{L^2(\Omega)\times L^2(\Omega)}\\
\text{subject to}&\quad\left(\begin{smallmatrix}\dot{x}(t)\\y(t)\end{smallmatrix}\right)=\left[\begin{smallmatrix}\phantom{-}\tF\&\tG\\[-1mm]\\-\tK\&\tL\end{smallmatrix}\right]\spvek{Hx(t)}{u(t)}, \quad
x(0) = x_0,\quad u \in \Uad.
\end{split}
\end{align}
In particular, the terminal weight is given by $\frac12\big(\|\mathbf{p}(T)\|_{L^2(\Omega)}^2+\|\mathbf{w}(T)-\mathbf{w}_f\|_{L^2(\Omega)}^2\big)$. 
The admissible controls are supposed to be evolving pointwisely in the set of all elements
$H^{-1/2}(\Gamma_1)$ which are represented by square integrable functions that evolve in a~box around zero. That is, for some $\underline{u},\overline{u}\in\R$ with $\underline{u}\leq\overline{u}$, we consider $\Uad \subset L^2([0,T];H^{-1/2}(\Gamma_1))$ with
\begin{equation}
\begin{aligned}
\Uad&=\setdef{u\in L^2([0,T];H^{-1/2}(\Gamma_1)}{u(t)\in U_{\rm ad}\text{ for a.e.\ }t\in[0,T]},\\
\text{where }U_{\rm ad}&=\setdef{u\in L^2(\Gamma)}{\underline{u}\leq u(\xi)\leq \overline{u}\text{ for a.e.\ }\xi\in \Gamma_1}\subset H^{-1/2}(\Gamma_1)
\end{aligned}\label{eq:Uadwave}
\end{equation}
In view of Remark~\ref{eq:modpHcont}, this leads to an optimal control problem 
\begin{align}\label{eq:waveopt_mod}
\begin{split}
\textrm{minimize}&\quad
\frac12\int_0^T\|{w(t)}\|_{L^2(\Omega)}^2\,{\rm d}t+\frac12\left\| \pvek{F_{\mathrm{disp,ext}}\changed{x(T)}-\spvek{0}{\mathbf{w}_f}}{\sqrt{H}\changed{x(T)}}\right\|^2_{L^2(\Omega)\times L^2(\Omega)\times X}\\
\text{subject to}&\quad\left(\begin{smallmatrix}\dot{x}(t)\\w(t)\end{smallmatrix}\right)=\left[\begin{smallmatrix}\tF\&\tG\\[-1mm]\\\tR\&\tS\end{smallmatrix}\right]\spvek{Hx(t)}{u(t)}, \quad
x(0) = x_0,\quad u \in \Uad,\quad .
\end{split}
\end{align}

The operator square root is simply given by the multiplication operator on $X={L^2(\Omega)}\times{L^2(\Omega;\R^2)}$
\[\sqrt{H}:X\to X, \quad\spvek{\mathbf{p}}{\mathbf{q}}\mapsto \spvek{\rho^{-1/2}\mathbf{p}}{\mathcal{T}^{1/2}\mathbf{q}}\]
For the initial value, we choose one corresponding to an initial velocity in $\grad H^1_{\Gamma_0}$, together with an initial force in $H(\operatorname{div},\Omega)$ with boundary trace in $U_{\rm ad}$ (as defined in \eqref{eq:Uadwave}), where we additionally assume that the initial force is initiated by the stress that corresponds to a~displacement. More precisely, we assume
\begin{align}
x_0=\spvek{\mathbf{p}_0}{\mathbf{q}_0}\in \rho H^1_{\Gamma_0}(\Omega)\times \mathcal{T}^{-1}H(\operatorname{div},\Omega)\text{ with }
\mathbf{q}_0\in\mathcal{T}^{-1}\grad H^1_{\Gamma_0}(\Omega)\text{ and }
\gamma_{N,\Gamma_1}\mathcal{T}\mathbf{q}_0\in U_{\rm ad}.\label{eq:waveinit}
\end{align}
To conclude existence of optimal controls and optimality conditions, we verify Assumptions~\ref{ass:ocp}. 
\begin{itemize}
	\item[\eqref{ass:ocp1}] As, by the argumentation \changed{after Assumption~\ref{ass:wave}}, $\sbvek{\tF\&\tG}{\tK\&\tL}\sbmat{H}00I$ is a~port-Hamiltonian system node, we can conclude from Proposition~\ref{prop:disseq} that \[\left[\begin{smallmatrix}\tF\&\tG\\[-1mm] \\\tR\&\tS
	\end{smallmatrix}\right]
	\sbmat{H}00{\Id_{U}}\]
	with $\tR\&\tS$ as in \eqref{eq:disswave} is a~system node.
	\item[\eqref{ass:ocp4}] This is shown in Proposition~\eqref{prop:Fdisp}\,\ref{prop:Fdispc}.
	\item[\eqref{ass:ocp2}] \changed{This is satisfied by assumption as $(0,\mathbf{w}_f,0,0)\in Z = L^2(\Omega)\times L^2(\Omega) \times L^2(\Omega) \times L^2(\Omega;\R^2)$ and as the initial value satisfies \eqref{eq:waveinit}.}
	\item[\eqref{ass:ocp3}] Convexity of $\Uad$ is clear. To show closedness of $\Uad$, in view of Remark~\ref{rem:optass} it suffices to show that ${U}_\mathrm{ad}$ as in \eqref{eq:Uadwave} is a~closed subset of $H^{-1/2}(\Gamma_1)$.
	Thus, consider a sequence $(u_n)\in {U}_\mathrm{ad}$ such that $u_n \to u\in H^{-1/2}(\Gamma_1)$. As $\Gamma_1$ has a~finite one-dimensional Lebesgue measure $|\Gamma_1|$, we have 
	\[\forall\,n\in\N:\quad\|u_n\|_{L^2(\Gamma_1)} \leq |\Gamma_1|^{\frac12}\max\{|\underline{u}|,|\overline{u}|\},\] the sequence $(u_n)$ is bounded in $L^2(\Gamma_1)$, it is, by considering suitable subsequences, no loss of generality to assume that $(u_n)$ is in $L^2(\Gamma_1)$ weakly convergent to some $\hat u \in L^2(\Gamma_1)$. 
	Due to convexity and as $U_{\rm ad}$ is closed in $L^2(\Gamma_1)$, we have $\hat u \in \Uad$ by the separation theorem. Further using the continuous embedding $L^2(\Gamma_1)\hookrightarrow H^{-1/2}(\Gamma_1)$, the limits coincide, i.e., $\hat u = u$. This shows that $u\in \changed{U_{\mathrm{ad}}}$.
	\item[\eqref{ass:ocp5}] As $\Gamma_1$ has finite one-dimensional Lebesgue measure, $U_{\rm ad}$ as in \eqref{eq:Uadwave} is a~bounded subset of $L^2(\Gamma_1)$, and thus also 
	a~bounded subset of $H^{-1/2}(\Gamma_1)$. Now, by using that $[0,T]$ is a~finite interval, we verify that $\Uad$ as in \eqref{eq:Uadwave} is a~bounded subset of $L^2([0,T];H^{-1/2}(\Gamma_1))$.
	\item[\eqref{ass:ocp6}] The condition \eqref{eq:waveinit} yields that there exists some $\hat{u}_0\in U_{\rm ad}$, such that 
	\[\spvek{x_0}{\hat{u}_0}\in\dom\left(\tF\&\tG\sbmat{H}00I\right).\] Now we may choose the constant input $\hat{u}\equiv \hat{u}_0\in \mathcal{U}_\mathrm{ad}$, and Proposition~\ref{prop:solex} gives rise to a classical trajectory $(x,u,y)$ in the sense of Definition~\ref{def:traj}. In particular, the terminal state is in $X$, and the output fulfills $y\in L^2([0,T];Y)$ for $Y=U^*$.
\end{itemize}



\begin{figure}[htb] 
	\centering
	\includegraphics[width=.52\linewidth]{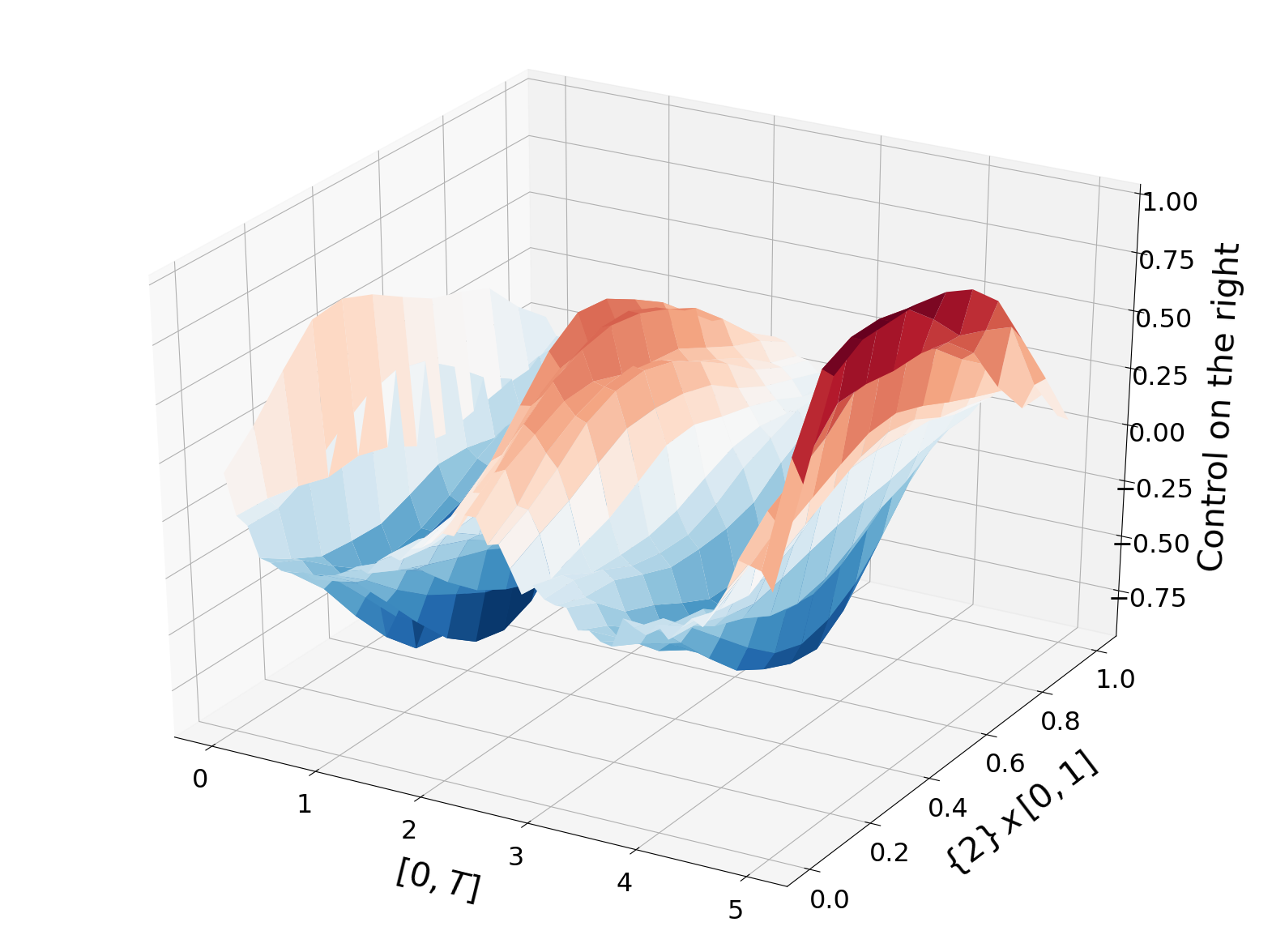}\hspace{.02\linewidth}
	\includegraphics[width=.45\linewidth]{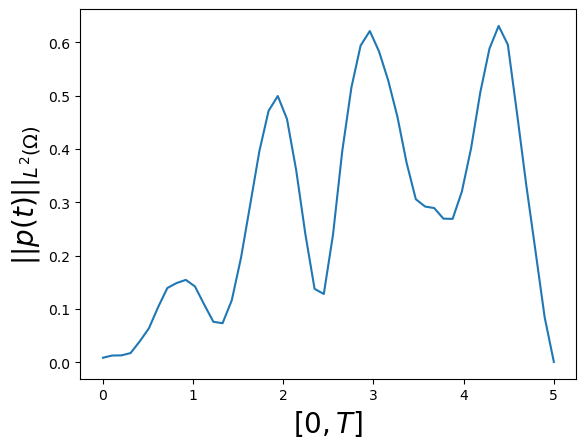}
	\caption{Control on the right boundary (left) and kinetic energy over time (right)}
	\label{fig:conmom}
\end{figure}

We now may conclude existence an optimal control using Theorem~\ref{thm:optcont}. 
\begin{cor}\label{cor:waveopt}
	Under Assumptions~\ref{ass:wave} on $\Omega\subset\R^2$, $\Gamma_0,\Gamma_1\subset\partial\Omega$, $\rho,\mathcal{T},d\in L^\infty(\Omega)$, and for $T>0$, $\mathbf{w}_f\in L^2(\Omega)$,
	and with $X$ as in \eqref{eq:Xwave}, $H$ as in \eqref{eq:Hwave}, $U$ as in \eqref{eq:Uwave}, $\tF\&\tG$, $\tK\&\tL$ as in \eqref{eq:wave_node}, the optimal control problem \changed{\eqref{eq:waveopt}} has a~solution. In other words, there exists an optimal control $u_{\rm opt}\in L^2([0,T];U)$ in the sense of Definition~\ref{def:cost}.
\end{cor}
We briefly comment on (possible) uniqueness of the optimal control.
\begin{rem}[Uniqueness of the optimal control problem \eqref{eq:waveopt}]
	\changed{Recall from Theorem~\ref{thm:optcont2} that uniqueness of an optimal control problem \eqref{eq:OCP} is, loosely speaking, equivalent to the property that the zero control is the only one which causes zero cost in conjunction with the trivial initial value. 
		For the problem \eqref{eq:waveopt}, this means, by invoking the reformulation \eqref{eq:waveopt_mod}, the following:
		If, for some $y\in L^2([0,T];U^*)$, 
		the wave equation \eqref{eq:waveeq} holds with $\mathbf{w}_0(\cdot)=0$, $\mathbf{v}_0(\cdot)=0$, $\mathbf{w}(\cdot,T)=0$, $\frac{\partial}{\partial t}\mathbf{w}(\cdot,T)=0$, and $d\frac{\partial}{\partial t}\mathbf{w}(\cdot,\cdot)\equiv0$, then $u:=\gamma_{N,\Gamma_1}\frac{\partial}{\partial t}\mathbf{w}$ has to vanish constantly.}
	This is, for instance, guaranteed, if the wave equation on $\Omega$ with input formed by the Dirichlet boundary of the velocity on $\Gamma_0$
	and ``distributed output'' formed by $d \tfrac{\partial}{\partial t}\mathbf{w}(\cdot,t)$ is observable. The latter issue has been addressed in previous work, such as in \cite{BLR92}, where conditions for observability are established. These conditions are derived from geometric considerations involving $\operatorname{supp}(d)$ (the support of $d$, representing the location of damping), $\Omega$ and $T$, illustrated through the concept of ``geodetic rays''.
\end{rem}
Again, we provide a numerical example using {FEniCS}~\cite{alnaes2015fenics} and dolfin-adjoint~\cite{mitusch2019dolfin}. To this end, we choose linear finite elements for the momentum variable and linear vector-valued finite elements for the stress. As a time-integrator, we choose the implicit midpoint rule. The parameters, including the desired target profile and the chosen finite element grid are depicted in Figure~\ref{fig:refwave}.\\
\begin{center}
	\begin{minipage}[]{.45\linewidth}
		\begin{align*}    
		x_0 &= 0, \ \mathcal{T}\equiv 1, \ \rho \equiv 1,\ d\equiv 0.05,\\
		T&=5,\ \mathbf{w}_f = \operatorname{dist}(\Gamma_0),\ \alpha_T = 10,\\
		\overline{u}&=1,\ \underline{u}=-1
		\end{align*} 
	\end{minipage}\hspace{.05\linewidth}
	\begin{minipage}[]{.45\linewidth}
		\vspace*{-.35cm}     
		\includegraphics[width=.9\columnwidth]{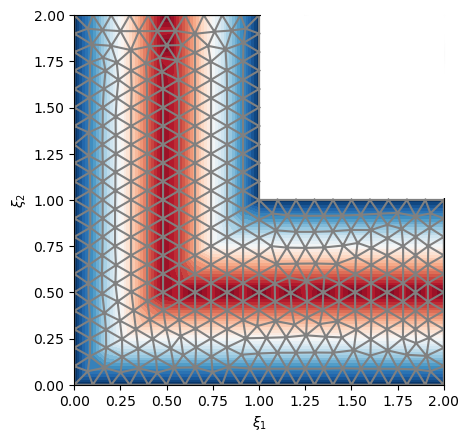}
	\end{minipage}
	\vspace{-.4cm}
	\captionof{figure}{\label{fig:refwave}Parameters (left) and target displacement (right).}
\end{center}~\\
In Figure~\ref{fig:conmom}, we depict the part of the optimal control on the right boundary $\{2\}\times [0,1]$ and the norm of the optimal momentum. In both figures, we observe a swing-up behavior of the optimal displacement and a small terminal momentum due to its penalization via $\|p(T)\|^2_{L^2(\Omega)}$ in the terminal weight. The swing-up behavior is necessary due to the control constraints limiting the force that can be applied at the boundary.

\begin{figure}[htb]
	\centering
	\includegraphics[width=.47\linewidth]{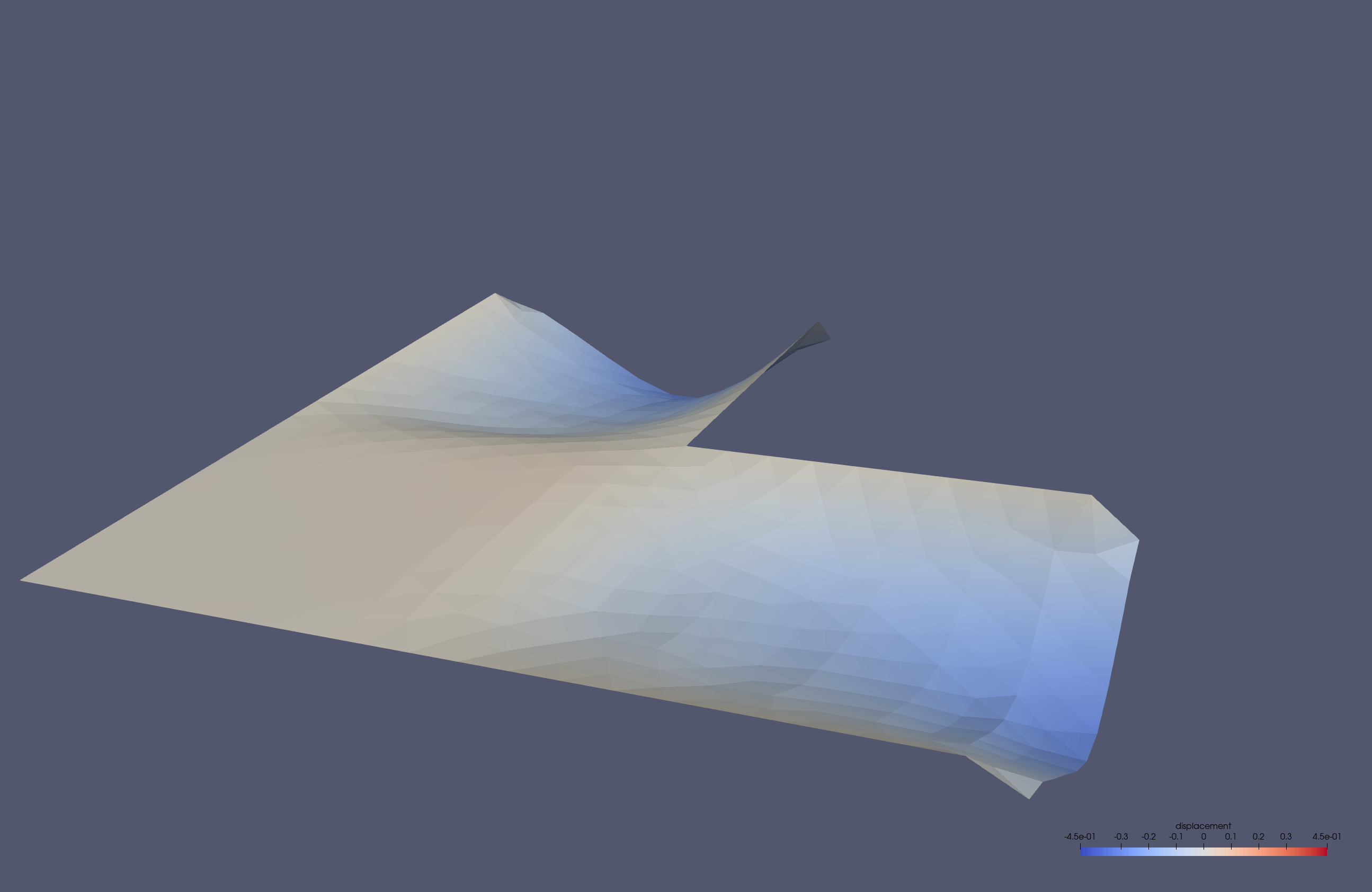}\hspace{.02\linewidth}
	\includegraphics[width=.47\linewidth]{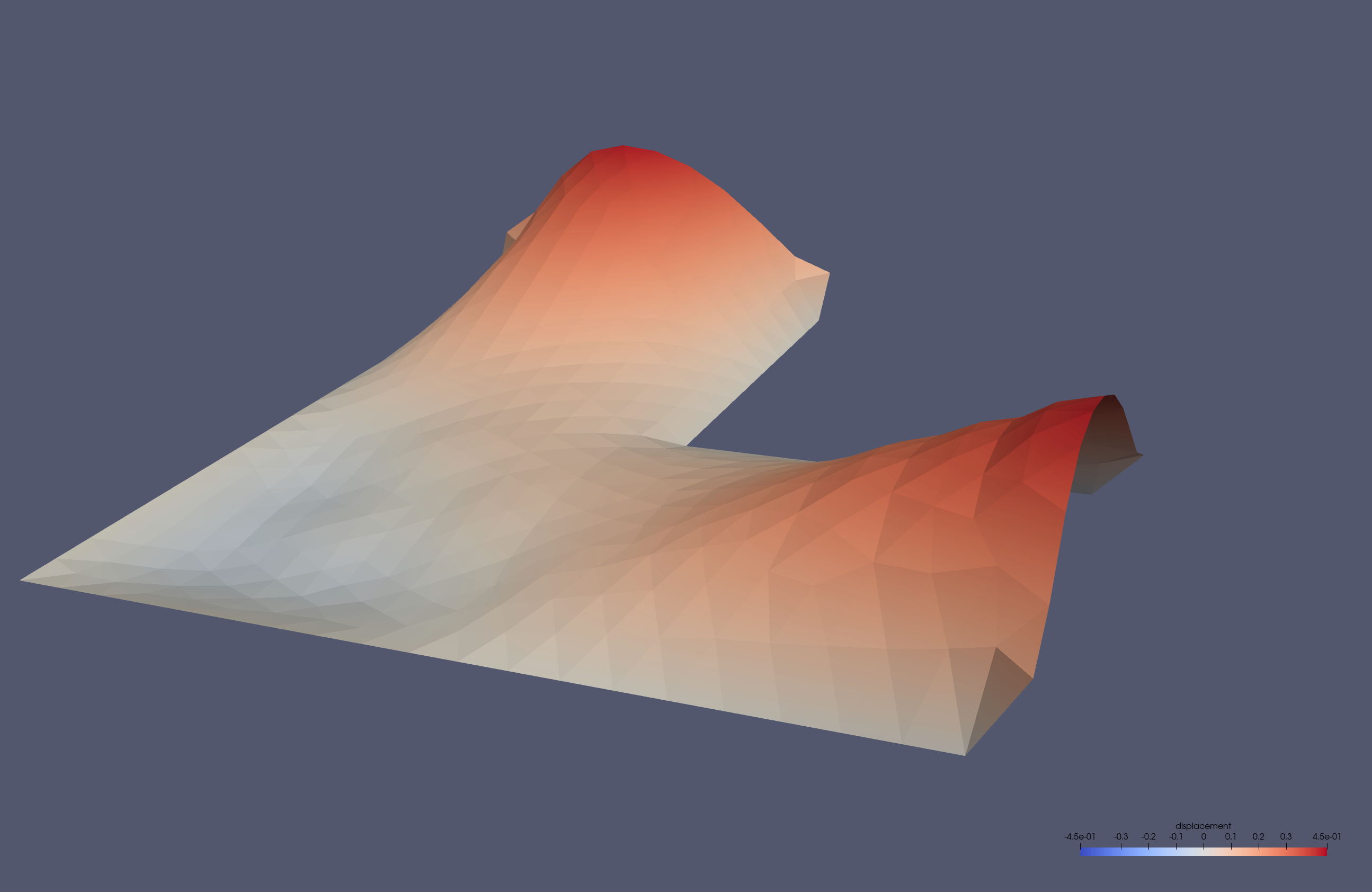}\\
	\includegraphics[width=.47\linewidth]{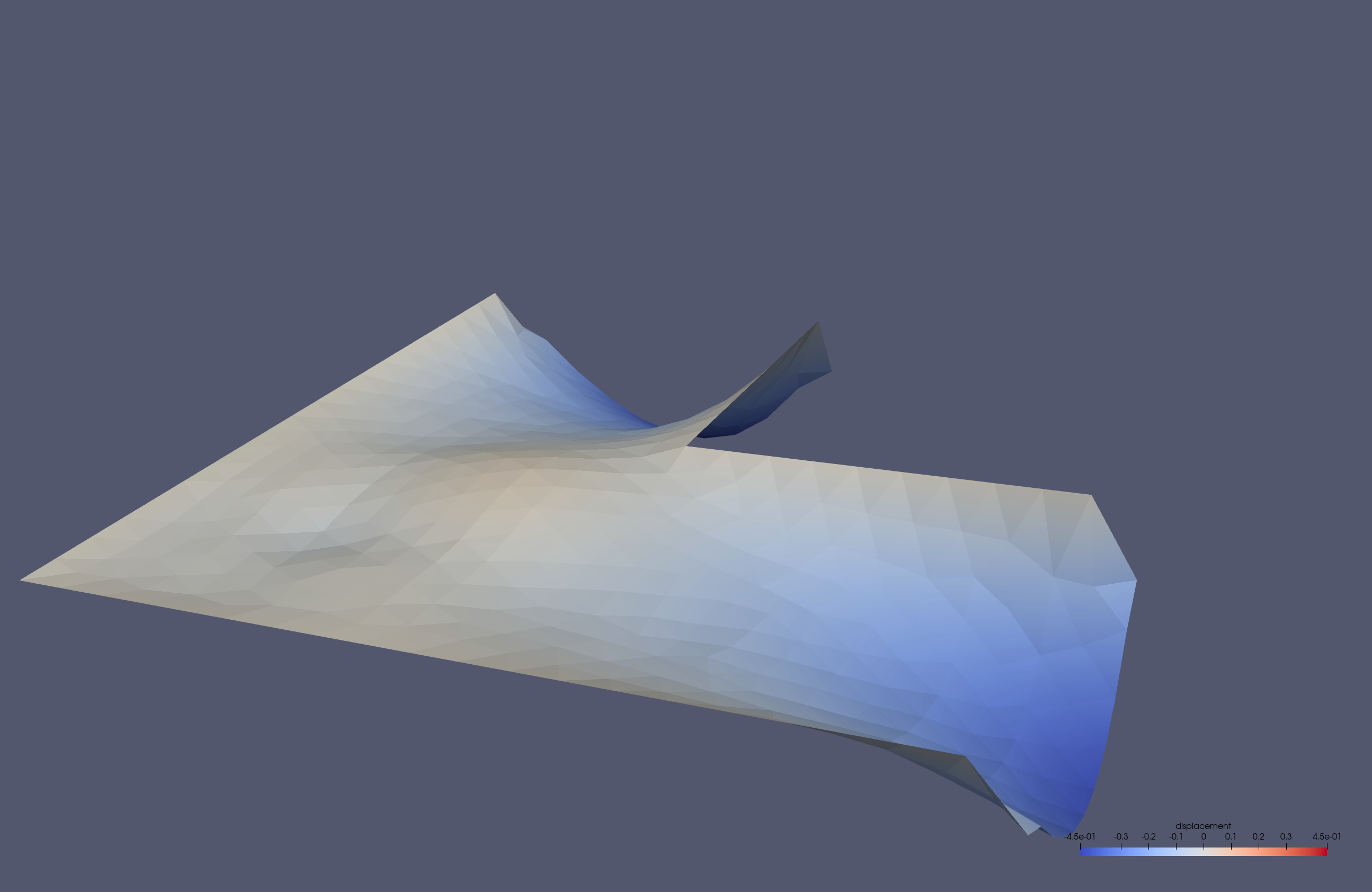}\hspace{.02\linewidth}
	\includegraphics[width=.47\linewidth]{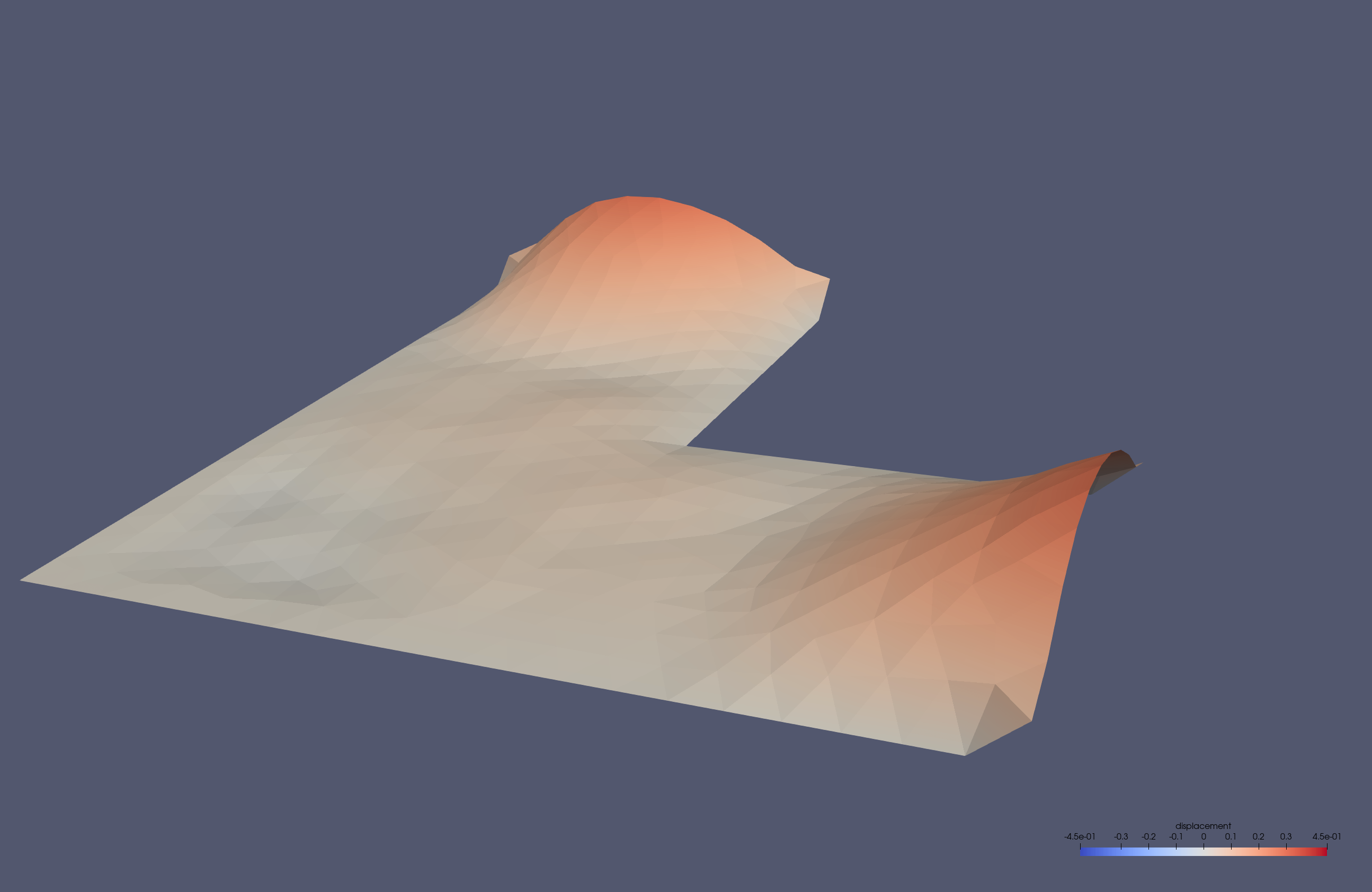}\\
	\includegraphics[width=.47\linewidth]{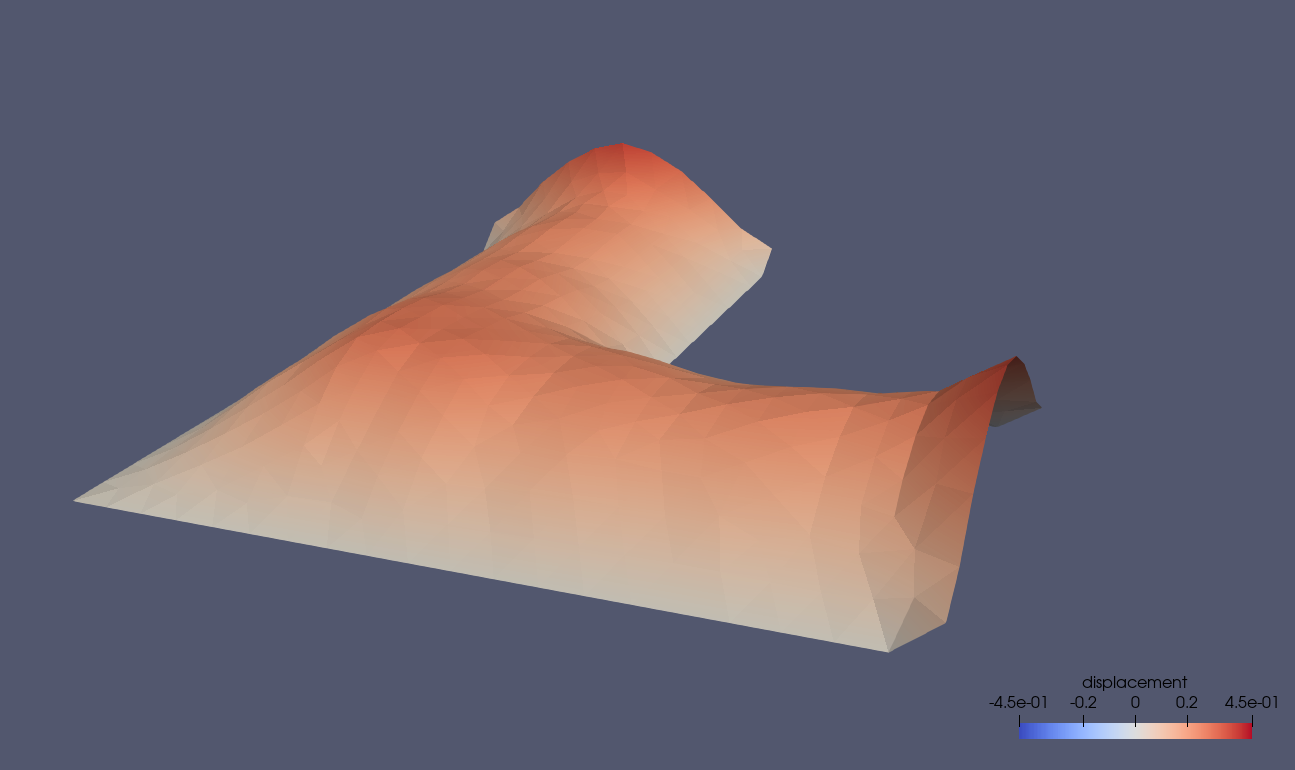}
	\caption{Snapshots of the displacement variable $\mathbf{w}(t)$ for times $t\in \{1.5,2.5,3.5,4.5,5\}$ (row-wise from top to bottom).}
	\label{fig:snapshots}
\end{figure}

\noindent Last, we provide snapshots of the displacement profile $\mathbf{w}(t)$ for different time instances $t\in [0,T]$ in Figure~\ref{fig:snapshots}\footnote{On the arXiv-page of this work, a video is available in the ancillary files section.}. As in the right plot of Figure~\ref{fig:conmom}, the swing-up behavior is clearly visible. Further, as can be seen in the in the last snapshot at the terminal time $t=T=5$, the terminal displacement $\mathbf{w}(T)$ is approximating the piecewise \changed{linear} reference signal $\mathbf{w}_f$ depicted on the right plot of Figure~\ref{fig:refwave}. This displacement is achieved with vanishing momentum variable, cf.\ the right plot of Figure~\ref{fig:conmom}.

\section{Conclusion}
In this work, we have proven existence of solutions and we have provided optimality conditions for linear-quadratic optimal control with abstract (and not necessarily well-posed) infinite-dimensional systems under control constraints and penalization of the terminal state and the $L^2$-norm of the output. Further, 
terminal state constraints have been discussed. We have introduced a~bunch of (not necessarily bounded) solution \changed{operators} using the theory of system nodes, which allowed us to show that under standard assumptions (such as convexity and closedness of the set of admissible inputs), optimal controls do exist. Further, we formulated an adjoint equation arising in the optimality condition by means of the adjoint system node. We further provided applications to port-Hamiltonian system nodes for which we have aimed for energy-optimal state transitions. Last, we have provided two numerical examples by means of Dirichlet boundary controlled heat equation in one spatial variable and an energy efficient control of a~port-Hamiltonian system formed by a~two-dimensional boundary-controlled wave equation on an L-shaped domain.


\bibliographystyle{abbrv}
\bibliography{references.bib}

\begin{landscape}
	\appendix
	\section{Summary of involved operators}\label{sec:ourops}
	\vspace{1cm}
	
	\begin{tabular}[h]{|c|c|c|c|c|c|c|}
		\hline
		
		\makecell[c]{\\\bf Symbol\\[3mm] } & {\bf Name} & {\bf Domain} & {\bf Target Space} & {\bf Argument} & {\bf Value} & {\bf Where defined?}\\
		\hhline{|=|=|=|=|=|=|=|}
		\makecell[c]{\\$\frakA_{-1}(T)$\\[3mm] }
		& \makecell[c]{semigroup \\on $X_{-1}$ at $T$} & $X_{-1}$ & $X_{-1}$ & $x_0$ & 
		\makecell[c]{$x(T)$ of \eqref{eq:ODEnode}\\ with $u=0$}  & Rem.\ \ref{rem:nodes} \\
		\hline
		\makecell[c]{\\${\frakB_T}$\\[3mm] } & input-to-state map & $L^2([0,T];U)$ & $X_{-1}$ & $u$ & 
		\makecell[c]{$x(T)$ of \eqref{eq:ODEnode}\\ with $x_0=0$}  & eq.\ \eqref{eq:genISmap} \\
		\hline
		\makecell[c]{\\${\frakC_T}$\\[3mm] }
		& state-to-output map & $X_{-1}$ & $H^{-2}_{0l}
		([0,T];Y)$ & $x_0$ & 
		\makecell[c]{$y$ of \eqref{eq:ODEnode}\\ with $u=0$}  & eq.\ \eqref{eq:genSOmap} \\
		\hline
		\makecell[c]{\\${\frakD_T}$\\[3mm] }
		& input-to-output map & $L^2([0,T];U)$ & $H^{-2}_{0l}([0,T];Y)$ & $u$ & 
		\makecell[c]{$y$ of \eqref{eq:ODEnode}\\ with $x_0=0$}  & eq.\ \eqref{eq:genIOmap} \\
		\hline
		\makecell[c]{\\$\frakT_{F,T}$\\[3mm] }
		& $F$-terminal value map & \makecell[c]{$\dom({\frakT_{F,T}})\subset$\\ $X_{-1}\times L^2([0,T];U)$} & $Z$ & $(x_0,u)$ & 
		{$Fx(T)$ of \eqref{eq:ODEnode}}  & Def.~\ref{def:FIS} \\
		\hline
		\makecell[c]{\\$\frakI_{F,T}$\\[3mm] }
		& $F$-input map & \makecell[c]{$\dom({\frakI_{F,T}})\subset$\\ $L^2([0,T];U)$} & $Z\times L^2([0,T];Y)$ & $u$ & 
		\makecell[c]{$(Fx(T),y)$ of \eqref{eq:ODEnode}\\ with $x_0=0$} & Def.~\ref{def:FIO} \\
		\hline
		\makecell[c]{\\$\frakO_{G,T}$\\[3mm] }
		& $G$-output map & \makecell[c]{$\dom({\frakO_{G,T}})\subset$\\ $Z\times L^2([0,T];U)$} & $ L^2([0,T];Y)$ & $(z,u)$ & 
		\makecell[c]{$y$ of \eqref{eq:ODEnode}\\ with $x_0=Gz$} & Def.~\ref{def:FIO} \\
		\hline\end{tabular}
	\captionof{table}{Involved operators}\label{tab:ourops}
\end{landscape}
\end{document}